\documentclass[11pt,twoside]{article}

\usepackage[maxalphanames=1,
  minalphanames=1,
  style=alphabetic,
  backend=bibtex,
  firstinits=true,
  url=true]{biblatex}

\AtBeginBibliography{\small}
\usepackage{fullpage}

\usepackage{epsf} \usepackage{fancyhdr} \usepackage{graphics}
\usepackage{graphicx} \usepackage{psfrag} \usepackage{microtype}

\usepackage{subfigure}
\usepackage{algorithmic} \usepackage{color,xcolor}
\usepackage{arydshln}

\usepackage[linesnumbered,ruled]{algorithm2e}%

\DontPrintSemicolon \usepackage{color}

\usepackage{amsthm} \usepackage{amsfonts} \usepackage{amsmath}
\usepackage{amssymb,bbm} \usepackage{comment} \usepackage{booktabs}

\usepackage{url}
\usepackage[colorlinks, linkcolor=magenta, citecolor=blue]{hyperref}

\usepackage[noabbrev]{cleveref}

\usepackage{nicefrac} \usepackage{comment}

\usepackage{chngpage}

 \usepackage{tabularx}%

\usepackage{enumitem}
\usepackage{booktabs}
\usepackage{caption}

\usepackage{bm,bbm}
\usepackage{mathtools}

\newcommand*\R[0]{\mathbb{R}}

\newcommand*\tr[0]{\text{tr}}

\newcommand{\real}{\ensuremath{\mathbb{R}}}

\newcommand{\eps}{\epsilon}

\newcommand{\order}[1]{\ensuremath{\mathcal{O}\parenth{#1}}}

\newcommand{\thetastar}{\ensuremath{\theta^*}}
\newcommand{\thetahat}{\ensuremath{\widehat{\theta}}}

\newcommand{\parenth}[1]{\left( #1 \right)}

\newcommand{\abss}[1]{\left| #1 \right |}

\newcommand{\vterm}{V}
\newcommand{\vtermtil}{\widetilde{V}}

\newcommand{\sphere}{\ensuremath{\mathbb{S}}}
\newcommand{\cN}{\ensuremath{\mathcal{N}}}

\newcommand{\mydefn}{\ensuremath{:=}}

\newcommand{\defn}{:=}

\newcommand{\matsnorm}[2]{|\!|\!| #1 | \! | \!|_{{#2}}}

\newcommand{\vecnorm}[2]{\| #1\|_{#2}}
\newcommand{\enorm}[1]{\vecnorm{#1}{2}} 
\newcommand{\rnorm}[2]{\| #1\|_{#2}}

\newcommand{\fronorm}[1]{\ensuremath{\matsnorm{#1}{\footnotesize{\mbox{fro}}}}}
\newcommand{\opnorm}[1]{\ensuremath{\matsnorm{#1}{\tiny{\mbox{op}}}}}

\newcommand{\inprod}[2]{\ensuremath{\langle #1 , \, #2 \rangle}}
\newcommand{\binprod}[2]{\ensuremath{\Big\langle #1 , \, #2 \Big\rangle}}

\newcommand{\vtinprod}[2]{\ensuremath{\langle #1 , \, #2 \rangle}}

\newcommand{\Exs}{\ensuremath{{\mathbb{E}}}}
\newcommand{\Prob}{\ensuremath{{\mathbb{P}}}}

\DeclareMathOperator{\var}{var}
\DeclareMathOperator{\cov}{cov}
\DeclareMathOperator{\trace}{trace}

\newtheoremstyle{named}{}{}{\itshape}{}{\bfseries}{.}{.5em}{\thmnote{#3's }#1}
\theoremstyle{named}

\theoremstyle{plain}

\newlength{\widebarargwidth}
\newlength{\widebarargheight}
\newlength{\widebarargdepth}
\DeclareRobustCommand{\widebar}[1]{%
  \settowidth{\widebarargwidth}{\ensuremath{#1}}%
  \settoheight{\widebarargheight}{\ensuremath{#1}}%
  \settodepth{\widebarargdepth}{\ensuremath{#1}}%
  \addtolength{\widebarargwidth}{-0.3\widebarargheight}%
  \addtolength{\widebarargwidth}{-0.3\widebarargdepth}%
  \makebox[0pt][l]{\hspace{0.3\widebarargheight}%
    \hspace{0.3\widebarargdepth}%
    \addtolength{\widebarargheight}{0.3ex}%
    \rule[\widebarargheight]{0.95\widebarargwidth}{0.1ex}}%
  {#1}}

\makeatletter
\long\def\@makecaption#1#2{
        \vskip 0.8ex
        \setbox\@tempboxa\hbox{\small {\bf #1:} #2}
        \parindent 1.5em  
        \dimen0=\hsize
        \advance\dimen0 by -3em
        \ifdim \wd\@tempboxa >\dimen0
                \hbox to \hsize{
                        \parindent 0em
                        \hfil
                        \parbox{\dimen0}{\def\baselinestretch{0.96}\small
                                {\bf #1.} #2
                                }
                        \hfil}
        \else \hbox to \hsize{\hfil \box\@tempboxa \hfil}
        \fi
        }
\makeatother

\long\def\comment#1{}
\definecolor{battleshipgrey}{rgb}{0.52, 0.52, 0.51}
\definecolor{darkgray}{rgb}{0.66, 0.66, 0.66}
\definecolor{darkgreen}{rgb}{0.0, 0.2, 0.13}
\definecolor{darkspringgreen}{rgb}{0.09, 0.45, 0.27}
\definecolor{dukeblue}{rgb}{0.0, 0.0, 0.61}
\definecolor{olivedrab7}{rgb}{0.24, 0.2, 0.12}
\definecolor{darkblue}{rgb}{0.0, 0.0, 0.55}
\definecolor{darkscarlet}{rgb}{0.34, 0.01, 0.1}
\definecolor{candyapplered}{rgb}{1.0, 0.03, 0.0}
\definecolor{ao(english)}{rgb}{0.0, 0.5, 0.0}
\definecolor{applegreen}{rgb}{0.55, 0.71, 0.0}

\newcommand{\widgraph}[2]{\includegraphics[keepaspectratio,width=#1]{#2}}

\newcommand{\E}{\mathbb E}

\def\ParDim{d}

\def\Numobs{n}
\def\Sam{Z}

\def\ZFun{{h}}

\newcommand{\mymatrix}[1]{\ensuremath{\mathbf{#1}}}

\def\GradFun{{\mymatrix{J}}}
\def\TargetFun{{\tau}}
\def\Par{{\theta}}
\def\TruePar{{\theta^*}}
\def\EstPar{{\widehat\theta}}

\def\ParSpace{{\Theta}}

\def\IntSec{{\eta}}
\def\IntPar{{\widetilde\theta}}

\def\ConParNorm{{\MyBou}}

\newcommand{\OrConb}{}
\def\OrConc{{2/3}}
\def\OrCond{{2/3}}
\def\InvOrConc{{3/2}}
\def\InvOrCond{{3/2}}
\def\OrPar{{\sigma}}

\def\OrPara{{\OrPar}}
\def\OrParb{{\OrPar}}
\def\OrParc{{\OrPar}}
\def\OrPard{{\OrPar}}

\newcommand{\MyBou}{\ensuremath{\rho}}

\def\LipZFa{\sigma}
\def\LipZFb{\sigma}

\newcommand{\PlainLip}{\ensuremath{L}}
\newcommand{\LipTara}{\ensuremath{\PlainLip}}
\newcommand{\LipTarb}{\ensuremath{\PlainLip}}
\newcommand{\LipTarc}{\ensuremath{\PlainLip}}

\def\Sma{{S}}

\newcommand{\LinApp}{\ensuremath{L}}

\def\BGradTar{{L}}
\def\BGradTartwo{{L}}

\newcommand{\Errfirst}{\ensuremath{\LinError_\numobs}}

\def\SupEm{Q}
\def\Direc{u}
\def\Direcb{v}
\def\Direcc{w}
\def\Direcd{{s}}

\def\Coverset{{\mathcal{M}}}
\def\Funclass{{\mathcal{F}}}
\def\Ccerr{{\dimdelta}}
\def\Factora{{b}}
\def\Factorb{{f}}

\def\GZInv{{\boldsymbol M}}
\def\GZ{{\tilde{\boldsymbol  M}}}

\newcommand{\ZF}{v}

\def\JackEst{{\widehat\TargetFun_{ \mathrm{jac}}}}

\def\Term{{T}}
\def\STerm{{\tilde {T}}}
\def\Res{{R}}

\newcommand{\specialpar}{\ensuremath{\nu}}
\def\Longc{{\specialpar}}

\def\Jack{{\mathrm{jac}}}

\def\JackVar{{\widehat{V}_{\Jack,\Numobs}}}

\def\EmpGZF{{\GradFun_{\Numobs}}}  

\newcommand{\subsup}[3]{ #1_{#2, -#3}}

\def\EmpGZFa{{ \subsup{\widetilde{\GradFun}}{\Numobs}{i}}}  
\def\EmpGZFb{{\widebar{\GradFun}_{\Numobs,-i}}}   

\def\EmpGZFc{{\widehat{\GradFun}_{\Numobs-1,-i}}}   

\def\EmpGZFcall{{\widehat{\GradFun}_{\Numobs,-i}}}

\newcommand{\OneDirScore}{\ensuremath{\ell}}

\usepackage{xcolor}

\newcommand{\numobs}{n}
\newcommand{\usedim}{d}

\newcommand{\sionebase}{\phi}
\newcommand{\sitwobase}{\psi}

\newcommand{\sione}{\sionebase_i}
\newcommand{\sitwo}{\sitwobase_i}

\newcommand{\AVEHAT}[1]{\widehat{#1}_\numobs}
\newcommand{\sioneave}{\AVEHAT{\sionebase}}
\newcommand{\sitwoave}{\AVEHAT{\sitwobase}}

\newcommand{\GradFunminusi}{{\widetilde\GradFun_{\Numobs,-i}}}

\newcommand{\taustar}{\tau^*}
\newcommand{\truepropscore}{\pi}

\newcommand{\betastar}{\beta^*}

\newcommand{\ivcoef}{\pi}
\newcommand{\iverror}{\eta}
\newcommand{\ivdim}{k}
\newcommand{\numtrial}{T}
\newcommand{\tind}{t}
\newcommand{\unbiased}{\widehat\TargetFun_{\text{unbiased}}}

\newcommand{\sitwowithout}{\psi}

\newcommand{\RegLoss}{F}
\newcommand{\Zfunexam}{f}

\newcommand{\ebias}{\widehat{\mathrm{Bias}}}
\newcommand{\emse}{\widehat{\mathrm{MSE}}}

\newcommand{\stdest}{\widehat{\sigma}}
\newcommand{\unbiasedname}{{\mathrm{unbiased}}}

\newcommand{\betaspace}{{\mathcal{B}}}
\newcommand{\smoothtrun}{{s}}

\newcommand{\newtruncate}{{T}}

\newcommand{\LooBound}{B}

\newcommand{\strongconvex}{\ensuremath{\gamma}}

\newcommand{\myassumption}[3]{
  \begin{enumerate}[label={\bf{(#1)}}]
  \item \label{#2} {#3}
  \end{enumerate}
 }

\newcommand{\coordinate}{e}

\newcommand{\indicator}{\mathbb{I}}

\newcommand{\tauhat}{\MYHAT{\avgtreat}}

\theoremstyle{definition} \newtheorem{example}{Example}
\newcommand{\goodendex}{\hfill $\clubsuit$}

\newtheorem{lemma}{Lemma}
\newtheorem{proposition}{Proposition}
\newtheorem{corollary}{Corollary} \newtheorem{assumption}{Assumption}

\usepackage{mathrsfs}

\usepackage{amsmath,amsfonts,amsthm,amssymb} \usepackage{graphicx,
  color, setspace}

\usepackage{comment}

\usepackage{float} \usepackage{indentfirst}

 \newtheorem{theorem}{Theorem}

\newcommand{\avgtreat}{\tau}

\newcommand{\Event}{\mathscr{E}}

\long\def\comment#1{}  \long\def\comment#1{}

\newcommand{\funcClass}{\mathcal{F}}

\newcommand{\Action}{A} \newcommand{\action}{a}
\newcommand{\Outcome}{Y} 
\newcommand{\State}{X}

\newcommand{\IdMat}{\ensuremath{\mathbf{I}}}

\newcommand{\lammax}{\ensuremath{\lambda_{\mbox{\tiny{max}}}}}

\renewcommand{\usedim}{{\ParDim}}

\newcommand{\noiseW}{W} \newcommand{\highbias}{B_\usedim}
\newcommand{\highorder}{\ensuremath{H}}

\newcommand{ \logistic}[2]{\pi (#1; #2)}

\usepackage{scalerel}

\newcommand{\tauhatipw}{\smallsuper{\tauhat}{IPW}}

\newcommand{\covar}{X} 

\newcommand{\act}{A}

\usepackage{mathrsfs}

\newcommand{\Sphere}[1]{\ensuremath{\mathbb{S}^{#1}}}

\newcommand{\Remainder}{\ensuremath{R}}

\newcommand{\Fclass}{\ensuremath{\mathcal{F}}}

\newcommand{\Bias}{\ensuremath{B}}
\newcommand{\estedpropscore}{\widehat{\pi}}

\renewcommand{\thetastar}{\ensuremath{{\theta^{*}}}}

\newcommand{\Amat}{\ensuremath{\mymatrix{A}}}
\newcommand{\Bmat}{\ensuremath{\mymatrix{B}}}

\newcommand{\Jmat}{\ensuremath{\mymatrix{J}}}
\newcommand{\Mmat}{\ensuremath{\mymatrix{M}}}

\newcommand{\Normal}{\ensuremath{\mathcal{N}}}

\newcommand{\newexponent}{\ensuremath{\delta}}

\newcommand{\covariate}{x}

\newcommand{\sumn}{\sum_{i=1}^{\numobs}}

\def\Hesdif{{\Mmat}} 

\newcommand{\Mhat}{\ensuremath{\widehat{\Mmat}_\numobs}}

\newcommand{\polyshort}{\ensuremath{C}}
\newcommand{\poly}{\polyshort(L, \sigma, \MyBou, \gamma)}

\newcommand{\polyshortprime}{\ensuremath{\polyshort'}}
\newcommand{\polyprime}{\polyshortprime(L,\sigma, \MyBou, \gamma)}

\newcommand{\MYHAT}[1]{\ensuremath{\widehat{#1}}}
\newcommand{\betahat}{\MYHAT{\beta}}

\newcommand{\ivcoeffstar}{\ivcoef^*}
\newcommand{\alphastar}{\ensuremath{\alpha^*}}

\newcommand{\DelJack}{\ensuremath{\smallsuper{\varepsilon_\numobs}{jack}}}
\newcommand{\DelPlug}{\ensuremath{\smallsuper{\varepsilon_\numobs}{plug}}}

\newcommand{\sigtil}{\ensuremath{\widetilde{\sigma}}}

\newcommand{\betatil}{\ensuremath{\widetilde{\beta}}}
\newcommand{\tautil}{\ensuremath{\widetilde{\tau}}}

\usepackage{scalerel}
\newcommand{\smallsuper}[2]{#1^{\scaleto{\mbox{#2}}{4pt}}}
\newcommand{\smallsub}[2]{#1_{\scaleto{\mbox{#2}}{4pt}}}

\newcommand{\MyIntSec}[1]{\ensuremath{\IntSec_{#1}}}

\newcommand{\Err}{\ensuremath{\widehat{\Delta}}}

\newcommand{\myusedim}{\ensuremath{\Ccerr}}

\newcommand{\LOO}[3]{\ensuremath{ {#1}_{#2}^{(-#3)}}}

\newcommand{\StochNoise}{\ensuremath{W}}
\newcommand{\StochNoisebase}{\tilde{W}}

\newcommand{\PopZfun}{\ensuremath{H}}

\newcommand{\SuffStat}{\ensuremath{T}}

\newcommand{\revdefn}{\ensuremath{=\!\colon}}

\newcommand{\simpleorpar}{\ensuremath{\sigma}}

\newcommand{\instrument}{\ensuremath{W}}

\newcommand{\Zfun}{\ZFun}

\newcommand{\pistar}{\pi^*}

\newcommand{\Pardim}{\ParDim}

\newcommand{\Pdist}{\ensuremath{\mathbb{P}}}

\newcommand{\Debias}{\ensuremath{D}}

\newcommand{\TermTil}{\ensuremath{\widetilde{T}}}

\newcommand{\LinError}{\ensuremath{\LinErrorPl}}
\newcommand{\LinErrorPl}{\ensuremath{E}}
\newcommand{\RemainZ}{\ensuremath{R_\ZFun}}
\newcommand{\RemainTau}{\ensuremath{R_\TargetFun}}
\newcommand{\RemainL}{\ensuremath{R_{\LinErrorPl}}}

\newcommand{\RemainLone}{\ensuremath{R_{\LinErrorPl \, 1}}}
\newcommand{\RemainLtwo}{\ensuremath{R_{\LinErrorPl \, 2}}}
\newcommand{\RemainLthree}{\ensuremath{R_{\LinErrorPl \, 3}}}

\newcommand{\Ustat}{\ensuremath{U_\numobs}}

\newcommand{\PlugDiff}{\ensuremath{D_{\mbox{\tiny{plug}}}}}

\newcommand{\dimdelta}{\ensuremath{\ParDim_\delta}}

\renewcommand{\Jmat}{\ensuremath{\GradFun_\TruePar}}

\newcommand{\UstatTil}{\ensuremath{\tilde{U}_\numobs}}
\newcommand{\DelHat}{\ensuremath{\widehat{\Delta}_\numobs}}

\newcommand{\dimdeltatil}{\ensuremath{\ParDim_{\numobs, \delta}}}

\newcommand{\Direca}{\ensuremath{\Direc}}

\newcommand{\Parspace}{\ensuremath{\ParSpace}}

\newcommand{\myunder}[1]{\noindent{\underline{#1}}}

\newcommand{\NewCover}{\mathcal{N}}

\newcommand{\Ytil}{\ensuremath{\widetilde{Y}}}

\newcommand{\ParEst}{\EstPar}

\newcommand{\PlainDelHat}{\ensuremath{\widehat{\Delta}}}

\newcommand{\LooDiff}{\ensuremath{\LOO{\PlainDelHat}{\numobs-1}{i}}}
\newcommand{\bigenorm}[1]{\ensuremath{ \Big \| #1 \Big \|_2}}

\newcommand{\VarTrue}{\ensuremath{\TrueVarPlain}_\numobs}
\newcommand{\TrueVar}{\VarTrue}

\newcommand{\TrueVarPlain}{\nu}
\newcommand{\VarTruePlain}{\TrueVarPlain}

\newcommand{\newfigdir}{figs/fig_new}

\newcommand{\highdimto}{\ensuremath{\rightsquigarrow}}

\newcommand{\convprob}{\ensuremath{\stackrel{\tiny{\mbox{prob.}}}{\longrightarrow}}}
\newcommand{\convdist}{\ensuremath{\stackrel{\tiny{\mbox{dist.}}}{\longrightarrow}}}

\newcommand{\QmatZ}{\ensuremath{{\mathbf{Q}^\ZFun}}}
\newcommand{\QmatTau}{\ensuremath{{\mathbf{Q}^\TargetFun}}}

\newcommand{\Qmat}{\ensuremath{\mathbf{Q}}}

\newcommand{\rexp}{\ensuremath{r}}
\newcommand{\rexplow}{\ensuremath{\underbar{\rexp}}}
\newcommand{\rexpup}{\ensuremath{\widebar{\rexp}}}

\newcommand{\highplug}{\ensuremath{\highorder^{\mbox{\tiny{plug}}}_\numobs}}
\newcommand{\highjack}{\ensuremath{\highorder^{\mbox{\tiny{jack}}}_\numobs}}

\newcommand{\noise}{\ensuremath{\varepsilon}}

\renewcommand{\Bmat}{\ensuremath{\mathbf{B}}}

\newcommand{\Smat}{\ensuremath{\mathbf{S}}}

\newcommand{\mytrans}{\ensuremath{\top}}

\begin{document}

\begin{center}
  {\bf{\LARGE{When is it worthwhile to jackknife?  \\ Breaking the
        quadratic barrier for $Z$-estimators}}}

\vspace*{.2in} {\large{
 \begin{tabular}{ccccc}
  Licong Lin $ ^{\dagger, \star}$ & Fangzhou Su $^{ \dagger, \star}$ &
  Wenlong Mou$^{ \diamond}$ & Peng Ding$^{\dagger}$ & Martin
  J. Wainwright$^{\dagger, \ddagger}$
 \end{tabular}
}

 \begin{tabular}{c}
 Department of Statistics$^\dagger$, UC Berkeley\\
 \end{tabular}

 \medskip
 
 \begin{tabular}{c}
   Departments of EECS \& Mathematics$^\ddagger$ \\ Lab for
   Information and Decision Systems \\ Statistics and Data Science
   Center \\ Massachusetts Institute of Technology
 \end{tabular}

 \begin{center}
   Department of Statistics, University of Toronto$^\diamond$
 \end{center}
}

\begin{abstract}
  Resampling methods are especially well-suited to inference with
  estimators that provide only ``black-box'' access.  Jackknife is a
  form of resampling, widely used for bias correction and variance
  estimation, that is well-understood under classical scaling where
  the sample size $\numobs$ grows for a fixed problem.  We study its
  behavior in application to estimating functionals using
  high-dimensional $Z$-estimators, allowing both the sample size
  $\Numobs$ and problem dimension $\ParDim$ to diverge.  We begin
  showing that the plug-in estimator based on the $Z$-estimate suffers
  from a quadratic breakdown: while it is $\sqrt{\numobs}$-consistent
  and asymptotically normal whenever $\numobs \gtrsim \ParDim^2$, it
  fails for a broad class of problems whenever $\numobs \lesssim
  \ParDim^2$.  We then show that under suitable regularity conditions,
  applying a jackknife correction yields an estimate that is
  $\sqrt{\numobs}$-consistent and asymptotically normal whenever
  $\Numobs\gtrsim \ParDim^{3/2}$. This provides strong motivation for
  the use of jackknife in high-dimensional problems where the
  dimension is moderate relative to sample size.  We illustrate
  consequences of our general theory for various specific
  $Z$-estimators, including non-linear functionals in linear models;
  generalized linear models; and the inverse propensity score weighting (IPW) estimate for the average
  treatment effect, among others.
  \let\thefootnote\relax\footnote{$^\star$LL and FS contributed
  equally to this work.}
\end{abstract}
\end{center}


\section{Introduction}
\label{SecIntro}

It is common in modern statistics to encounter estimators that are
defined implicitly by a computational routine: while it is possible to
compute the estimate for any dataset used as input, details of the
estimator (and the underlying statistical model) are unknown.  Such
access is known as ``black box'', since the estimator is available
only as an input-output device.  Resampling methods such as the
jackknife and bootstrap are especially well-suited to inference with
black-box estimators, since they require only successive computations
of the estimate on suitably perturbed datasets.

With this broader perspective in mind, the main contribution of this
paper is develop insight into the behavior of the jackknife for
black-box inference, with particular emphasis on its finite sample
properties and dependence on the problem dimension.  We study this
general problem in the context of estimating functionals using a
$Z$-estimator.  The class of $Z$-estimators is very broad, with
applications across a wide range of empirical disciplines and a rich
theory (e.g., see the books~\cite{bickel1993efficient,vaart1996weak,
  van2000asymptotic,
  kosorok2008introduction,tsiatis2006semiparametric} for further
background).  At the population level, any parametric $Z$-estimator is
based upon representing a parameter of interest $\TruePar \in
\R^{\ParDim}$ as the solution of some set of moment equations---say
\mbox{$\Exs_{\Sam}[ \ZFun(\Sam,\TruePar) ] = 0$,} where $\Zfun(\Sam,
\TruePar) \in \real^\usedim$, and the expectation is taken over $\Sam$
drawn according to some unknown distribution $\Pdist$.  Given samples
$\{\Sam_i \}_{i=1}^\numobs$ drawn from $\Pdist$, we can use them to
approximate the population level expectation, which leads to the
empirical analog
\begin{align}
\label{eq:z_estimate_1}
 \frac{1}{\Numobs} \sum_{i=1}^\Numobs \ZFun(\Sam_i, \Par) & = 0.
\end{align}
A solution $\EstPar_\numobs$ to these empirical moment equations is
known as a \emph{$Z$-estimate} of $\TruePar$.  The class of
$Z$-estimators includes many well-known procedures as special cases;
see~\Cref{SecExamples} for some concrete examples that we analyze in
this paper.

\subsection{From plug-in estimate to jackknife correction}

Suppose that our goal is to estimate the value $\TargetFun(\TruePar)$
of a known function \mbox{$\TargetFun: \R^{\ParDim} \mapsto \R$}
evaluated at $\TruePar$.  Such functional estimation problems arise in
various contexts, with particular cases including inference for linear
functionals in regression problems
(e.g.,~\cite{portnoy1984asymptotic,portnoy1985asymptotic,mammen1989asymptotics,lin2023semi});
inference for linear functionals in exponential
families~\cite{portnoy1988asymptotic}; estimating the treatment effect
and related quantities in causal inference
(e.g.,~\cite{robins1994estimation,robins2008higher,robins2009quadratic,hahn1998role,hirano2003efficient,su2023estimated});
inference problems involving instrumental
variables~\cite{angrist1999jackknife,cattaneo2018alternative}; as well
as other semi-parametric estimation problems
(e.g.,~\cite{bickel1993efficient,tsiatis2006semiparametric,chernozhukov2018double,kennedy2022semiparametric}).

\paragraph{Plug-in estimation and the quadratic barrier:}

A standard approach to estimating $\TargetFun(\TruePar)$ is via the
\emph{plug-in principle:} first compute the $Z$-estimate
$\EstPar_\numobs$ from the fixed point
equation~\eqref{eq:z_estimate_1}, and then form the plug-in estimate
$\TargetFun(\EstPar_\numobs)$.  In a classical analysis, one first
uses regularity conditions to argue that $\EstPar_\numobs$ is
consistent and asymptotically normal; then, as long as the functional
$\TargetFun$ is sufficiently smooth, applying the delta method yields
related guarantees for the error $\sqrt{\numobs} \big
(\TargetFun(\EstPar_\numobs) - \TargetFun(\TruePar) \big)$.

However, this favorable behavior can break down dramatically if we
undertake a more demanding analysis in which the parameter dimension
$\ParDim$ is allowed to grow with the sample size $\Numobs$.  Notably,
there are various particular cases for which the plug-in estimate is
known to exhibit a \emph{quadratic barrier} in the dimension
$\usedim$, meaning that the condition $\usedim^2/\numobs \rightarrow
0$ is necessary and sufficient for $\sqrt{\numobs}$-consistency and
normality.  For instance, in a seminal paper on the maximum likelihood estimator (MLE) in exponential
families, Portnoy~\cite{portnoy1988asymptotic} demonstrated that the
plug-estimator $\TargetFun(\EstPar_\numobs)$ based on the MLE need not
be $\sqrt{\numobs}$-consistent if $\usedim^2/\numobs$ stays bounded
away from zero.  More recently, Su et al.~\cite{su2023estimated}
exhibited the same quadratic barrier for the  inverse propensity score weighting (IPW) estimator of the
average treatment effect in causal inference.

\begin{figure}[h]
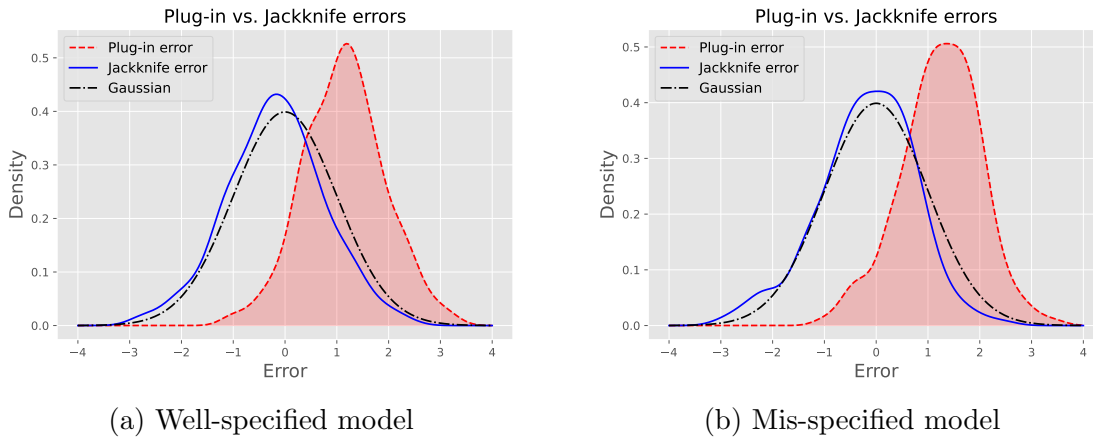

  \begin{center}
    \begin{tabular}{ccc}
      \widgraph{0.45\textwidth}{\newfigdir/fig_quad_gauss_numtrials500_n400_d20}
      &&
      \widgraph{0.45\textwidth}{\newfigdir/fig_quad_qshift_numtrials500_n400_d20}
      \\
      (a) Well-specified model && (b) Mis-specified model
    \end{tabular}
  \end{center}
  \caption{Illustration of the quadratic barrier for the plug-in
    estimator when used for estimating a non-linear functional in
    linear regression (see~\Cref{ExaQuadLinear} for details) with
    $\numobs = 400$ samples in dimension $\usedim = 20$.  Each panel
    shows histograms of the $\sqrt{\numobs}$-rescaled error for the
    plug-in estimator (red dashed lines, shaded) and the jackknife
    corrected-estimator (blue solid line).  The standard Gaussian
    density is shown in black dash-dotted lines for comparison.  (a)
    Results for a well-specified linear model.  (b) Results for a
    mis-specified linear model.  In both cases, the plug-in estimator
    exhibits a large positive bias.}
  \label{FigQuadSims}
\end{figure}

An initial contribution of this paper is characterizing a broad class
of problems\footnote{To be clear, this class does not include
estimation of linear functionals in linear models, where better
scalings are possible
(e.g.,~\cite{portnoy1984asymptotic,portnoy1985asymptotic,mammen1989asymptotics}). Another
line of work imposes sparsity constraints for linear regression to
obtain favorable dimensional scaling
(e.g.,~\cite{javanmard2014confidence,van2014asymptotically,zhang2014confidence,zhang2017simultaneous}).}
for which the quadratic breakdown of plug-in estimation is
unavoidable; see~\Cref{CorHDAsymptoticPlugIn} for a precise statement.
\Cref{FigQuadSims} provides a graphical illustration of quadratic
breakdown when the plug-in method is used to estimate a non-linear
functional in linear regression, as described in more detail
in~\Cref{ExaQuadLinear}.  As shown in these plots, for a problem with
sample size $\numobs = 400$ in dimension $\usedim = 20$ (so that
$\numobs = \usedim^2$), the plug-in estimator exhibits a large bias so
that the Gaussian approximation is highly inaccurate.

\paragraph{Debiasing approaches:}
This undesirable behavior is caused by a bias term that dominates the
stochastic fluctuations (more precisely, see the term $\highbias$ in
equation~\ref{EqnDefnStochNoise}).  Accordingly, it is natural to
consider methods that seek to estimate and remove this bias.  There is
now a substantial literature on such debiasing methods, including
methods tailored for generalized linear
models~\cite{cordeiro1991bias,kosmidis2021jeffreys}, as well as for
maximum likelihood estimators~\cite{cox1968general,firth1993bias}.
While we focus on a high-dimensional parametric setting, in the
semi-parametric literature, various forms of debiasing are based on
the von Mises expansion
(e.g.,~\cite{robins1994estimation,robins2008higher,robins2009quadratic,su2023estimated});
see Kennedy~\cite{kennedy2022semiparametric} for a survey.

Another broad class of approaches are based on resampling methods,
including the jackknife and
bootstrap~\cite{efron1982jackknife,efron1994introduction,shao2012jackknife}.
As opposed to being tailored to the specific details of the model
and/or estimator, these approaches are generic in nature, and require
only ``black box'' access---that is, one has access to a computer
program that returns the estimate for any given dataset, but all other
details of the model and estimator are unknown.  There is an evolving
literature on both jackknife debiasing
(e.g.,~\cite{parr1983note,hahn2004jackknife,shao2012jackknife,dhaene2015split,cattaneo2019two,hahn2022efficient}),
as well as bootstrap-based approaches
(e.g.,~\cite{efron1992bootstrap,efron1994introduction,hall2013bootstrap,karoui2016can,cattaneo2018kernel,cavaliere2024bootstrap}.)
With the speed of modern computing, black-box procedures are becoming
increasingly common in statistical practice, so that shedding insight
into their behavior, particularly in high-dimensional settings, is
important.

\paragraph{Jackknife debiasing breaks the quadratic barrier:}

Our second set of results provides insight into the high-dimensional
behavior of the jackknife when used both for debiasing an estimator,
and for estimating variances so as to form confidence intervals. The
use of jackknife for debiasing is simple to describe: it ``corrects''
the original estimator by subtracting a bias estimate based on the
jackknife's leave-one-out formalism.  More precisely, in application
to the $Z$-formulation and functional estimation problem formulated
here, the jackknife corrected estimate takes the form
\begin{subequations}
\begin{align}
\label{eq:jacknife_est_def}
\JackEst & \defn \TargetFun(\EstPar_\numobs) -
\frac{\Numobs-1}{\Numobs} \sum_{i=1}^{\Numobs} \big \{
\TargetFun(\LOO{\EstPar}{\numobs-1}{i}) - \TargetFun(\EstPar_\numobs)
\big \},
\end{align}
where for each $i = 1, \ldots, \numobs$, the leave-one-out (LOO)
estimate $\LOO{\EstPar}{\numobs-1}{i}$ is a solution of the estimating
equations
\begin{align}
  \label{EqnDefnLOO}
 \frac{1}{\Numobs-1} \sum_{j\neq i} \ZFun(\Sam_j, \theta) & =0.
\end{align}
\end{subequations}
Note that this $\LOO{\EstPar}{\numobs-1}{i}$ is simply the
$Z$-estimate obtained by removing the $i$-th sample, hence the
``leave-one-out'' terminology.

In the classical setting (fixed dimension $\usedim$ with sample size
$\numobs$ increasing), jackknife debiasing is well-understood, and
known to have mixed effects.  In terms of mean-squared error (MSE),
under certain smoothness conditions, it removes squared bias terms up
to terms of the order $\numobs^{-2}$.  At the same time, this bias
removal can come at a price: jackknife debiasing can inflate the
variance, and thereby lead to an estimator with higher MSE.  In
contrast, the analysis of this paper establishes that for functional
estimation under high-dimensional scaling---that is, when the
dimension $\usedim$ increases with sample size---then jackknife
debiasing can yield significant improvements over the standard
plug-in.  In particular, there is a broad class of problems for which
plug-in suffers from the quadratic barrier (requiring
$\usedim^2/\numobs \rightarrow 0$) whereas---\emph{in sharp
contrast}---jackknife debiasing yields $\sqrt{\numobs}$-consistency
and asymptotic normality under the milder scaling
$\usedim^{3/2}/\numobs \rightarrow 0$.  See~\Cref{ThmJackknife}
and~\Cref{CorHDAsymptoticDebias} for precise statements of this
guarantee, and~\Cref{FigQuadSims} for a graphical demonstration of the
good behavior of the jackknife approach for a problem with $\numobs =
\usedim^2$.

Our distributional theory (cf. Corollaries~\ref{CorHDAsymptoticPlugIn}
and~\ref{CorHDAsymptoticDebias}) provides guarantees in terms of an
unknown population variance.  Thus, fully data-dependent confidence
intervals (CIs) also require a consistent estimator of the variance.
The jackknife can also be used for variance
estimation~\cite{efron1982jackknife,shao2012jackknife}, and our second
result on the jackknife (\Cref{PropJackVarEst}) proves that this
variance estimate is consistent under the scaling
$\usedim^{3/2}/\numobs \rightarrow 0$.  Coupled with our earlier
results, we obtain CIs with guaranteed coverage in this
high-dimensional regime.  We note that jackknife variance estimation
is known to be inconsistent when $\usedim/\numobs$ converges to a
non-zero constant~\cite{karoui2016can}, so that there is a narrow
window for potential improvement to our guarantee.


\paragraph{Notation:}
We make use of $k^{th}$-order tensors $T$, which can be viewed as
multi-dimensional arrays in $(\R^{\ParDim})^{\otimes k}$.  For two
$k$-th order tensors $A$ and $B$ of matching dimensions, we define
their inner product \mbox{$\inprod{A}{B} \defn \sum_{i_1, \ldots, i_k
    \in [\ParDim]} A_{i_1 i_2, \ldots, i_k} B_{i_1 i_2, \ldots,
    i_k}$.}  When $k = 1$, then $\inprod{\cdot}{\cdot}$ is equivalent
to the standard Euclidean inner product of two vectors in
$\R^{\ParDim}$, whereas for $k = 2$, it is equivalent to the trace
inner product on matrix space.  For an $\ell$-th order tensor $A$ and
an $k$-th order tensor $B$ such that $\ell \leq k$, with slight abuse
of notation, we define $\vtinprod{A}{B}$ to be a $(k-\ell)$-th order
tensor with
\begin{align*}
  (\vtinprod{A}{B})_{i_{\ell + 1},\ldots, i_{k}} \defn \sum_{i_1,
    \ldots, i_\ell \in [\ParDim]} A_{i_{1}, \ldots, i_{\ell}}
  B_{i_{1}, \ldots, i_\ell, i_{\ell + 1}, \ldots, i_{k}}.
\end{align*}

For $u_1,\ldots,u_k\in \R^{\ParDim}$, we let $U=u_1\otimes\cdots
\otimes u_k$ be the $k$-th order tensor with elements
$U_{i_1,i_2,\ldots,i_k} = \prod_{j=1}^k u_{j,i_j}$.  For any $k$-th
order tensor $T\in(\R^{\ParDim})^{\otimes k}$, we define its operator
norm $\opnorm{T} \defn \sup_{u_1,\ldots,u_k:\enorm{u_i} = 1}
\inprod{T}{u_1 \otimes \cdots \otimes u_k}$, and write
$T[u_1,\ldots,u_ \ell] \defn \vtinprod{T}{u_1 \otimes \cdots \otimes
  u_\ell }$ for any $\ell \geq 1$. For any twice continuously
differentiable vector-valued function $f(\Par):\R^{\ParDim} \mapsto
\R^{\ParDim}$, its gradient $\nabla_\Par f(\Par)$ and Hessian $\nabla_\Par^2
f(\Par)$ are second and third order tensors with $(\nabla_\Par
f(\Par))_{ij} = \partial_{x_j} f_i(\Par)$ and $(\nabla_\Par^2
f(\Par))_{ijk} = \partial_{x_k}\partial_{x_j} f_i(\Par)$ for any
$i,j,k\in[\ParDim]$, respectively. We occasionally omit the subscript $\Par$ in $\nabla_\Par$ when it is clear from the context.

In this paper, we use $c$ (or $c'$) to denote universal constants. We
use $\polyshort \equiv \polyshort(X,Y,Z,\ldots)$ (or $\polyshortprime
\equiv \polyshortprime(X,Y,Z,\ldots)$) to denote constants that depend
polynomially in $X$, $Y$, $Z$, $\ldots$ but are \emph{independent} of
$\Numobs,\ParDim$. \mbox{ We allow the values of $c$ (or $c'$) and
  $\polyshort$ (or $\polyshortprime$) to vary from place to place. }
In addition, we use $\lesssim$ (or $\gtrsim$) to denote inequalities
that hold up to logarithmic factors. For two sequences
$a_n,b_n\to\infty$, we use $a_n/b_n\highdimto0$ to denote $a_n^{1 +
  \delta}/b_n \rightarrow 0$ for some positive $\delta>0.$


\section{Main results}
\label{sec:main-results}

In this section, we present our main results and discuss some of their
consequences. Our focus is non-asymptotic and asymptotic guarantees of
the plug-in and the jackknife estimators. Our main results are
finite-sample bounds on the $\sqrt{\numobs}$-rescaled errors
\begin{align*}
\DelPlug \defn \sqrt{\numobs} \big \{ \tau(\thetahat_\numobs) -
\tau(\thetastar) \}, \quad \mbox{and} \quad \DelJack \defn
\sqrt{\numobs} \big \{ \JackEst - \tau(\thetastar) \}
\end{align*}
of the plug-in and jackknife corrected
estimators~\eqref{eq:jacknife_est_def}, respectively, along with
asymptotic guarantees derived from these finite-sample results.

We begin in~\Cref{SecExamples} by specifying the conditions that
underlie our analysis, along with discussion of various motivating
examples that illustrate these conditions.  \Cref{SecQuant} is devoted
to our main results, with finite-sample guarantees given
in~\Cref{ThmPlugin} and~\Cref{ThmJackknife}.
In~\Cref{subsec:asymptotic}, we use our finite-sample results to study
the asymptotics of both the plug-in and jackknife-adjusted estimators
under a form of high-dimensional scaling, in which we allow the
problem parameters to scale with the sample size.  This analysis
provides a clear justification for the use of the jackknife-corrected
estimator in certain high-dimensional regimes.  Finally, in
\Cref{SecCI}, we provide guarantees for jackknife variance estimation,
along with the confidence intervals whose validity follows from these
guarantees.


\subsection{Regularity conditions and examples}
\label{SecExamples}

Throughout the paper, we assume that the the parameter space
$\ParSpace$ is a convex, bounded and open set in $\R^{\ParDim}$.
Given the boundedness assumption, there is some finite $\ConParNorm
\equiv \ConParNorm(\ParSpace)$ such that
\begin{align}
\label{EqnDefnConParNorm}
\enorm{\Par} \leq \ConParNorm \quad \mbox{for all $\Par
  \in \ParSpace$,}
\end{align}
and the radius $\ConParNorm$ appears in various pre-factors in our
non-asymptotic results.

\paragraph{Smoothness conditions:}
From classical theory on the jackknife (e.g., Shao and
Tu~\cite[Section 2.3]{shao2012jackknife}), the jackknife need not
behave well unless the $Z$-function $\ZFun$ and the target functional
$\TargetFun$ have appropriate smoothness properties.  In this paper,
we impose the following Lipschitz conditions:
\myassumption{$L$-smooth}{ass:smooth_function}{The function $\tau:
  \real^\ParDim \rightarrow \real$ is three times continuously
  differentiable, and there exists a constant $\PlainLip$ such that
\begin{subequations}  
  \begin{align}
 \abss{\TargetFun(\Par_1) - \TargetFun(\Par_2)} \leq \LipTara
 \enorm{\Par_1 - \Par_2}, \qquad & \vecnorm{
   \nabla_{\Par}\TargetFun(\Par_1) - \nabla_{\Par}
   \TargetFun(\Par_2)}{2} \leq \LipTarb\enorm{\Par_1 - \Par_2}, \;
 \mbox{and} \qquad \\
\opnorm{ \nabla_{\Par}^2 \TargetFun(\Par_1) - \nabla_{\Par}^2
  \TargetFun(\Par_2)} & \leq \LipTarc\enorm{\Par_1 - \Par_2}.
  \end{align}
\end{subequations}
for all $\Par_1, \Par_2 \in \ParSpace$.  }

\noindent \paragraph{Tail conditions:} In addition, we require certain
tail conditions on the $Z$-function $\ZFun$ and its derivatives.  For
a parameter $\alpha > 0$, we state these conditions in terms of the
\emph{$\alpha$-Orlicz norm}
\begin{align}
\label{EqnOrlicz}  
\rnorm{X}{\psi_{\alpha}} & \defn \inf \big \{ \lambda>0 \mid
\E[\psi_{\alpha}(|X|/\lambda)] \leq 1 \big \} \qquad \mbox{where
  $\psi_\alpha(x) \defn \exp(x^\alpha)-1$.}
\end{align}
For $\alpha = 1$, we obtain sub-exponential random variables, whereas
setting $\alpha = 2$ corresponds to sub-Gaussian random variables.  We
say that a random vector $Y \in \real^{\ParDim}$ is
$\sigma$-sub-exponential if $\rnorm{\inprod{u}{Y}}{\psi_1} \leq
\sigma$ for all unit-norm vectors $u$, and write this condition as
$\rnorm{Y}{\psi_1} \leq \sigma$. We generalize this notion to
higher-order tensors in the natural way.

\myassumption{$\sigma$-Tail}{ass:tail}{Suppose that: (a) For each
  $\Sam$, the function $\theta \mapsto \ZFun(\Sam, \theta)$ is three
  times continuously differentiable; (b) there exists $\sigma <
  \infty$ such that
  \begin{align*}
    \rnorm{\Zfun(Z, \theta)}{\psi_1} \stackrel{(i)}{\leq} \sigma,
    \quad \mbox{and} \quad \rnorm{\nabla_\theta \Zfun(Z,
      \theta)}{\psi_1} \stackrel{(ii)}{\leq} \sigma \qquad \mbox{for
      each $\theta \in \ParSpace$,}
  \end{align*}
  and (c) for all unit vectors $\Direc, \Direcb, \Direcc, \Direcd \in
  \Sphere{\ParDim-1}$:
  \begin{align*}    
    \rnorm{\sup_{\Par\in \ParSpace}| \nabla^2_{\Par} \ZFun(\Sam,
      \Par)[\Direc,\Direcb,\Direcc]|}{\psi_{\OrConc}} \leq
    \simpleorpar \quad \mbox{and} \quad {\sup_{\Par \in \ParSpace}|
      \nabla^3_{\Par} \ZFun(\Sam,
      \Par)[\Direc,\Direcb,\Direcc,\Direcd]|} \leq \Factora(\Sam,
    \Direc) \Factorb(\Sam, \Direc,\Direcb,\Direcc,\Direcd)
  \end{align*}
  for some functions $\Factora$ and $\Factorb$ such that
  $\rnorm{\Factora(\Sam, \Direc)}{\psi_1} \leq 1$ and
  $\rnorm{\Factorb(\Sam,
    \Direc,\Direcb,\Direcc,\Direcd)}{\psi_{\OrCond}} \leq
  \simpleorpar$.  }

\noindent Condition~\ref{ass:tail} is closely related to but slightly
stronger than imposing Lipschitz conditions\footnote{More precisely,
from the sub-exponential condition (ii) in part (b) on the Jacobian
$\nabla \Zfun(\Sam, \Par)$, we have $\sup_{\theta \in \Par}
\opnorm{\Exs \nabla \Zfun(\Sam, \Par)} \leq c \sigma$, from which a
Taylor series argument guarantees that
\mbox{$\vecnorm{\Exs[\ZFun(\Sam, \Par_1)]- \Exs[\ZFun(\Sam,
      \Par_2)]}{2} \leq c\OrParb\vecnorm{\Par_1 - \Par_2}{2}$.}
Similar arguments apply to the derivatives of $\Zfun$.} on the
population function $H(\Par) \defn \Exs[\Zfun(\Sam, \Par)]$. \\

Our final requirement is that the $Z$-estimate is consistent in a
certain sense.  This \emph{convergence condition} involves the
population $Z$-function $\PopZfun(\Par) \defn \Exs[\ZFun(\Sam,
  \Par)]$, along with its Jacobian
\begin{align}
\label{EqnDefnJmat}
\Jmat \defn \nabla \PopZfun(\TruePar) \; \stackrel{\diamond}{=} \;
\Exs \big[ \nabla \ZFun(\Sam, \TruePar) \big].
\end{align}
evaluated at the true parameter.  In asserting the equivalence
$(\diamond)$, we have made use of standard arguments to interchange
the order of integration and differentiation.

\myassumption{$\strongconvex$-Conv}{ass:convergence}{Suppose that the
  minimum singular value of $\Jmat$ is lower bounded as
  $\sigma_{\min}(\Jmat) \geq \strongconvex$ for some $\strongconvex >
  0$, and that the $Z$-estimate $\EstPar_{\Numobs}$ from
  equation~\eqref{eq:z_estimate_1} is consistent in the sense that
\begin{align}
    \enorm{\EstPar_{\Numobs}-\TruePar} & \leq \polyshort \sqrt{\log
      \Numobs}\sqrt{\frac{\ParDim + \log(1/\delta)}{\Numobs}} \qquad
    \mbox{with probability at least $1-\delta$}
\end{align}
for some $\polyshort = \poly$. }

A sufficient condition for~\ref{ass:convergence} is the following \emph{concavity condition}.
\myassumption{Con}{ass:concavity}{ The matrix $-\nabla \PopZfun(\Par)$
  is positive semidefinite for all $\Par \in \ParSpace$ and strictly
  positive definite for $\TruePar$---namely, we have $-\nabla
  \PopZfun(\TruePar) \succeq \strongconvex\IdMat_{\ParDim}$ for some
  constant $\strongconvex>0$.  }
  We refer to Appendix~\ref{sec:suff_cond_of_conv} for further details.\\


\noindent Let us consider some canonical examples of $Z$-estimators to
illustrate these conditions:
\begin{example}[$Z$-estimators for linear models]
\label{ExaLinearRegression}  
Given a covariate vector $X \in \real^\usedim$ and an associated
response variable $Y \in \real$, the goal of linear regression is to
determine the optimal linear predictor $x \mapsto \widehat{y} =
\inprod{x}{\thetastar}$, where optimality is typically measured in
terms of the mean-squared error (MSE).  A broad class of methods can
be formulated as $Z$-estimators.  More precisely, given $\numobs$
i.i.d. pairs $\{ (\covar_i, Y_i) \}_{i=1}^\numobs$, consider an
estimate defined by the fixed point equation
\begin{align}
\label{EqnLinearZ}  
\frac{1}{\numobs} \sumn \covar_i \Zfunexam(Y_i-
\inprod{\covar_i}{\theta}) = 0,
\end{align}
where $\Zfunexam: \real \rightarrow \real$ is a given function.  For
instance, with the choice $\Zfunexam(t) = t$, the $Z$-estimator
defines the optimality condition associated with the ordinary least
squares estimator.  More generally, if we compute the regression
solution using a differentiable loss function $\RegLoss: \real
\rightarrow \real$, then the zero-gradient conditions associated with
an optimum are of the form~\eqref{EqnLinearZ} with $\Zfunexam(t) \defn
\RegLoss'(t)$.  We recover least-squares with the loss function
$\RegLoss(t) = t^2/2$.

The family of estimators~\eqref{EqnLinearZ} are a special case of our
general set-up with the $Z$-function
\begin{align}
\ZFun(\Sam_i, \Par) & \defn \covar_i \Zfunexam(Y_i-
\inprod{\covar_i}{\Par}).
\end{align}
Assumption~\ref{ass:tail}~and~\ref{ass:concavity} are satisfied if we
impose suitable conditions on $\Zfunexam$ and its derivatives, as well
as tail conditions on the pair $(\covar, Y)$.  In the case of ordinary
linear regression ($\Zfunexam(t) = t$), we can compute the derivatives
$\nabla_\theta \ZFun(\Sam_i, \Par) = - \covar_i \covar_i^{\mytrans}$
and $\nabla^2_\theta \Zfun(\Sam_i, \Par) = 0$, so that it is
relatively easy to verify condition~\ref{ass:tail} when the pair
$(\covar, Y)$ are sub-Gaussian.  Moreover, there are a variety of
settings under which linear regression satisfies the convergence
condition~\ref{ass:convergence}.

For more general choices of the function $\Zfunexam$, we need to
impose conditions on its deriviatives; see~\Cref{sec:linear_verify}
for more details of these calculations.  As a sidenote, our smoothness
conditions do not apply to the Huber loss~\cite{huber1973robust},
commonly used to enforce robustness, since it is not second-order
continuously differentiable.  However, this issue can be side-stepped
via the use of a smooth version known as pseudo-Huber
loss~\cite{charbonnier1994two}, which satisfies our conditions.

\goodendex
\end{example}

\vspace*{0.1in}

\begin{example}[Logistic regression and GLMs]
\label{ExaGLM}
For a binary response $Y \in \{0,1\}$, suppose that we model the log
odds as a linear function of a covariate vector $\covariate \in
\real^\usedim$---that is $\log \frac{\Prob_\thetastar(Y = 1 \mid X =
  x)}{\Prob_\thetastar(Y = 0 \mid X = x)} = \inprod{\thetastar}{x}$.
The maximum likelihood estimator (MLE) based on samples $\{(X_i, Y_i) \}_{i=1}^\numobs$ is defined by
the stationary condition
\begin{align}
\label{EqnLogisticRegression}  
\frac{1}{\numobs} \sum_{i=1}^\numobs X_i \big \{ Y_i - \varphi \big(
\inprod{\theta}{X_i} \big) \big \} = 0 \qquad \mbox{where} \quad
\varphi(t) \defn \frac{e^t}{1 + e^t}.
\end{align}
Thus, logistic regression is a special case of a $Z$-estimator with
$\Zfun(\Sam; \theta) \defn X \big \{ Y - \varphi
\big(\inprod{\theta}{X} \big) \big \}$.  More generally, the MLE for
any generalized linear model (GLM) with canonical link function is
equivalent to a $Z$-estimator of the
form~\eqref{EqnLogisticRegression} with
\begin{align*}
  \varphi(t) = \Phi'(t) \quad \mbox{where} \quad \Phi(t) \defn \log
  \int_{\mathcal{Y}} e^{t y} \mu(dy) \qquad \mbox{for some measure
    $\mu$ on $\mathcal{Y}$.}
\end{align*}
We recover logistic regression with $\Phi(t) = \log(1 + e^t)$ and
$\mu$ being the counting measure on $\mathcal{Y} = \{0, 1\}$.
\goodendex
\end{example}

\begin{example}[ IPW estimator for the average treatment effect]
\label{ExaIPW}
There are various application domains---including medical trials and
public policy---in which it is of interest to assess the effect of
applying some ``treatment''.  This problem is challenging when the
data is \emph{not} based on randomized controlled trials, but instead
is observational in nature.

More concretely, let $\Action \in \{0, 1 \}$ be a binary indicator
variable for ``treated'' versus ``non-treated'' status, and let
$\State \in \real^{\usedim-1}$ be a vector of covariates.\footnote{Our
reason for choosing dimension $\usedim - 1$ will become clear
momentarily.} We use $\Outcome(a)$ to denote a potential outcome
associated with treatment status $a \in \{0,1 \}$, and define the
outcome function $\mu(\State, \action) \defn \Exs \big[
  \Outcome(\action) \mid \State \big]$.  Letting $\pi(x) =
\Prob[\Action = 1 \mid \State = x]$ denote the conditional probablity
of treatment, also known as the propensity function.  Suppose that we
observe triples $(\Action, \State, \Outcome)$ with $\Action \sim
\pi(\State)$, and outcome $\Outcome = \Outcome(\Action)$.  Our goal is
to estimate the average treatment effect $\taustar \defn \Exs[
  \Outcome(1) - \Outcome(0) ]$, or ATE for short.

Under standard unconfoundedness assumptions, the ATE is identified by
the population relation $\taustar = \Exs \left[ \frac{ \Action
    \Outcome}{\truepropscore(\State)} - \frac{(1 - \Action)
    \Outcome}{1 - \truepropscore(\State)} \right]$.  Thus, given
i.i.d. samples $\{(\Action_i, X_i, Y_i)\}_{i=1}^\numobs$, we can
consider the inverse propensity score weighting (IPW)
estimator~\cite{horvitz1952generalization}, which takes the form
  \begin{align}
\label{EqnDefnIPW}    
    \tauhatipw = \frac{1}{\numobs}\sum_{i=1}^\numobs \Biggr[ \frac{
        \Action_i \Outcome_i}{\estedpropscore(\State_i)} - \frac{(1 -
        \Action_i) \Outcome_i}{1 - \estedpropscore(\State_i)} \Biggr],
\end{align}
where $\estedpropscore$ is an estimate of the propensity score.

Different IPW variants are obtained by imposing modeling assumptions
on the propensity score, and deriving suitable estimators based on
these assumptions.  For example, we might adopt the logistic model
described in~\Cref{ExaGLM}; with $\varphi$ denoting the logistic
function from equation~\eqref{EqnLogisticRegression}, write $\pi(x) =
\varphi\big(\inprod{\betastar}{x} \big)$ for some unknown parameter vector
$\betastar \in \real^{\usedim-1}$.  With this choice, the
zero-gradient condition defining the MLE $\betahat \in
\real^{\usedim-1}$ takes the form given in
equation~\eqref{EqnLogisticRegression}.

We can describe this two-stage procedure as a special case of a
$Z$-estimator as follows.  Introduce the full parameter vector
$\TruePar = ( \betastar,\taustar) \in \R^{\ParDim}$.  Our goal is to
estimate the first co-ordinate of $\TruePar$, which corresponds to the
linear functional $\TargetFun(\TruePar) \defn
\inprod{e_\ParDim}{\TruePar}$, where $e_\ParDim \in \real^\ParDim$
denotes the standard basis vector with one in the last entry.  Then
the estimate $\thetahat \defn ( \betahat,\tauhatipw)$ corresponds to
the $Z$-estimate defined by the function
\begin{align}
\label{eq:zfun_ipw_example}  
\ZFun(\Sam, \Par) & \defn
\begin{bmatrix}
\covar \big \{ \act - \varphi(\inprod{\beta}{\State}) \big \} \\
\frac{\Action \Outcome}{ \varphi (\inprod{\beta}{\State})} - \frac{(1
  - \Action) \Outcome}{1 - \varphi (\inprod{\beta}{\State}) } - \tau
\end{bmatrix} \qquad \mbox{with $\theta = ( \beta,\tau)$.}
\end{align}
See~\Cref{sec:ipw_verify} for verification of conditions under which
this $Z$-estimator satisfies our regularity conditions.
\goodendex
\end{example}

\newcommand{\xmacro}{{2}}
\newcommand{\wmacro}{{k}}
\newcommand{\xcoordinate}{{\mathrm X}}
\newcommand{\xvector}{{ X}}
\newcommand{\ivinceptpar}{{\alpha}}
\newcommand{\ivtargetpar}{{\beta}}
\newcommand{\ivinceptparstar}{{\alpha^*}}
\newcommand{\ivtargetparstar}{{\beta^*}}
\newcommand{\ivinceptparest}{{\widehat\alpha}}
\newcommand{\ivtargetparest}{{\widehat\beta}}
\newcommand{\ivwithoutincept}{{\mathrm W}}
\newcommand{\TSLS}{\operatorname{TSLS}}

\begin{example}[Instrumental variables estimate]
\label{ExamIV}
Consider a linear regression model of the form
\begin{subequations}
\begin{align}
\label{EqnOverOne}
Y & = \ivinceptparstar + \xvector \ivtargetparstar + \noise
\end{align}
with scalar response $Y \in \real$, covariate $\xvector \in \real$ and
noise variable $\noise$.  In the classical setting, the covariate is
uncorrelated with the noise (i.e., $\Exs[ \xvector \varepsilon] = 0$),
in which case it is said to be exogeneous.  When this uncorrelatedness
fails to hold, the covariate is said to be \emph{endogeneous}; with an
endogeneous covariate, standard estimates of the parameters
$\thetastar \defn (\alphastar, \betastar) \in \real^2$, such as
ordinary least-squares, are no longer consistent.  Endogeneous linear
models arise in various settings, including omitted variable bias
(arising when relevant covariate is omitted from a well-specified
linear model); regression problems with missing or noisy observations
of covariates; and in the class of TD estimators used in reinforcement
learning (e.g.,~\cite{KhaPanRuaWaiJor20,SutBar18}).

An instrument $\instrument \in \real^\wmacro$ (with some dimension
$\wmacro \geq 2$) is random vector such that $\Exs[\instrument
  \begin{bmatrix} 1 & \xvector 
\end{bmatrix}] \neq 0$ and $\Exs[\instrument \noise] = 0$.
In terms of the shorthand $\theta \defn (\alpha, \beta) \in \real^2$,
we can use the instrument to define the $Z$-estimator
  \begin{align*}
  \Zfun(\Sam, \theta) & = \instrument Y - \instrument (\alpha + X
  \beta) \qquad \mbox{where $Z = (X, \instrument, Y)$.}
  \end{align*}
  Recalling that $\thetastar = (\alphastar, \betastar)$ are the true
  parameters, note that we have $\Exs[h(\Sam, \thetastar)] = 0$ by our
  assumptions on the instrument.

In the \emph{over-specified setting}, this $Z$-estimator will not have
a unique solution, and a standard approach is to specify a unique
estimate via two-stage least-squares.  In particular, without loss of
generality (by projecting $\xvector$ onto the span of the
instruments), we can always write
\begin{align}
\label{EqnOverTwo}
  \xvector = \inprod{\instrument}{\ivcoeffstar} + \iverror
\end{align}
some vector $\ivcoeffstar\in\R^{\wmacro}$ and zero-mean noise variable
such that $\Exs[\instrument \iverror] = 0$. We then define the
augmented parameter vector $\thetastar = (\alphastar, \betastar,
\ivcoeffstar) \in \real^{\usedim}$ where $\usedim \defn \wmacro + 2$,
and consider the $Z$-estimator defined by the function
\begin{align}
\label{EqnOverThree}  
  \smallsub{\ZFun}{TSLS}(\Sam, \Par) \defn \begin{bmatrix} \instrument
    \xvector - \instrument\instrument^\top \ivcoef\\ \ivcoef^\top
    \instrument Y - \ivcoef^\top \instrument(\xvector \ivtargetpar +
    \ivinceptpar) \\ Y - (\xvector \ivtargetpar +
    \ivinceptpar) \end{bmatrix} \in \real^{\usedim} \qquad \mbox{where
    $\Sam = (\xvector, \instrument, Y)$.}
\end{align}
\end{subequations}
This is a linear $Z$-estimator, and it is straightforward to verify
that our conditions hold under appropriate sub-Gaussian conditions on
$(\xvector, \ivwithoutincept)$.  \hfill \goodendex
\end{example}


\subsection{Non-asymptotic bounds and their consequences}
\label{SecQuant}

Having set up and motivated the class of $Z$-estimators, we now turn
to some analysis, beginning with the plug-in estimator
$\tau(\thetahat_\numobs)$.

\subsubsection{Quantitative bounds for the plug-in estimator}
\label{SecQuantPlugin}
The main result of this section is to provide an explicit and
non-asymptotic bound on the difference between the
$\sqrt{\numobs}$-rescaled plug-in error \mbox{$\DelPlug \defn
  \sqrt{\numobs} \big \{ \tau(\thetahat_\numobs) - \tau(\thetastar)
  \big \}$} and a sum of the form $\StochNoise_\numobs + \frac{1}{
  \sqrt{\numobs}} \highbias$.  Here $\StochNoise_\numobs$ is a
zero-mean random variable, whereas $\highbias$ is a deterministic bias
term.  Although the bias term $\highbias$ does not depend on the
sample size, it typically increases as a function of the dimension
$\usedim$, as reflected in our choice of notation.

Let us now specify the stochastic noise and bias terms.  In order to
do so, recall from equation~\eqref{EqnDefnJmat} that the matrix
$\GradFun_\TruePar \defn \Exs \big[ \nabla_\theta \ZFun(\Sam,
  \TruePar) \big]$ is assumed to be invertible. We then define the
$\ParDim$-dimensional matrix\footnote{To be clear, the tensor product
$\vtinprod{\IntSec_\TruePar}{ \nabla^2_\Par \ZFun(\Sam,\TruePar)}$ in
equation~\eqref{EqnDefnMmat} is a $\ParDim \times \ParDim$ matrix with
\mbox{$(\vtinprod{\IntSec_\TruePar}{ \nabla^2_\Par
    \ZFun(\Sam,\TruePar)})_{ij} = \sum_{k=1}^\ParDim
  (\IntSec_{\TruePar})_{k} (\nabla^2_\Par
  \ZFun_k(\Sam,\TruePar))_{ij}$.}} $\Mmat$ and vector
$\MyIntSec{\TruePar} \in \real^\ParDim$ as
\begin{subequations}
\begin{align}
\label{EqnDefnMmat}  
\MyIntSec{\TruePar} \defn - \GradFun_\TruePar^{-1} \nabla
\TargetFun(\TruePar), \quad \mbox{and} \quad
\Mmat \defn \nabla^2\TargetFun(\TruePar) - \Exs
      [\vtinprod{\IntSec_\TruePar}{ \nabla^2_\Par
          \ZFun(\Sam,\TruePar)}].
\end{align}
Moreover, we introduce the random vectors
\begin{align}
\label{EqnSidefn}  
\sionebase(\Sam) \defn - \Big[ \GradFun_\TruePar - \nabla_{\Par}
  \ZFun(\Sam,\TruePar) \Big] \MyIntSec{\TruePar} \quad \mbox{and}
\quad
  \sitwobase(\Sam) \defn \GradFun_{\TruePar}^{-1} \ZFun(\Sam,
  \TruePar),
\end{align}
\end{subequations}
and observe that $\sitwobase(\Sam)$ is zero-mean by definition of the
$Z$-estimator, whereas $\sionebase(\Sam)$ is zero-mean by the
definition of $\Jmat$.

With this notation, the \emph{stochastic noise term}
$\StochNoise_\numobs$ and \emph{deterministic bias term} $\highbias$
are given by
\begin{subequations}
\begin{align}
\label{EqnDefnStochNoise}    
\StochNoise_\numobs \defn \frac{1}{ \sqrt{\Numobs}} \sum_{i=1}^\Numobs
\inprod{\IntSec_\TruePar}{\ZFun(\Sam_i, \TruePar)}, \quad \mbox{and}
\quad \highbias \defn \Exs \Big[
  \inprod{\sionebase(\Sam)}{\sitwobase(\Sam)} + \tfrac{1}{2}
  \binprod{\sitwobase(\Sam)}{\Mmat\sitwobase(\Sam)} \Big].
\end{align}
Observe that we have
\begin{align}
\label{EqnDefnVarTrue}
\VarTruePlain \defn \var \big(\StochNoise_\numobs \big) & = (\nabla
\TargetFun(\TruePar))^{\mytrans} \Jmat^{-\mytrans}
\cov_\Prob(\ZFun(\Sam, \TruePar)) \Jmat^{-1} \nabla
\TargetFun(\TruePar),
\end{align}
\end{subequations}
so that the covariance of $\StochNoise_\numobs$ is consistent with the
classical asymptotic prediction.

Our results apply to all sample sizes above a certain minimal level,
which depends on a user-defined quantity $\delta \in (0,1)$ that
controls the probability of failure.  In particular, we assume that
there exists a constant $\polyshort = \poly$ such that
\begin{align}
\label{eq:linear_approx_Z_est_sample_ahead}
\Numobs/ \log^3(\Numobs/\delta) \geq \polyshort (\ParDim + \log
(\Numobs/\delta))^{3/2}.
\end{align}
\begin{theorem}
\label{ThmPlugin}
Given a sample size $\Numobs$ satisfying the lower
bound~\eqref{eq:linear_approx_Z_est_sample_ahead}, and under the
conditions~\ref{ass:smooth_function}, ~\ref{ass:tail} and
~\ref{ass:convergence}, there is a constant $\polyshortprime =
\polyprime$ such that
\begin{align}
\label{EqnPlugin}
\Biggr| \sqrt{\numobs} \big \{ \TargetFun(\EstPar_\numobs) -
\TargetFun(\TruePar) \big \} - \big \{ \StochNoise_\numobs + \frac{1}{
  \sqrt{\numobs}} \highbias \big \} \Biggr| & \leq \polyshortprime \,
\Big\{ \sqrt{ \frac{\ParDim}{\numobs}} + \frac{\ParDim^{3/2}}{\numobs}
\Big\} \log^{5}(\numobs/\delta)
\end{align}
with probability at least $1 - \delta$.
\end{theorem}
\noindent See~\Cref{SecProofThmPlugin} for the proof
of~\Cref{ThmPlugin}. \\

Let us make a few comments on this claim.  Note that the RHS of the
bound~\eqref{EqnPlugin} converges to zero as along as
$\frac{\Pardim^{3/2}}{\numobs} \highdimto 0$.  Consequently, in this
regime, \Cref{ThmPlugin} guarantees that the $\sqrt{\numobs}$-rescaled
plug-in error is well-approximated by the sum $\StochNoise_\numobs +
\frac{1}{\sqrt{\numobs}} \highbias$.  From the
definition~\eqref{EqnDefnStochNoise} and the given assumptions, note
that we are guaranteed that $\StochNoise_\numobs \convdist \cN (0,
\TrueVarPlain^2)$, where the variance $\TrueVarPlain^2$ was defined in
equation~\eqref{EqnDefnVarTrue}.  Thus, when the bias term is
negligible, the bound~\eqref{EqnPlugin} recovers the behavior expected
based on standard asymptotics and the delta method
(e.g.,~\cite{van2000asymptotic}).  In the classical analysis (with
dimension $\usedim$ and the problem sequence fixed), it immediately
follows that $\highbias/\sqrt{\numobs} \rightarrow 0$, so that the
bias is negligible.

When we allow the dimension $\usedim$ to scale, the central question
to address is the scaling of the bias term $\highbias$.  There are
various settings in which it scales linearly in dimension, with one
simple example given below.  Under this linear scaling, we have
$\highbias/\sqrt{\numobs} \asymp \Pardim/\sqrt{\numobs}$, so that we
require $\Pardim^2/\numobs \rightarrow 0$ in order for the rescaled
plug-in error to be $\sqrt{\numobs}$-consistent. We make this
intuition precise in~\Cref{CorHDAsymptoticPlugIn} to follow
in~\Cref{subsec:asymptotic}. \\

\noindent Let us consider a simple example in which the bias term
scales linearly in the dimension.

\begin{example}[Quadratic functionals in linear regression]
\label{ExaQuadLinear}
Recall the problem of linear regression introduced
in~\Cref{ExaLinearRegression} with $\Sam = (X,Y)$.  Beginning with the
case of a well-specified linear model, we can write $Y =
\inprod{X}{\thetastar} + \noise$ for some zero-mean noise $\noise$
such that $\Exs[\noise \mid X] = 0$.  Let us consider the behavior of
the ordinary-least squares estimator (obtained
from~\Cref{ExaLinearRegression} with $\Zfunexam(t) = t$), and the
problem of estimating some non-linear functional.  A standard example
is a quadratic functional $\tau(\theta) \defn \frac{1}{2}
\theta^{\mytrans} \Qmat \theta$ for some known matrix $\Qmat$; as
particular case, the choice $\Qmat = \IdMat$ corresponds to the
problem of estimating the squared norm of the regression vector.  This
is of interest in hypothesis testing, among other applications.

In the $Z$-formulation of ordinary least squares, we have $\Zfun(\Sam,
\theta) = X (Y - \inprod{X}{\theta})$, from which some straightforward
calculations yield $\Jmat = \Exs[ X X^{\mytrans}]$, $\Mmat = \nabla^2
\TargetFun(\thetastar)$, and
\begin{align*}
\sionebase(\Sam) = \big[ \Jmat - X X^{\mytrans} \big] \Jmat^{-1}
\nabla \TargetFun(\thetastar), \quad \mbox{and} \quad \sitwobase(\Sam)
= \Jmat^{-1} X \noise.
\end{align*}
From the conditional independence $\Exs[\noise \mid X] = 0$, it
follows that $\Exs \big[ \inprod{\sionebase(\Sam)}{\sitwobase(\Sam)}
  \big] = 0$.  Putting together the pieces, the bias term takes the
form
\begin{align}
\label{EqnQuadBias}  
\highbias & = \tfrac{1}{2} \Exs \Big[ \sitwobase^{\mytrans}(\Sam)
  \nabla^2 \TargetFun(\thetastar) \sitwobase(\Sam) \Big] \; = \;
\tfrac{1}{2} \trace \Big( \Exs[\noise^2 X X^{\mytrans}] \,
\Jmat^{-\mytrans} \nabla^2 \TargetFun(\thetastar) \Jmat^{-1} \Big).
\end{align}
Consequently, as long as $\nabla^2 \TargetFun(\thetastar) \neq 0$, the
bias term is the trace of a $\usedim$-dimensional matrix: it scales
linearly in dimension in a generic setting.\footnote{Of course, one
could construct examples where this matrix is near low-rank, in which
the scaling could be sub-linear in $\usedim$, but this is not the
generic case.}

For future reference, it is also worthwhile considering the effect of
mis-specification in the linear model, in which case we have
$\Exs[\noise X] = 0$, but $\Exs[\noise \mid X] \neq 0$.  In this case,
a little calculation yields
\begin{align*}
  \Exs \big[ \inprod{\sionebase(\Sam)}{\sitwobase(\Sam)} \big] = \Exs
  \Big[\noise \big( X^{\mytrans} \Jmat^{-1} \big[ \Jmat - X
      X^{\mytrans} \big] \Jmat^{-1} \nabla \TargetFun(\thetastar)
    \big) \Big].
\end{align*}
When $\Exs[\noise \mid X ] \neq 0$, this term need not be non-zero, so
that there can be other bias components---in addition to the previous
term~\eqref{EqnQuadBias}---under mis-specification.  \goodendex
\end{example}

\subsubsection{Quantitative bounds for the jackknife estimator}
\label{SecQuantJack}

Thus far, we have argued that the plug-in estimator is
$\sqrt{\numobs}$-consistent and agrees with the asymptotic prediction
as long as $\Pardim^2/\numobs \rightarrow 0$.  In this section, we
show that jackknife correction yields a strict improvement, in that
the favorable behavior is retained as long as, up to logarithmic
factors, the ratio $\ParDim^{3/2}/\numobs$ converges to zero.  We do
so by giving a non-asymptotic bound between the
$\sqrt{\numobs}$-rescaled jackknife error $\DelJack \defn
\sqrt{\numobs}\{\JackEst - \TargetFun(\TruePar) \}$ and the stochastic
noise term $\StochNoise_\numobs$ from
equation~\eqref{EqnDefnStochNoise}.
\begin{theorem}
\label{ThmJackknife}
Under the conditions of~\Cref{ThmPlugin}, there is a constant
$\polyshortprime = \polyprime$ such that the jackknife corrected
estimate $\JackEst$ satisfies the bound
\begin{align}
\label{EqnJackknife}
\Big| \sqrt{\numobs}\{\JackEst - \TargetFun(\TruePar) \} -
\StochNoise_\numobs \Big| & \leq \polyshortprime \ \; \Big \{
\sqrt{\frac{\ParDim}{\numobs}} + \frac{\ParDim^{3/2}}{\numobs} \Big \}
\log^{13/2}(\Numobs/\delta)
\end{align}
 with probability at least $1 - \delta$.
\end{theorem}
\noindent See~\Cref{SecProofThmJackknife} for the proof
of~\Cref{ThmJackknife}. \\

A few remarks regarding the claim are in order.
Theorem~\ref{ThmJackknife} shows that when $\ParDim^{3/2}/\numobs
\highdimto 0$, then the jackknife error $\DelJack$ is
well-approximated by the stochastic noise term $\StochNoise_\numobs$
alone---the bias term $\highbias$ no longer appears.  Thus, the
jackknife procedure no longer suffers from the quadratic barrier that
characterizes the plug-in approach.  This provides strong motivation
for its use in regimes where $\ParDim^2/\numobs$ does not converge to
zero, and our later simulation results confirm its benefits in
practice.

At a high level, as shown in the proof, the bias term $\highbias$
arises as the non-zero mean of $U$-statistic when the error is
expanded.  In the classical analysis (e.g., ~\cite[Theorem
  2.1]{efron1982jackknife}, \cite[Example 1.3]{shao2012jackknife}), it
is relatively straightforward to show how the jackknife removes this
quantity.  In contrast, establishing the non-asymptotic
guarantee~\eqref{EqnJackknife} requires a significantly more delicate
analysis, including some new technical results of broader interest
that we develop.

In particular, let us highlight here a key result---essential to the
argument and challenging to prove---that gives non-asymptotic bounds
on various functions of leave-one-out differences
\begin{align}
  \label{EqnDefnLooDiff}
\LooDiff & \defn \EstPar_\numobs - \LOO{\EstPar}{\numobs-1}{i} \qquad
\mbox{for $i \in [\numobs] \defn \{1, \ldots, \numobs \}$.}
\end{align}

\begin{lemma}
  \label{LemJackknifeKey}
Suppose that conditions~\ref{ass:tail} and~\ref{ass:convergence} are
in force, and the sample size $\Numobs$ satisfies the lower
bound~\eqref{eq:linear_approx_Z_est_sample_ahead}. Introduce the
shorthand $\dimdelta\defn\ParDim+\log(1/\delta)$. Then each of the
following statements holds with probability at least $1 - \delta$:
\begin{subequations}
\begin{align}
\label{eq:LemJackknifeKey_eq1}    
 \enorm{\LooDiff} \; \leq \; \frac{c}{\strongconvex \Numobs}
 \enorm{\ZFun(\Sam_i, \EstPar_\numobs)} & \; \leq \; \polyshort \frac{
   \sqrt{\dimdelta} \log^2(\Numobs/\delta)}{\Numobs} \qquad \mbox{for
   all $i \in [\numobs]$,} \\
\label{eq:LemJackknifeKey_eq3}
 \bigenorm{\LooDiff + \tfrac{\GradFun_{\TruePar}^{-1} \ZFun(\Sam_i,
     \TruePar)}{\Numobs-1}} & \; \leq \; \polyshort
 \log^3(\Numobs/\delta) \frac{\dimdelta}{\Numobs^{3/2}} \qquad
 \mbox{for all $i \in [\numobs]$, and} \\
\label{eq:LemJackknifeKey_eq4}
\enorm{\sum_{i=1}^\Numobs \LooDiff } & \; \leq \; \polyshort
\log^4(\Numobs/\delta) \frac{{\dimdelta}}{\Numobs} \qquad \mbox{for
  all $i \in [\numobs]$.}
\end{align}
\end{subequations}
\end{lemma}
\noindent We refer the reader to~\Cref{sec:proof_LemJackknifeKey} for
the proof.  Note that the functional $\TargetFun$ plays no role in
this lemma, so that it applies more generally to $Z$-estimators, and
so may be of independent interest.


\subsection{High-dimensional asymptotic normality}
\label{subsec:asymptotic}

Thus far, we have established non-asymptotic bounds on the
approximation of the plug-in and jackknife errors in term of the
stochastic noise component $\StochNoise_\numobs$, and for the plug-in
estimator, the bias term $\highbias$.  These bounds allow us to study
the asymptotics of these estimators under a form of high-dimensional
scaling, in which we allow the problem dimension $\usedim$, as well as
all the related parameters $(L, \sigma, \MyBou, 1/\gamma)$ to grow as
a function of $\numobs$.  Although all of these parameters depend on
$\numobs$, we elect to omit this dependence so as to keep the notation
streamlined.

So as to make our distributional statements clear, we do retain the
potential $\numobs$-dependence of the variance $\TrueVar^2$ defined in
equation~\eqref{EqnDefnVarTrue}.  We provide distributional statements
about the rescaled error $\sqrt{\numobs} \frac{(\tauhat -
  \tau(\thetastar))}{\TrueVar}$, where $\tauhat$ is either the plug-in
or the jackknife-corrected estimate. \\

Given that we are allowing problem parameters to scale with sample
size, we require certain regularity conditions on the problem
sequence:
\paragraph{Regular sequences:}  We say that problem sequence
is \emph{regular} if: (a) \Cref{ThmPlugin,ThmJackknife} hold with a
function $\polyprime$ such that
\begin{subequations}
\begin{align}
\label{EqnRegularScaling}  
\lim_{\numobs \rightarrow +\infty} \polyprime \; \numobs^{ -t} = 0
\qquad \mbox{for any $t > 0$,}
\end{align}
and the variance sequence $\VarTrue^2 = \var(\StochNoise_\numobs)$
satisfies
\begin{align}
\label{EqnVarScale}  
 \lim \inf_{\numobs \rightarrow + \infty} \VarTrue^2 > 0, \quad
 \mbox{and} \quad \lim \sup_{\numobs \rightarrow + \infty} \VarTrue^2
 < + \infty.
\end{align}
\end{subequations}
These conditions are mild, and serve to eliminate pathological
instances.\\

In stating our results, we make use of the shorthand $a_n/b_n
\highdimto 0$ to mean that there exists some $\delta > 0$ such that
$a_n^{1 + \delta}/b_n \rightarrow 0$.

\begin{corollary}[High-dimensional asymptotics for plug-in]
\label{CorHDAsymptoticPlugIn}
Under the conditions of~\Cref{ThmPlugin}, for any regular sequence of
problems, the plug-in estimator $\TargetFun(\EstPar_\numobs)$ has the
following properties:
\begin{itemize}
\item[(i)] {\underline{Asymptotic normality:}} Its error satisfies
\begin{subequations}
\begin{align}
 \sqrt{\numobs} \frac{(\TargetFun(\EstPar_\numobs) -
   \TargetFun(\TruePar))}{\TrueVar} \convdist \Normal(0,1) \qquad
 \mbox{whenever} \quad \frac{\ParDim^2}{\numobs} \highdimto 0.
\end{align}
\item[(ii)] {{\underline{Failure:}}} Suppose that $\ParDim^2/\numobs
  \rightarrow \infty$ while $\ParDim^{3/2}/\numobs \highdimto 0$ and
  $\highbias/\ParDim \not \longrightarrow 0$. Then asymptotic
  normality fails to hold---viz.:
\begin{align}
 \sqrt{\numobs} \frac{ (\TargetFun(\EstPar_\Numobs) -
   \TargetFun(\TruePar))}{\VarTrue} \not \convdist \Normal(0,1).
\end{align}
\end{subequations}
\end{itemize}
\end{corollary}
\noindent See~\Cref{SecProofCors} for the proof. \\

In brief, \Cref{CorHDAsymptoticPlugIn} asserts that, under certain
regularity conditions, the scaling of the ratio $\ParDim^2/\numobs$
specifies the \emph{boundary of asymptotic normality} for the plug-in
estimator.  Of particular interest is that it is guaranteed to fail
when $\usedim^2/\numobs \rightarrow \infty$ and the bias term
$\highbias$ is \emph{non-vanishing} in the sense that
$\highbias/\usedim$ does \emph{not} converge to zero.  This result
generalizes previous results by Portnoy~\cite{portnoy1988asymptotic},
who analyzed the special case of exponential families. \\

We now turn to an analogous high-dimensional guarantee for the
jackknife-corrected estimate $\JackEst$:
\begin{corollary}[High-dimensional asymptotics for $\JackEst$]
\label{CorHDAsymptoticDebias}
Under the conditions of~\Cref{ThmJackknife}, for any regular sequence
of problems, the jackknife estimator satisfies
\begin{align}
 \sqrt{\numobs} \frac{(\JackEst - \TargetFun(\TruePar))}{\TrueVar}
 \convdist \Normal(0,1) \qquad \mbox{whenever} \quad
 \frac{\ParDim^{3/2}}{\numobs} \highdimto 0.
\end{align}
\end{corollary}
\noindent See~\Cref{SecProofCors} for the proof.\\

\Cref{CorHDAsymptoticDebias} provides a crisp summary of the benefits
of the jackknife estimator relative to the plug-in: asymptotic
normality is retained under the milder requirement
$\usedim^{3/2}/\numobs \highdimto 0$, a setting in which the plug-in
estimator typically fails to be asymptotically normal.


\subsection{Jackknife variance estimation and confidence intervals}
\label{SecCI}

Observe that the distributional statements in
Corollaries~\ref{CorHDAsymptoticPlugIn}
and~\ref{CorHDAsymptoticDebias} involve the true variance $\TrueVar^2
= \var(\StochNoise_\numobs)$, which is an unknown quantity.  In order
to obtain distributional guarantees that are useful in practice, we
need to replace this population quantity with a consistent estimate.

The jackknife can again be used to provide such an estimate: in
particular, the jackknife estimate of
variance~\cite{efron1982jackknife,shao2012jackknife} is given by
\begin{align}
\label{EqnJackVarEst}
\JackVar^2 \defn \frac{\numobs-1}{\numobs} \sumn
\left[\TargetFun(\LOO{\EstPar}{\numobs-1}{i}) - \frac{ \sumn
    \TargetFun(\LOO{\EstPar}{\numobs-1}{i})}{\numobs} \right]^2.
\end{align}
Classical results provide consistency guarantees for this estimate
with the dimension fixed; here we establish a stronger guarantee that
holds with high-dimensional scaling.

\begin{proposition}
\label{PropJackVarEst}  
Under the conditions of~\Cref{CorHDAsymptoticDebias} (so that
$\usedim^{3/2}/\numobs \highdimto 0$), the jackknife variance estimate
$\JackVar$ is consistent---viz.  \mbox{$\big(\sqrt{\numobs} \JackVar -
  \VarTrue \big) \convprob 0$.}
\end{proposition}
\noindent The proof, given in~\Cref{SecProofPropJackVarEst}, makes
essential use of the previously stated~\Cref{LemJackknifeKey}
from~\Cref{SecQuantJack}, which provides non-asymptotic bounds on the
leave-one-out differences~\eqref{EqnDefnLooDiff}, and functions
thereof.  Tight control of such quantities is required so as to
guarantee consistency of $\JackVar$ under the scaling
$\usedim^{3/2}/\numobs \highdimto 0$, so that we can make inferential
statements about the jackknife estimator in the regime where it
strictly dominates the plug-in. \\

It is interesting to wonder to what extent the requirement
$\usedim^{3/2}/\numobs \highdimto 0$ could be loosened while still
guaranteeing consistency of the jackknife variance estimate.  As noted
earlier, even for the linear model, jackknife variance estimation is
known to be inconsistent if $\usedim/\numobs$ converges to a non-zero
constant~\cite{karoui2016can}, so that there is a small window for
improvement.

\paragraph{Confidence intervals:}  Let us conclude by summarizing
some inferential consequences of our corollaries and
\Cref{PropJackVarEst}.  These results, in conjunction with Slutsky's
lemma, ensure that
\begin{subequations}
\begin{align}
\frac{\TargetFun(\EstPar_\numobs) - \TargetFun(\TruePar)}{\JackVar}
\convdist \Normal(0, 1) \qquad \mbox{when $\usedim^2/\numobs
  \highdimto 0$ and}, \\ \mbox{and} \quad \frac{\JackEst -
  \TargetFun(\TruePar)}{\JackVar} \convdist \Normal(0, 1) \qquad
\mbox{when $\usedim^{3/2}/\numobs \highdimto 0$.}
\end{align}
\end{subequations}
These distributional statements can be used to generate confidence
intervals (CIs) with asymptotic guarantees on their levels in the
standard way, and we explore the properties of these CIs in our
numerical experiments given in~\Cref{sec:simulation}.


\section{Numerical studies}
\label{sec:simulation}

In this section, we present the results of some numerical studies
designed to assess the finite-sample performance of the jackknife
estimator $\JackEst$ relative to the plug-in estimator
$\TargetFun(\EstPar_\numobs)$.  In addition, depending on the problem
class under study, we also compare $\JackEst$ to \emph{non-black-box}
estimators that are specifically tailored to the problem structure
(e.g., Jeffrey's penalty method~\cite{firth1993bias} for MLEs).  In
general, we expect that such tailored procedures might outperform a
generic procedure such as the jackknife in some settings, and our goal
is to characterize any such differences.

Given our interest in finite-sample and high-dimensional results, we
consider two different regimes of sample size in our studies:
\begin{enumerate}
\item[(A)] Fixed sample size $\numobs$: Fixing the sample size
  $\numobs$, we vary the problem dimension $\ParDim=[\numobs^{r}]$ for
  a range of exponents $r \in [0.1, 0.9]$.
\item[(B)] Varying sample size $\numobs$: We vary the sample size
  $\numobs$, with the dimension scaling as $\usedim = [\numobs^{r}]$
  with fixed exponent $r = 2/3$.
\end{enumerate}
For the fixed-$\numobs$ scenario, the problem becomes harder as the
exponent $r$ increases. Moreover, based on the predictions of
\Cref{ThmPlugin,ThmJackknife}, the performance of the
jackknife-corrected estimate $\JackEst$ should remain valid for all
exponents $r < 2/3$, while the performance of the plug-in estimator
$\TargetFun(\EstPar_\numobs)$ might deteriorate for exponents $r >
1/2$.  Similarly, for the varying-$\numobs$ scenario, the error of the
plug-in estimator $\TargetFun(\EstPar_\numobs)$ might increase as
$\numobs$ increases, while the error of $\JackEst$ should remain at a
constant level.

For all simulations, we estimate the performance by averages over
$\numtrial = 1000$ Monte Carlo trials. For each trial $\tind = 1,
\ldots, \numtrial$, we generate $\numobs$ samples $\{\Sam_{i,
  \tind}\}_{i=1}^n$ and compute the estimates
$\tauhat_{\numobs,\tind}$ and their variance estimates
$\stdest^2_{\numobs,\tind}$. We estimate the variance of $\JackEst$
using the jackknife variance estimator, and the variance of
$\TargetFun(\EstPar_\numobs)$ via a plug-in estimate of the variance
$\VarTrue^2$ as defined in equation~\eqref{EqnDefnVarTrue}.  To
evaluate the performance of estimators, we compute the empirical bias
and empirical mean squared error (MSE)
\begin{align*}
\ebias & \defn | \frac{1}{\numtrial} \sum_{\tind=1}^{\numtrial}
\tauhat_{\numobs,\tind} - \taustar| \quad \mbox{and} \quad \emse \defn
\frac{1}{\numtrial} \sum_{\tind=1}^\numtrial (\tauhat_{\numobs,\tind}
- \taustar)^2,
\end{align*}
where $\taustar \defn \TargetFun(\TruePar)$ is the target value.

We also evaluate the inferential guarantee of estimators in terms of
their empirical $95\%$ coverage probability. Namely, for each point
and variance estimate, we use the normal approximation to the
$t$-statistic $(\tauhat_{\numobs,\tind} -
\taustar)/\stdest_{\numobs,\tind}$ with $\Normal(0,1)$ and construct
the approximate $95 \%$ confidence interval $[\tauhat_{\numobs,\tind}
  -1.96 \stdest_{\numobs, \tind},\tauhat_{\numobs, \tind} +1.96
  \stdest_{\numobs, \tind}]$ for $\taustar$. The empirical coverage
probability is defined as the portion of trials in which $\taustar$
falls into the confidence interval. We also report the average $95\%$
confidence interval length as $\frac{1}{\numobs}
\sum_{\tind=1}^\numtrial 3.92 \widehat\sigma_{\numobs,\tind}$.

In order to facilitate comparisons with the statement of
Theorem~\ref{ThmPlugin}---in particular, see
equation~\eqref{EqnPlugin}---we scale the bias and MSE, respectively,
by the pre-factors $\sqrt{\numobs}$ and $\numobs$ for the
varying-$\numobs$ regime.  As $\numobs$ increases,
equation~\eqref{EqnPlugin} shows that the bias of
$\TargetFun(\EstPar_\numobs)$ increases at a rate of $\order{\ParDim /
  \sqrt{\Numobs}}$, while the bias of $\JackEst$ remains
constant. Similarly, we expect the MSE of
$\TargetFun(\EstPar_\numobs)$ to grow, whereas the MSE of $\JackEst$
should remain stable.


\subsection{Quadratic functionals in linear regression}
\label{SecQuadSim}

We begin by considering the problem of estimating the value of a
quadratic functional $\TargetFun(\Par) = \Par^{\mytrans} \Amat \Par$
in the setting of linear regression; recall our previous discussion
from~\Cref{ExaQuadLinear} for motivation and set-up,
and~\Cref{FigQuadSims} for some motivating empirical results. For $i =
1, \ldots, \numobs$, we generate a standard Gaussian covariate vector
$\State_{i} \stackrel{i.i.d}{\sim} \Normal(0,\IdMat)$ and response
variable variable $y_i = \State_i^{\mytrans} \TruePar +
\varepsilon_i$, where $\varepsilon_{i}\stackrel{i.i.d}{\sim}
\Normal(0, \sigma^2)$ with $\sigma = 1$.  With the unknown parameter
vector $\TruePar \defn (1, \ldots, 1)/ \sqrt{\ParDim}$, we consider
the problem of estimating the quadratic functional $\TargetFun(\Par) =
\|\Par\|_2^2$, corresponding to a particular case with $\Amat \defn
\IdMat$. Observe that $\TargetFun(\TruePar) = 1$ by construction.

We evaluate three popular methods to estimate
$\TargetFun(\TruePar)$. The first method is the plug-in method, in
which the OLS estimate $\EstPar_\numobs$ of $\TruePar$ is used as a
plug-in to the quadratic form. The second method is the jackknife
estimator $\JackEst$.  Third, we compare to a specialized method from
Kline et al.~\cite{kline2020leave}; it is tailored to the problem at
hand, using a particular form of debiasing that exploits knowledge of
the model.  More precisely, define the variance estimator
$\widehat\sigma_i^2 = y_i(y_i - \covar_i^{\mytrans}
\LOO{\EstPar}{\numobs-1}{i})$ and the sample covariance matrix
$\Smat_{xx} \defn \frac{1}{\numobs} \sumn \covar_i
\covar_i^{\mytrans}$.  Kline et al. study the estimate
\begin{align*}
\unbiased & \defn \EstPar_\numobs^{\mytrans} \Amat \EstPar_\numobs -
\frac{1}{\numobs^2}\sumn \covar_i^{\mytrans} \Smat_{xx}^{-1} \Amat
\Smat_{xx}^{-1} \covar_i \widehat\sigma_i^2,
\end{align*}
along with the variance estimate\footnote{We apply with maximum with
$0$, since the estimate $\widehat\sigma_i^2$ can take negative
values.} given by \mbox{$\stdest^2_{\unbiasedname} \defn \max\left\{
  \frac{1}{\numobs} \sumn ( \nabla\TargetFun(\EstPar_\numobs)
  \Smat^{-1}_{xx}\covar_i)^2 \widehat\sigma_i^2, \; 0 \right \}$.}

\begin{figure}[ht!]
\begin{tabular}{cccc}
  \widgraph{0.22\textwidth}{\newfigdir/quad_bias_n400_p127_seed6666_rep_1000_info_fixed-n_runnum_1}
  &
  \widgraph{0.224\textwidth}{\newfigdir/quad_mse_n400_p127_seed6666_rep_1000_info_fixed-n_runnum_1}
  &
  \widgraph{0.219\textwidth}{\newfigdir/quad_coverage_n400_p127_seed6666_rep_1000_info_fixed-n_runnum_1}
  &
  \widgraph{0.224\textwidth}{\newfigdir/quad_length_n400_p127_seed6666_rep_1000_info_fixed-n_runnum_1}
\end{tabular}
\caption{Plots of the bias, MSE, coverage, and coverage length for
  three different estimators $\TargetFun(\EstPar_\numobs)$,
  $\JackEst$, and $\unbiased$ in~\Cref{SecQuadSim} for linear
  regression.  With sample size $\numobs = 400$, we set the dimension
  $\ParDim = [\numobs^{r}]$ for exponents $r \in [0.1, 0.9]$, and show
  plots with the exponent $r$ on the horizontal axis.  The points (on
  the curves in each plot) are obtained by taking a Monte Carlo
  average over $\numtrial=1000$ independent trials; error bars denote
  $\pm1$ of the standard deviation.  }
\label{fig:first-simulation}
\end{figure}

For the fixed-$\numobs$ case, we set the sample size $\numobs = 400$
and vary the problem dimension $\ParDim$.~\Cref{fig:first-simulation}
shows that when $\numobs \gtrsim \ParDim^2$---or equivalently,
$\ParDim \lesssim \numobs^{1/2}$---all three estimators
$\TargetFun(\EstPar_\numobs)$, $\JackEst$ and $\unbiased$ have similar
MSE.  However, when $\ParDim$ grows and exceeds $\numobs^{1/2}$, both
the bias and the MSE of the plug-in estimate
$\TargetFun(\EstPar_\numobs)$ increase rapidly.  Both the Kline et
al. estimate $\unbiased$ and jackknife debiasing mitigate this issue
in the regime $\ParDim \lesssim \numobs^{2/3}$; the jackknife estimate
exhibits some increase in bias/MSE in the regime $\ParDim \gtrsim
\numobs^{2/3}$.  Similar observations apply to the behavior of the
coverage probability: all three estimators have empirical coverage
probabilities close to the target $95\%$ when $\ParDim \lesssim
\numobs^{1/2}$. The plug-in method exhibits undercoverage in the
regime $\ParDim \geq \numobs^{1/2}$, whereas both the jackknife and
Kline method preserve coverage in the regime $\ParDim \lesssim
\numobs^{2/3}$ for the jackknife.  Again, the Kline approach has
better performance than the jackknife in the regime $\ParDim \gtrsim
\numobs^{2/3}$, which is out of the regime for which our jackknife
theory agrees.  Overall, the simulation results agree well
with~\Cref{CorHDAsymptoticPlugIn,CorHDAsymptoticDebias}, which
predicts that $\TargetFun(\EstPar_\numobs)$ has bad coverage or
unreasonable coverage length after $\ParDim \ge \numobs^{1/2}$, while
$\JackEst$ has valid coverage even when $\numobs^{1/2} \leq \ParDim
\lesssim \numobs^{2/3}$.

For the case of varying sample size $\numobs$ (in which the pair
$(\numobs,\ParDim)$ scales as $\ParDim = [\Numobs^{2/3}]$),
Figure~\ref{fig:second-simulation-growing} shows that both the bias
and MSE of $\TargetFun(\EstPar_\numobs)$ increase, and the coverage of
$\TargetFun(\EstPar_\numobs)$ drops significantly as $\numobs$
grows. In contrast, the bias, MSE and coverage of $\JackEst$ remains
stable as $\numobs$ grows. Again, these observations are consistent
with our theoretical predictions.

\newcommand{\quadwidth}{0.22\textwidth}
\begin{figure}[ht!]
  \begin{tabular}{cccc}
\widgraph{0.215\textwidth}{\newfigdir/quad_bias_n10240_p471_seed6666_rep_1000_info_vary-n_runnum_1}
&
\widgraph{\quadwidth}{\newfigdir/quad_mse_n10240_p471_seed6666_rep_1000_info_vary-n_runnum_1}
&
\widgraph{\quadwidth}{\newfigdir/quad_coverage_n10240_p471_seed6666_rep_1000_info_vary-n_runnum_1}
&
\widgraph{\quadwidth}{\newfigdir/quad_length_n10240_p471_seed6666_rep_1000_info_vary-n_runnum_1}
\\
(a) & (b) & (c) & (d)
\end{tabular}
\caption{Plots of the bias, MSE, coverage, and coverage length for
  three different estimators $\TargetFun(\EstPar_\numobs)$,
  $\JackEst$, $\unbiased$ in~\Cref{SecQuadSim} when $\ParDim =
  [\Numobs^{2/3}]$ and $n = 2^s \times 320$ for $s = 0,\ldots, 5$. The
  plots are obtained by averaging over $\numtrial=1000$ trials and the
  error bars denote the $\pm1$ standard deviation. }
\label{fig:second-simulation-growing}
\end{figure}


\subsection{Logistic regression}
\label{SecLogSim}

We now turn to the problem of logistic regression, for which we
evaluate the performance of the Maximum Likelihood Estimator (MLE),
the jackknife estimator derived from the MLE, and Jeffrey's penalty
method~\cite{firth1993bias}, which maximizes a penalized form of
likelihood tailored to the particular problem (i.e., it is a
non-black-box approach).

For logistic regression, we generate $\numobs$ i.i.d. samples $(X_i,
y_i)_{i=1}^\numobs$, where $X \in \real^{\usedim -1}$ is a covariate
vector, and the conditional distribution of the response $y_i \in
\{0,1\}$ given $X = X_i$ is given by
\begin{align*}
\mathbb{P}(y_i = 1 \mid X = X_i) = \frac{\exp(\alphastar +
  X_i^{\mytrans} \betastar)}{1 + \exp(\alphastar + X_i^{\mytrans}
  \betastar)}.
\end{align*}
where $\alphastar \in \real$ and $\betastar \in \real^{\usedim-1}$ are
unknown parameters.  We draw i.i.d. covariate vectors $X_i \sim
\mathcal{N}(0, \IdMat_{\usedim -1})$ and generate the responses
$\{y_i\}_{i=1}^\numobs$ from the model with $\alphastar = 1$ and
\mbox{$\betastar = (1/ \sqrt{\usedim}, \ldots, 1/ \sqrt{\usedim}) \in
  \R^{\ParDim - 1}$.}  We consider the problem of estimating the
linear functional $\TargetFun(\alphastar, \betastar) \defn
\alphastar$.

In the fixed-$\numobs$ scenario, similar
to~\Cref{SecQuadSim,SecWeakIVSim}, \Cref{FigLogSimGrowd} illustrates
that the performance of $\TargetFun(\EstPar)$ deteriorates when
$\ParDim \geq \numobs^{1/2}$, while $\JackEst$ remains valid when
$\ParDim \lesssim \numobs^{2/3}$, maintaining a coverage probability
close to 95\%.  Also, we observe in~\Cref{FigLogSimGrown} that the
bias and MSE increases while the coverage probability decreases for
the plug-in estimator, when $(\numobs,\ParDim)$ grow with
$\ParDim=[\numobs^{2/3}]$.  We remark that when $\usedim \gtrsim
\numobs^{2/3}$, which is out of the scope of our theoretical
discussions, then $\JackEst$ becomes more unstable than other two
estimators.
\newcommand{\logwidth}{0.22\textwidth}
\begin{figure}[htp]
\begin{tabular}{cccc}  
\widgraph{\logwidth}{\newfigdir/logistic_bias_n400_p127_alpha1.0_seed6666_rep_1000_info_fixed-n_runnum_1}
&
\widgraph{\logwidth}{\newfigdir/logistic_mse_n400_p127_alpha1.0_seed6666_rep_1000_info_fixed-n_runnum_1}
&
\widgraph{\logwidth}{\newfigdir/logistic_coverage_n400_p127_alpha1.0_seed6666_rep_1000_info_fixed-n_runnum_1}
&
\widgraph{\logwidth}{\newfigdir/logistic_length_n400_p127_alpha1.0_seed6666_rep_1000_info_fixed-n_runnum_1}
\\
(a) & (b) & (c) & (d)
\end{tabular}
\caption{Plots of the bias, MSE, coverage, and coverage length for
  $\TargetFun(\EstPar_\numobs)$, $\JackEst$ and
  $\tauhat_{\text{jeffreys}}$ in~\Cref{SecLogSim} for logistic
  regression.  With fixed sample size $\numobs=400$, we set problem
  dimension $\ParDim=[\numobs^{r}]$ for exponents $r \in [0.1, 0.9]$.
  The results are obtained by averaging over $\numtrial=1000$ trials
  and the error bars denote the $\pm1$ standard deviation. }
\label{FigLogSimGrowd}
\end{figure}

\newcommand{\logtwowidth}{0.22\textwidth}
\begin{figure}[ht!]
  \begin{center}
\begin{tabular}{cccc}    
\widgraph{0.215\textwidth}{\newfigdir/logistic_bias_n10240_p471_alpha1.0_seed6666_rep_1000_info_vary-n_runnum_1}&
\widgraph{\logtwowidth}{\newfigdir/logistic_mse_n10240_p471_alpha1.0_seed6666_rep_1000_info_vary-n_runnum_1}&
\widgraph{\logtwowidth}{\newfigdir/logistic_coverage_n10240_p471_alpha1.0_seed6666_rep_1000_info_vary-n_runnum_1}&
\widgraph{\logtwowidth}{\newfigdir/logistic_length_n10240_p471_alpha1.0_seed6666_rep_1000_info_vary-n_runnum_1}\\
(a) & (b) & (c) & (d)
\end{tabular}
\caption{Plots of the bias, MSE, coverage, and coverage length for for
  $\TargetFun(\EstPar_\numobs)$, $\JackEst$ and
  $\tauhat_{\text{jeffreys}}$ in~\Cref{SecLogSim} when $\ParDim =
  [\Numobs^{2/3}]$ and $n = 2^s \times 320$ for $s = 0,\ldots, 5$. The
  plots are obtained by averaging over $\numtrial = 1000$ trials and
  the error bars denote the $\pm1$ standard deviation.}
\label{FigLogSimGrown}
  \end{center}
\end{figure}


\subsection{IV estimators with growing number of instruments} 
\label{SecWeakIVSim}

Recall from the class of instrumental variable (IV) methods from
Example~\ref{ExamIV}, and in particular the defining
equations~\eqref{EqnOverOne} and~\eqref{EqnOverTwo}.  Angrist et
al.~\cite{angrist1999jackknife} proposed two debiasing procedures
tailored to the two-stage least squares (TSLS)
estimate~\eqref{EqnOverThree}, referred to as JIVE1 and JIVE2
respectively.  In this simulation study, we compare these
procedures---which are tailored specifically to IV problem---to our
black-box form of jackknife correction.

We performed simulations by generating $\numobs$ i.i.d. samples $(W_i,
X_i, Y_i)$ from an ensemble of the following type.  We choose each
$W_i$ from a uniform distribution over $\ivdim$ discrete outcomes (so
as to mimic the use of various types of indicator variables in
instrumental analysis).  We then generated each covariate $X_i$ using
equation~\eqref{EqnOverTwo} with the parameter vector \mbox{$\pistar
  \defn \begin{bmatrix} 0 & \frac{2}{\ivdim} & \frac{4}{\ivdim} &
    \ldots & \frac{2(\ivdim-1)}{\ivdim}
\end{bmatrix} \in \real^\ivdim$.}
We generated the response $Y_i$ via equation~\eqref{EqnOverOne} with
parameter $(\alphastar, \betastar = (0,1)$.  As for the noise terms
$(\varepsilon_i, \eta_i)$ in equations~\eqref{EqnOverOne}
and~\eqref{EqnOverTwo} respectively, we sampled them from a zero-mean
Gaussian distribution with equal variances \mbox{$\var(\varepsilon_i)
  = \var(\eta_i) = 0.25$} and covariance $\cov(\epsilon_i, \eta_i) =
0.2$. We focused on the problem of estimating the component
$\betastar$ of the full parameter vector $\thetastar = (\alphastar,
\betastar, \pistar) \in \real^{\ivdim + 2}$.

For ensembles of this type, our goal was to study the dependence on
the sample size $\numobs$ and instrument dimension $\usedim \defn
\ivdim + 2$.  We implemented the TSLS estimate~\eqref{EqnOverThree},
and compared it to the JIVE1 and JIVE2 estimators, as well as our
jackknife-corrected version of the TSLS estimate.  We estimate the
variance of TSLS, JIVE1, JIVE2 using the bootstrap with $T = 500$
bootstrap samples, and the variance of the jackknife using the
variance estimate $\JackVar^2 $ from equation~\eqref{EqnJackVarEst}.

\begin{figure}[htp]
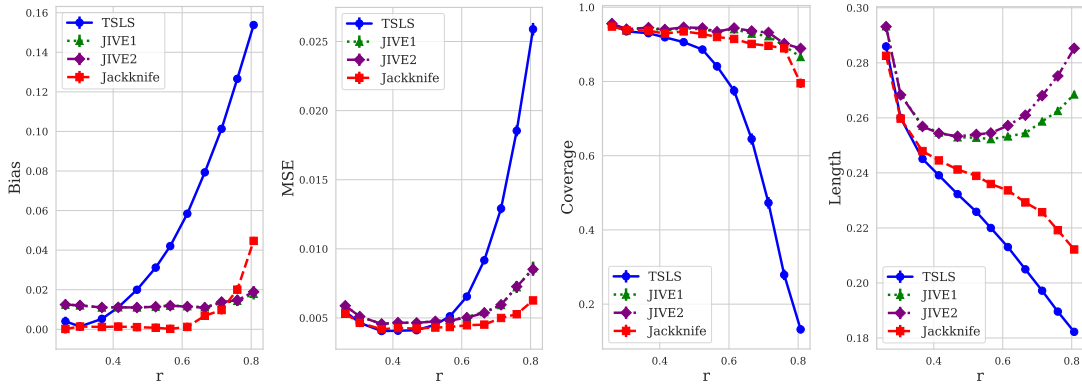

\widgraph{0.225\textwidth}{\newfigdir/iv_bias_n200_p72_seed8888_rep_1000_info_fixed-n_runnum_1}
\widgraph{0.23\textwidth}{\newfigdir/iv_mse_n200_p72_seed8888_rep_1000_info_fixed-n_runnum_1}
\widgraph{0.22\textwidth}{\newfigdir/iv_coverage_n200_p72_seed8888_rep_1000_info_fixed-n_runnum_1}
\widgraph{0.225\textwidth}{\newfigdir/iv_length_n200_p72_seed8888_rep_1000_info_fixed-n_runnum_1}
\caption{Plots of the bias, MSE, coverage, and coverage length for
  TSLS, the jackknife, JIVE1 and JIVE2 in~\Cref{SecWeakIVSim} for
  instrumental variables estimation.  We choose $\numobs=200$ and let
  $\ParDim=[\numobs^{r}]$ for exponents $r \in [0.2, 0.9]$. The
  results are obtained by averaging over $\numtrial=1000$ trials and
  the error bars denote the $\pm1$ standard deviation.}
\label{FigWeakIVSimGrowd} 
\end{figure}

First, we consider the fixed sample size scenario with $\numobs = 200$
and the dimension $\ParDim$ growing.  As shown
in~\Cref{FigWeakIVSimGrowd}, all four estimators have similar
mean-squared error (MSE) when $\ParDim \lesssim \numobs^{1/2}$.  This
should be contrasted with the regime $\ParDim \geq \numobs^{1/2}$,
where the bias and MSE of the plug-in estimate based on TSLS are much
larger than those of $\tauhat_{\text{JIVE1}}$,
$\tauhat_{\text{JIVE2}}$ and $\JackEst$. As for the coverage
probabilities, they are close to $95\%$ for all four estimators when
$\ParDim \lesssim \numobs^{1/2}$. In contrast, the coverage
probabilities of $\JackEst$, $\tauhat_{\text{JIVE1}}$ and
$\tauhat_{\text{JIVE2}}$ maintains when $\numobs^{1/2} \leq \ParDim
\lesssim \numobs^{2/3}$ while the coverage probability of
$\TargetFun(\EstPar_\numobs)$ decreases rapidly, as predicted by our
theory.  Moreover, when $\ParDim \geq \numobs^{0.7}$, although the
jackknife exhibits a larger bias than JIVE1 and JIVE2, it achieves
similar coverage and has a smaller MSE and confidence interval length
compared to JIVE1 and JIVE2.

\begin{figure}[ht!]
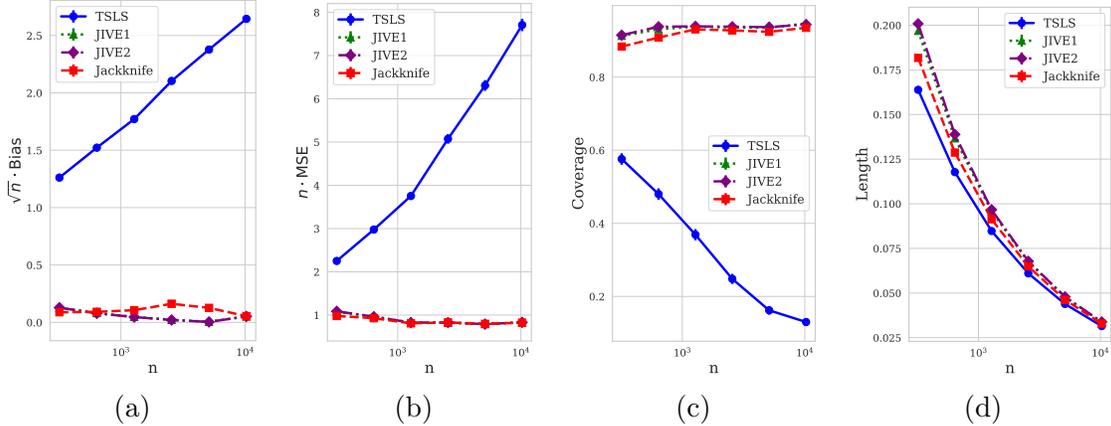

\begin{tabular}{cccc}
\widgraph{0.218\textwidth}{\newfigdir/iv_bias_n10240_p471_seed8888_rep_1000_info_vary-n_runnum_1}
&
\widgraph{0.208\textwidth}{\newfigdir/iv_mse_n10240_p471_seed8888_rep_1000_info_vary-n_runnum_1}
&
\widgraph{0.215\textwidth}{\newfigdir/iv_coverage_n10240_p471_seed8888_rep_1000_info_vary-n_runnum_1}
&
\widgraph{0.225\textwidth}{\newfigdir/iv_length_n10240_p471_seed8888_rep_1000_info_vary-n_runnum_1}
\\
(a) & (b) & (c) & (d)
\end{tabular}
\caption{Plots of the bias, MSE, coverage, and coverage length for
  TSLS, the jackknife, JIVE1 and JIVE2 in~\Cref{SecWeakIVSim} when
  $\ParDim = [\Numobs^{2/3}]$ and $n = 2^s \times 320$ for $s =
  0,\ldots, 5$. The plots are obtained by averaging over
  $\numtrial=1000$ trials and the error bars denote the $\pm1$
  standard deviation. }
\label{FigWeakIVSimGrown} 
\end{figure}

Next we compare the four estimators under the varied-$\numobs$
scenario where the pair $(\numobs,\ParDim)$ grow according to the
scaling $\ParDim = [\Numobs^{2/3}]$.  \Cref{FigWeakIVSimGrown} reveals
that the bias and MSE of $\TargetFun(\EstPar_\numobs)$ increase, and
the coverage of $\TargetFun(\EstPar_\numobs)$ declines significantly
as $\numobs$ increases, whereas the bias, MSE, and coverage of
$\JackEst$ remain approximately constant.  Similar
to~\Cref{SecQuadSim}, simulation results for instrumental variables
estimation agree with our theoretical predictions.  Additionally, the
jackknife estimator exhibits comparable coverage and achieves similar
or even lower MSE and confidence interval length than JIVE1 and JIVE2
across a broad range of $(\numobs,\ParDim)$ values.

\newcommand{\ivwidth}{0.21\textwidth}


\section{Proofs}
\label{sec:proofs}

We now turn to the proofs of our two theorems,
with~\Cref{SecProofThmPlugin,SecProofThmJackknife} devoted to the
proofs of~\Cref{ThmPlugin,ThmJackknife}, respectively.


\subsection{Proof of~\Cref{ThmPlugin}}
\label{SecProofThmPlugin}

We begin with a high-level overview of the proof structure, and also
highlighting some technical results that are of independent interest.
Recall that our theorem provides a non-asymptotic bound on the
difference $\PlugDiff \defn \Big( \TargetFun(\EstPar_\numobs) -
\TargetFun(\TruePar) \Big) - \frac{1}{\sqrt{\numobs}} \left
(\noiseW_{\numobs} + \frac{\highbias}{\sqrt{\numobs}} \right)$.  At
the core of our proof is a decomposition of this difference in terms
of a centered $U$-statistic, and four additional remainder
terms---namely, we write
\begin{align*}
\sqrt{\numobs} \PlugDiff = \Ustat + \sqrt{\numobs} \Big \{ \RemainL +
\RemainZ + \RemainTau \Big \},
\end{align*}
where the $U$-statistic $\Ustat$ is defined below in
equation~\eqref{EqnDefnUstat}, and we describe the three remainder
terms $\RemainL$, $\RemainZ$ and $\RemainTau$ in sequence.

Recalling the definition~\eqref{EqnDefnJmat} of $\Jmat$, the remainder
term $\RemainL$ depends on the \emph{first-order error} in the linear
approximation---namely,
\begin{align}
  \label{EqnDefnLinError}
\LinError_\numobs \defn \EstPar_\numobs - \TruePar - \frac{1}{\numobs}
\sumn \GradFun_{\TruePar}^{-1} \ZFun(\Sam_i, \TruePar)
\end{align}
obtained from a first-order Taylor series expansion of the
$Z$-estimating equations.  A central result in the proofs of both our
results is a non-asymptotic upper bound on $\|\LinError_\numobs\|_2$,
stated below as~\Cref{LemLinErrorBound}; see~\Cref{SecRemainL} for
this result, along with the definition~\eqref{EqnDefnRemainL} of
$\RemainL$.

The final terms to be controlled are:
the remainder terms $\RemainTau$ and $\RemainZ$, defined below in
equations~\eqref{EqnDefnRemainTau} and~\eqref{EqnDefnRemainZ}, which
arise from second-order terms in Taylor series expansions of
$\TargetFun$ and $\Zfun$, respectively.


\subsubsection{Controlling the $U$-statistic}
\label{SecUstat}
We begin by defining and controlling the centered $U$-statistic that
appears in our analysis. Recall the functions $\sionebase$ and
$\sitwobase$ from equation~\eqref{EqnSidefn}, and introduce the
empirical means
\begin{subequations}
\begin{align}
  \sioneave \defn \frac{1}{\numobs} \sumn \sionebase(\Sam_i) & =
  \frac{1}{\numobs} \sumn [\nabla \TargetFun(\TruePar) - \nabla_{\Par}
    \ZFun(\Sam_i,\TruePar)^{\mytrans} \MyIntSec{\TruePar}], \quad
  \mbox{and}\label{eq:def_phi} \\
\sitwoave \defn \frac{1}{\numobs} \sumn \sitwobase(\Sam_i) & =
\frac{1}{\numobs}\sumn \GradFun_{\TruePar}^{-1} \ZFun(\Sam_i,
\TruePar).\label{eq:def_psi}
\end{align}
Using this notation, we then define the $U$-statistic
\begin{align}
\label{EqnDefnUstat}  
\Ustat & \defn \sqrt{\numobs} \Big \{ \inprod{\sioneave}{ \sitwoave} +
\sitwoave^{\mytrans} \Mmat \sitwoave/2 \Big \} -
\frac{\highbias}{\sqrt{\numobs}},
\end{align}
\end{subequations}
where the matrix $\Mmat$ and bias term $\highbias$ were previously
defined in equations~\eqref{EqnDefnMmat}
and~\eqref{EqnDefnStochNoise}, respectively.  By construction, the
random variable $\Ustat$ is a $U$-statistic, and it can be verified to
have zero-mean.  It remains to establish a high-probability bound on
it, and the following result provides such a guarantee:
\begin{lemma}[Non-asymptotic control of $U$-statistic]
\label{LemUstatBound}
Under the conditions of~\Cref{ThmPlugin}, there is a pre-factor
$\polyshort = \poly$ such that
\begin{align}
\label{EqnUstatBound}
|\Ustat| & \leq \polyshort \sqrt{ \frac{\ParDim}{\numobs}}
\log^{5}(\numobs/\delta)
\end{align}
with probability at least $1 - \delta$.
\end{lemma}
\noindent See~\Cref{SecProofLemUstatBound} for the proof.


\subsubsection{Controlling the remainder term $\RemainL$}
\label{SecRemainL}
Next we define the remainder term
\begin{subequations}
\begin{align}
\label{EqnDefnRemainL}
\RemainL & \defn \underbrace{\inprod{\sioneave}{\Errfirst}}_{\revdefn
  \RemainLone} + \underbrace{\Errfirst^{\mytrans} \Mhat \Errfirst/2 +
  \sitwoave^{\mytrans} \Mhat \Errfirst}_{\revdefn \RemainLtwo} +
\underbrace{\sitwoave^{\mytrans} (\Mhat - \Mmat)
  \sitwoave/2}_{\RemainLthree},
\end{align}
where the first-order error $\LinError_\numobs$ was defined in
equation~\eqref{EqnDefnLinError}, and the random matrix $\Mhat$ is
given by $\Mhat \defn \nabla^2\TargetFun(\TruePar) -
\frac{1}{\Numobs}\sum_{i=1}^\Numobs \vtinprod{\IntSec_\TruePar}{
  \nabla^2_\Par \ZFun(\Sam_i, \TruePar)}$.  Note that by definition,
we have $\Exs[\Mhat] = \Mmat$, where the matrix $\Mmat$ was defined
previously~\eqref{EqnDefnMmat}.  The primary claim of this section is
that
\begin{align}
  \label{EqnRemainLBound}
|\RemainL| & \leq \polyshort \Big(\frac{\dimdelta}{\numobs}
\Big)^{3/2} \log^2(\numobs/\delta), \qquad \mbox{where $\dimdelta
  \defn \ParDim + \log(1/\delta)$,}
\end{align}
\end{subequations}
with probability at least $1 - \delta$. \\

In order to establish the bound~\eqref{EqnRemainLBound}, we begin with
a bound on the first-order error $\LinError_\numobs$.
\begin{lemma}[First-order approximation error]
\label{LemLinErrorBound}
Under the conditions of~\Cref{ThmPlugin} but with a relaxed sample
size condition
$\numobs/\log^3(\numobs/\delta)\geq\polyshortprime\cdot\Ccerr$, there
is a pre-factor $\polyshort = \poly$ such that
\begin{align}
\label{EqnLinErrorBound}  
  \| \LinError_\numobs\|_2 & \leq \polyshort \,
  \frac{\dimdelta}{\numobs} \log \numobs
\end{align}
with probability at least $1 - \delta$.
\end{lemma}
\noindent See~\Cref{SecProofLemLinErrorBound} for the proof.  We
suspect that this bound could be of more general utility in the
non-asymptotic analysis of $Z$-estimators. (Notably, it bounds the
full error vector $\thetahat - \thetastar$, as opposed to error
defined by the functional $\TargetFun$.) \\

The remaining ingredients in our analysis of $\RemainL$ are
high-probability bounds on the quantities $\|\sioneave\|_2$,
$\|\sitwoave\|_2$, and $\opnorm{\Mhat - \Mmat}$.  These terms can be
controlled via empirical process theory; in particular, we claim that
\begin{subequations}
\begin{align}
\label{EqnSiBound}  
\max \{ \|\sioneave\|_2, \|\sitwoave\|_2 \} & \leq \polyshort \Big \{
\sqrt{ \frac{\dimdelta}{\numobs}} + \frac{\dimdelta \log
  \numobs}{\numobs} \Big \} \qquad \mbox{and} \\
\label{EqnMBound}
\opnorm{\Mhat - \Mmat} & \leq \polyshort\Big( \sqrt{
  \frac{\Ccerr\log(\numobs/\delta)}{\Numobs}} + \frac{\dimdelta
  \log^{5/2} (\Numobs/\delta)}{\Numobs} \Big),
\end{align}
\end{subequations}
where each bound holds with probability at least $1 - \delta$.  Here,
as in the bound~\eqref{EqnLinErrorBound}, we state the bounds in terms
of a pre-factor $\polyshort = \poly$, whose value may change from
line-to-line.  We prove these bounds at the end of this section, using
auxiliary results on the suprema of empirical processes
from~\Cref{SecEmpProcess}.

Equipped with these auxiliary results, we are ready to prove the
claim~\eqref{EqnRemainLBound}.  Beginning with $\RemainLone$, we have
\begin{align}
|\RemainLone| \; \stackrel{(i)}{\leq} \; \|\sioneave\|_2
\|\LinError_\numobs\|_2 & \stackrel{(ii)}{\leq} \polyshort^2 \Big \{
\sqrt{\frac{\dimdelta}{\numobs}} + \frac{\dimdelta}{\numobs} \log
\numobs \Big \} \; \Big \{ \frac{\dimdelta}{\numobs} \log \numobs \Big
\} \notag \\
\label{EqnLoneBound}
& \stackrel{(iii)}{\leq} \polyshort' \Big(\frac{\dimdelta}{\numobs}
\Big)^{3/2} \log^2 \numobs,
\end{align}
where step (i) follows from the Cauchy--Schwarz inequality; step (ii)
follows by applying the bound on $\|\LinError_\numobs\|_2$ from
equation~\eqref{EqnLinErrorBound} and the bound~\eqref{EqnSiBound} on
$\|\sioneave\|_2$; and step (iii) follows since $\numobs \geq
\dimdelta$ by assumption.

Now turning to $\RemainLtwo$, we again apply the Cauchy--Schwarz
inequality along with the definition of the matrix operator norm,
thereby finding that
\begin{align*}
|\RemainLtwo| & \leq \tfrac{1}{2} \opnorm{\Mhat} \Big \{
\|\LinError_\numobs\|_2^2 + 2\|\sitwoave\|_2 \|\LinError_\numobs\|_2
\Big \} \notag \\
& \leq \tfrac{1}{2} \Big \{ \opnorm{\Mmat} + \opnorm{\Mhat - \Mmat}
\Big \} \; \Big \{ \|\LinError_\numobs\|_2^2 + 2\|\sitwoave\|_2
\|\LinError_\numobs\|_2 \Big \} \notag \\
& \leq \tfrac{1}{2} \Big \{ C(L, \strongconvex) + \opnorm{\Mhat -
  \Mmat} \Big \} \; \Big \{ \|\LinError_\numobs\|_2^2 +
2\|\sitwoave\|_2 \|\LinError_\numobs\|_2 \Big \},
\end{align*}
where the final step uses the fact that $\opnorm{\Mmat} \leq C(L,
\strongconvex) \defn L + \frac{L}{\strongconvex}$ from
condition~\ref{ass:tail}.

Under our sample size condition, the bound~\eqref{EqnMBound} implies
that $\opnorm{\Mhat - \Mmat} \leq \polyshort'$, where $\polyshort'$ is
another pre-factor.  Combining this fact with the
bounds~\eqref{EqnLinErrorBound} and~\eqref{EqnSiBound}, we see that
\begin{align}
|\RemainLtwo| & \leq \tilde{\polyshort} \Biggr \{
\Big(\frac{\dimdelta}{\numobs} \log(\numobs) \Big)^2 +
\Big(\frac{\dimdelta}{\numobs} \log(\numobs) \Big)
\Big(\sqrt{\frac{\dimdelta}{\numobs}} +
\frac{\dimdelta\log\numobs}{\numobs} \Big) \Biggr \} \notag \\
\label{EqnLtwoBound}
& \leq \tilde{\polyshort} \Big(\frac{\dimdelta}{\numobs} \Big)^{3/2}
\log^2(\numobs).
\end{align}
where the value of $\tilde{\polyshort}$ may change from line-to-line.

As for $\RemainLthree$, we have
\begin{align}
|\RemainLthree| & \leq \tfrac{1}{2} \opnorm{\Mhat - \Mmat} \; \|
\sitwoave\|_2^2 \leq \tilde{C} \Big \{
\sqrt{\frac{\dimdelta}{\numobs}} + \frac{\dimdelta \log
  \numobs}{\numobs} \Big \}^2 \cdot \Big( \sqrt{
  \frac{\Ccerr\log(\numobs/\delta)}{\Numobs}} + \frac{\dimdelta
  \log^{5/2} (\Numobs/\delta)}{\Numobs} \Big)\label{EqnthreeBound}
\end{align} for another pre-factor $\tilde{C}$,
where we have applied the bounds~\eqref{EqnSiBound}
and~\eqref{EqnMBound}.

Finally, combining the three
bounds~\eqref{EqnLoneBound},~\eqref{EqnLtwoBound},~\eqref{EqnthreeBound}
and the sample size condition yields the
claim~\eqref{EqnRemainLBound}. \\


\myunder{Proof of the bounds~\eqref{EqnSiBound}
  and~\eqref{EqnMBound}:} equation~\eqref{EqnSiBound} follows from the
first claim of Lemma~\ref{lm:op_norm_zfun} and the fact that
$\opnorm{\GradFun^{-1}},\vecnorm{\IntSec_\TruePar}{2}\leq \polyshort$
for some constant $\polyshort=\poly.$ Equation~\eqref{EqnMBound}
follows from the second claim of Lemma~\ref{lm:op_norm_zfun} with
$\Direc=\IntSec_\TruePar/\vecnorm{\IntSec_\TruePar}{2}$, combined with
the boundedness of $\vecnorm{\IntSec_\TruePar}{2}$.


\subsubsection{Controlling the Taylor series errors $\RemainTau$ and $\RemainZ$}
\label{SecRemainTau}

Using the shorthand $\DelHat \defn \ParEst_\numobs - \TruePar$, we
define the following two error terms
\begin{subequations}
\begin{align}
\label{EqnDefnRemainTau}    
\RemainTau & \defn \Err_\numobs^{\mytrans} \left \{ \int_0^1(1-s)[
  \nabla^2 \TargetFun \big(\TruePar + s \Err_\numobs \big) - \nabla^2
  \TargetFun(\TruePar)] \textup{d}s \right \} \Err_\numobs, \quad
\mbox{and} \\
\label{EqnDefnRemainZ}
\RemainZ & \defn - \Err_\numobs^{\mytrans} \left \{ \frac{1}{\Numobs}
\sum_{i=1}^\Numobs \left[ \int_{0}^1 (1-s)
  [\vtinprod{\IntSec_\TruePar}{ \nabla^2_\Par \ZFun \big(\Sam_i,
      \TruePar + s \Err_\numobs \big)- \nabla^2_\Par \ZFun(\Sam_i,
      \TruePar)}] \textup{d}s \right] \right \} \Err_\numobs.
\end{align}
\end{subequations}
As should be clear from their form, they arise from the integral error
representation in a Taylor series expansion.  With this notation, we
have:
\begin{lemma}[Bounds on third-order errors]
\label{LemRemainTau}
Under the conditions of~\Cref{ThmPlugin}, we have
\begin{subequations}
\begin{align}
\label{EqnRemainTauBound}      
|\RemainTau| & \leq \polyshort ( \frac{\myusedim}{\Numobs})^{3/2}
\log^2 \numobs, \quad \mbox{and} \quad \\
\label{EqnRemainZBound}
|\RemainZ| & \leq \polyshort \Big(\frac{\dimdelta}{\Numobs}\Big)^{3/2}
\log^{2} (\Numobs/\delta)
\end{align}
\end{subequations}
with probability at least $1 - \delta$.
\end{lemma}
\noindent See~\Cref{sec:proof_LemRemainTau} for the proof.


\subsubsection{Putting together the pieces}

Finally, we combine the pieces.  We begin with a formal statement of
the decomposition result that underlies our analysis:
\begin{lemma}
\label{LemPieces}
With the $U$-statistic from equation~\eqref{EqnDefnUstat} and
remainder terms $(\RemainL, \RemainTau, \RemainZ)$ from
equations~\eqref{EqnDefnRemainL},~\eqref{EqnDefnRemainTau}~and~\eqref{EqnDefnRemainZ}
respectively, we have
\begin{align}
\label{EqnPieces}  
\PlugDiff = \frac{\Ustat}{\sqrt{\numobs}} + \Big \{ \RemainL +
\RemainZ + \RemainTau \Big \}.
\end{align}
\end{lemma}
\noindent See~\Cref{SecProofLemPieces} for the proof of this claim.

Given the decomposition~\eqref{EqnPieces}, we can now complete the
proof of the bound~\eqref{EqnPlugin} stated in~\Cref{ThmPlugin}. In
particular, by the triangle inequality, we have
\begin{align*}
\sqrt{\numobs} |\PlugDiff| & \leq |\Ustat| + \sqrt{\numobs} \Big \{
|\RemainL| +|\RemainZ| + |\RemainTau| \Big \} \\
& \overset{(i)}{\leq}
\polyshort\Big(\sqrt{\frac{\ParDim}{\numobs}}\log^5(\numobs/\delta)+\frac{\dimdelta^{3/2}}{\Numobs}
\log^{2} (\Numobs/\delta)\Big) \leq \polyshort \,\Big\{ \sqrt{
  \frac{\ParDim}{\numobs}} + \frac{\ParDim^{3/2}}{\numobs} \Big\}
\log^{5}(\numobs/\delta)
\end{align*}
where inequality (i) follows by applying~\Cref{LemUstatBound} to bound
$|\Ustat|$; equation~\eqref{EqnRemainLBound} to bound $|\RemainL|$;
and~\Cref{LemRemainTau} to bound $|\RemainTau|$ and $|\RemainZ|$.


\subsection{Proof of~\Cref{ThmJackknife}}
\label{SecProofThmJackknife}

We now turn to the analysis of the error in the jackknife-corrected
estimate.  Throughout the proof, we use the shorthand $\DelHat =
\ParEst_\numobs - \TruePar$.  We begin outlining the proof strategy,
before stating the key lemmas and combining the pieces.


\subsubsection{Proof strategy}

We begin by deriving a useful expression for the error in the
jackknife-corrected estimator $\JackEst$.  From~\Cref{LemPieces}, we
have the decomposition
\begin{subequations}
\begin{align}
\label{EqnPieceTwo}  
\TargetFun(\EstPar_\numobs) - \TargetFun(\TruePar) -
\frac{1}{\sqrt{\numobs}} \noiseW_{\numobs} & = \underbrace{\Big \{
  \frac{\highbias}{\numobs} + \frac{\Ustat}{\sqrt{\numobs}} + \RemainL
  \Big \}}_{\revdefn \Term} + \RemainZ + \RemainTau.
\end{align}
From the definition~\eqref{eq:jacknife_est_def}, the jackknife
estimator modifies the plug-in estimator by subtracting a bias
estimate $\Debias_\numobs$; more precisely, we have
\begin{align}
\label{EqnPieceOne}
\JackEst = \TargetFun(\EstPar_\numobs) - \underbrace{
  \frac{\Numobs-1}{\Numobs}\sum_{i=1}^{\Numobs}
  (\TargetFun(\LOO{\EstPar}{\numobs-1}{i}) -
  \TargetFun(\EstPar_\numobs))}_{\revdefn \Debias_\numobs}.
\end{align}
Combining equations~\eqref{EqnPieceOne} and~\eqref{EqnPieceTwo}, we
find that
\begin{align}
\label{EqnJackDecomp}    
\JackEst - \TargetFun(\TruePar) - \frac{\noiseW_{\numobs}}{
  \sqrt{\numobs}} & = \Big \{ \Term - \Debias_\numobs \Big \} +
\RemainTau + \RemainZ.
\end{align}
\end{subequations}
Note that \Cref{LemRemainTau} from the proof of~\Cref{ThmPlugin} gives
high probability bounds on the terms $\RemainTau$ and $\RemainZ$.

Thus, the remaining challenge is to control the difference $\Term -
\Debias_\numobs$.  In~\Cref{SecProofLemPieces}, as part of
proving~\Cref{LemPieces}, we showed that $\Term = \Term_2 + \Term_3$,
where
\begin{align*}
\Term_2 & \defn \Big[\nabla \TargetFun(\TruePar)^{\mytrans} -
  \frac{1}{\Numobs}\sum_{i=1}^\Numobs \IntSec_\TruePar^{\mytrans}
  \nabla_\Par \ZFun(\Sam_i, \TruePar) \Big] \cdot \Err_\numobs, \quad
\mbox{and} \\
\Term_3 & \defn \frac{1}{2} \Err_\numobs^{\mytrans} \left \{
\nabla^2\TargetFun(\TruePar) - \frac{1}{\Numobs} \sum_{i=1}^\Numobs
\vtinprod{\IntSec_\TruePar}{ \nabla^2_\Par \ZFun(\Sam_i, \TruePar)}
\right \} \Err_\numobs.
\end{align*}
With this representation in hand, our strategy is to decompose the
debiasing term as a sum $\Debias_\numobs = \TermTil_2 + \TermTil_3 +
\TermTil_R$, where, the differences $\TermTil_j - \Term_j$ are
well-controlled for each $j \in \{2, 3 \}$, and the overall remainder
term $\TermTil_R$ can also be bounded.  More precisely, recall the
shorthand $\LOO{\PlainDelHat}{\numobs-1}{i} \defn
(\LOO{\EstPar}{\numobs-1}{i} - \EstPar_\numobs)$, corresponding to the
$i^{th}$-LOO error.  With this notation, we define
\begin{subequations}
  \begin{align}
\label{EqnDefnTermTil2}    
\TermTil_2 & \defn \sum_{i=1}^\Numobs \left \{ \Big[
  \frac{\Numobs-1}{\Numobs}
  \nabla\TargetFun(\EstPar_\numobs)^{\mytrans}- \frac{1}{\Numobs}
  \sum_{j\neq i}\IntSec_\TruePar^{\mytrans} \nabla_\Par
  \ZFun(\Sam_j,\EstPar_\numobs) \Big] \LOO{\PlainDelHat}{\numobs-1}{i}
\right \}, \quad \mbox{and} \\
\label{EqnDefnTermTil3}    
\TermTil_3 & \defn \frac{1}{2} \sum_{i=1}^\Numobs
(\LOO{\PlainDelHat}{\numobs-1}{i})^{\mytrans} \Big[
  \frac{\Numobs-1}{\Numobs} \nabla^2\TargetFun(\EstPar_\numobs)-
  \frac{1}{\Numobs}\sum_{j\neq i} \vtinprod{\IntSec_\TruePar}{
    \nabla^2_\Par \ZFun(\Sam_j,\EstPar_\numobs)}\Big]
\LOO{\PlainDelHat}{\numobs-1}{i}.
  \end{align}
\end{subequations}
Note the clear parallels to the definitions of $\TermTil_j$ and those
of $\Term_j$, for $j \in \{1, 2 \}$.

Overall, with this set-up, we see that there are two main pieces to
the remainder of the proof: (a) bounding the remainder term; and (b)
bounding each of the differences $\TermTil_j - \Term_j$ for $j \in
\{1, 2\}$.

\subsubsection{Bounding the remainder term}

The remainder term in our decomposition of $\Debias_\numobs$ is given
by $\TermTil_R \defn \Debias_\numobs - (\TermTil_2 + \TermTil_3)$.
Our first result gives a high probability bound on it:
\begin{lemma}
\label{LemDebiasRemainder}  
Under the conditions of~\Cref{ThmJackknife}, we have the bound
\begin{align*}
  |\TermTil_R| & \leq \polyshort \frac{\dimdelta^{3/2}
    \log^{6.5}(\Numobs/\delta)}{\Numobs^2}.
\end{align*}
with probability at least $1 - \delta$.
\end{lemma}
\noindent See~\Cref{AppLemDebiasRemainder} for the proof.


\subsubsection{Bounds on the differences $\TermTil_j - \Term_j$}

\noindent The next step in our proof---and the most challenging---is
to bound the differences $\TermTil_j - \Term_j$.
\begin{lemma}
  \label{LemHard}
Under the conditions of~\Cref{ThmJackknife}, we have the bounds
\begin{subequations}
    \begin{align}
\label{EqnHardTwo}      
|\TermTil_2 - \Term_2| & \leq \polyshort \Big(
\log^{4.5}(\Numobs/\delta) \frac{\dimdelta^{3/2}}{\Numobs^{3/2}} +
\log^{5}(\Numobs/\delta) \frac{ \sqrt{\ParDim}}{\Numobs}\Big), \quad
\mbox{and} \\
\label{EqnHardThree}        
  |\TermTil_3 - \Term_3| & \leq \polyshort \log^{5}(\Numobs/\delta)
  \frac{\dimdelta^{3/2}}{\Numobs^{3/2}},
\end{align}
\end{subequations}
where each statement holds with probability at least $1 - \delta$.
\end{lemma}
\noindent See~\Cref{AppLemHardTwo,AppLemHardThree} for the proofs of
the two bounds stated in this lemma.


\subsubsection{Combining the ingredients}

Recall the decomposition~\eqref{EqnJackDecomp}, as well as as the
representation $\Debias_\numobs = \TermTil_R + \TermTil_2 +
\TermTil_3$.  Combined with the triangle inequality, we have
\begin{align*}
\widebar E_\numobs \defn \Big| \sqrt{\numobs} \big(\JackEst -
\TargetFun(\TruePar) \big) - \noiseW_{\numobs} \Big| & =
\sqrt{\numobs} \Big | \Term - \Debias_\numobs + \RemainTau + \RemainZ
\Big| \\
& \leq \sqrt{\numobs} \Big \{ \underbrace{|\TermTil_R|}_{(a)} +
\underbrace{|\Term_2 - \TermTil_2| + |\Term_3 - \TermTil_3|}_{(b)} +
\underbrace{|\RemainTau| + |\RemainZ|}_{(c)} \Big \}.
\end{align*}
We now bound terms (a), (b) and (c) using
Lemmas~\ref{LemDebiasRemainder}, ~\ref{LemHard}
and~\ref{LemRemainTau}, respectively.  Doing so yields
\begin{align*}
  \widebar E_\numobs & \leq \sqrt{\numobs} \Biggr \{ \polyshort
  \frac{\dimdelta^{3/2} \log^{6.5}(\Numobs/\delta)}{\Numobs^2} +
  \polyshort \Big( \log^{5}(\Numobs/\delta)
  \frac{\dimdelta^{3/2}}{\Numobs^{3/2}} + \log^{5}(\Numobs/\delta)
  \frac{ \sqrt{\ParDim}}{\Numobs}\Big) + \polyshort
  \Big(\frac{\dimdelta}{\Numobs}\Big)^{3/2} \log^{2} (\Numobs/\delta)
  \Biggr \},\\ & \leq \polyshort\log^5(\numobs/\delta)
  \frac{\sqrt\ParDim}{\sqrt{\numobs}} +
  \polyshort\log^{5}(\numobs/\delta) \frac{\dimdelta^{3/2}}{\numobs} +
  \polyshort\log^{6.5}(\numobs/\delta)
  \frac{\dimdelta^{3/2}}{\numobs^{3/2}} \\ &\leq \polyshort
  \Big(\sqrt{\frac{\ParDim}\numobs}+\frac{\ParDim^{3/2}}{\numobs}\Big)\log^{13/2}(\numobs/\delta).
\end{align*}
which completes the proof of~\Cref{ThmJackknife}.


\subsection{Proof of Corollaries}
\label{SecProofCors}

We now prove the high-dimensional asymptotic claims from
Corollaries~\ref{CorHDAsymptoticPlugIn}
and~\ref{CorHDAsymptoticDebias}. Define the exponents $\rexplow
\mydefn \lim \inf \limits_{\numobs \rightarrow + \infty} \frac{\log
  \numobs}{\log \usedim}$ and $\rexpup \mydefn \lim \sup
\limits_{\numobs \rightarrow + \infty} \frac{\log \numobs}{\log
  \usedim}$, and note that~\Cref{ThmPlugin,ThmJackknife} guarantee
decompositions of the form
\begin{align}
\label{EqnCLTDecomp}  
\sqrt{\numobs} \big(\TargetFun(\EstPar_\numobs) - \TargetFun(\TruePar)
\big) \stackrel{(a)}{=} \StochNoise_\numobs + \frac{\highbias}{
  \sqrt{\numobs}} + \highplug, \quad \mbox{and} \quad \sqrt{\numobs}
\big(\JackEst - \TargetFun(\TruePar) \big) \stackrel{(b)}{=}
\StochNoise_\numobs + \highjack,
\end{align}
where the remainder terms $\highplug$ and $\highjack$ converge to zero
in probability for sequences with $\rexplow > 3/2$.

Our proof is based on establishing the following two auxilary claims.
First, under the conditions
of~\Cref{CorHDAsymptoticPlugIn}~(or~\ref{CorHDAsymptoticDebias}), we
have
\begin{subequations}
  \begin{align}
    \label{eq:CLT}
 \frac{\StochNoise_\numobs}{\TrueVar} \convdist \Normal (0, 1),
  \end{align}
where the variance $\TrueVar^2$ was previously
defined~\eqref{EqnDefnVarTrue}.  Second, under the conditions
of~\Cref{CorHDAsymptoticPlugIn}, we have
\begin{align}
\label{eq:highbias-asymptotics}
\frac{|\highbias|}{ \sqrt{\numobs}} \longrightarrow
\begin{cases} 0  & \mbox{if 
    $\rexplow > 2$, and} \\
  \infty & \mbox{if $\rexplow, \rexpup \in (3/2, 2)$}.
\end{cases}
\end{align}
\end{subequations}
We return to prove these auxiliary claims shortly.  Taking them as
given, let us complete the proof of the corollaries, beginning with
the high-dimensional asymptotics of the plug-in estimator.


\paragraph{Proof of~\Cref{CorHDAsymptoticPlugIn}:}

It follows from the variance scaling~\eqref{EqnVarScale} that
$\highplug/\VarTrue \convprob 0$, and moreover, \Cref{ThmPlugin}
guarantees that $\highplug \convprob 0$ when $\rexplow > 3/2$.  We
combine these facts with the decomposition~\eqref{EqnCLTDecomp}(a) and
the two auxiliary claims~\eqref{eq:CLT}
and~\eqref{eq:highbias-asymptotics}; doing so and applying Slutsky's
theorem yields that the ratio
\begin{align*}
 \sqrt{\numobs} \frac{(\TargetFun(\EstPar_\numobs) -
   \TargetFun(\TruePar))}{\TrueVar} \convdist \Normal(0,1) \qquad
 \mbox{when $\rexplow > 2$}
\end{align*}
and it diverges when $\rexplow, \rexpup \in (3/2, 2)$, as claimed.

\paragraph{Proof of~\Cref{CorHDAsymptoticDebias}:}
Similarly, \Cref{ThmJackknife} implies $\highjack \convprob 0$ when
$\rexplow > 3/2$.  We combine with this fact with the
decomposition~\eqref{EqnCLTDecomp}(b); applying Slutsky's theorem with
the auxiliary claims~\eqref{eq:CLT} yields
\begin{align*}
 \sqrt{\numobs} \frac{(\JackEst - \TargetFun(\TruePar))}{\TrueVar}
 \convdist \Normal(0,1), \quad \mbox{when $\rexplow > 3/2$,}
\end{align*}
as claimed. \\

\noindent It remains to establish the two auxiliary claims, namely
equations~\eqref{eq:CLT} and~\eqref{eq:highbias-asymptotics}.


\subsubsection{Proof of~\cref{eq:CLT}}

We verify this claim using the Lyapunov conditions for the central
limit theorem.  We have $\StochNoise_\numobs = \frac{1}{
  \sqrt{\numobs}} \sum_{i=1}^\numobs \StochNoisebase_i$, where the
random variables $\StochNoisebase_i \defn \inprod{\IntSec_\TruePar}
{\ZFun(\Sam_i, \TruePar)}$ are i.i.d. zero-mean variables with
$\Exs[\StochNoisebase_i^3] \leq \poly$.  Moreover, we note that $\lim
\inf_{\numobs \rightarrow + \infty} \TrueVar > 0$ by assumption.
Finally, we can write
\begin{align*}
\lim_{\numobs \rightarrow \infty} \frac{\sum_{i=1}^\numobs
  \Exs[\StochNoisebase_i^3]}{(\sum_{i=1}^\numobs
  \Exs[\StochNoisebase_i]^2 )^{3/2}} = \lim_{n\rightarrow \infty}
\frac{\numobs \Exs[\StochNoisebase_i^3]}{\numobs^{3/2}
  (\Exs[\StochNoisebase_i^2])^{3/2}} \leq \lim_{\numobs \rightarrow
  \infty} \frac{\poly\numobs}{\numobs^{3/2}} = 0,
\end{align*}
where the last inequality uses condition ~\eqref{EqnRegularScaling}.


\subsubsection{Proof of~\cref{eq:highbias-asymptotics}}

Under the condition $\rexplow, \rexpup\in (3/2, 2)$, the limit
$\highbias/\ParDim \not \longrightarrow 0$ implies that
$\frac{|\highbias|}{\sqrt{\numobs}} \rightarrow \infty$.  Recalling
the definition~\eqref{EqnDefnStochNoise} of $\Bias_\ParDim$, it can be
verified using condition~\ref{ass:tail} that $|\Bias_\ParDim| \leq
\poly \ParDim$.  When $\rexplow > 2$, we are guaranteed the existence
of some $\newexponent \in (0,1)$ such that $\ParDim^2/\numobs^{1
  -\newexponent} \to 0$.  Using these facts, we write
\begin{align*}
\lim_{\numobs \to \infty} \frac{|\Bias_\ParDim|}{\sqrt{\numobs}} &
\leq \lim_{\numobs \to \infty} \frac{\poly \ParDim}{\sqrt{\numobs}} \\
& \leq \underbrace{\lim_{\numobs \to \infty}
  \frac{\poly}{\sqrt{\numobs^{\newexponent}}}}_{ = 0} \cdot
\underbrace{\lim_{\numobs \to \infty}
  \sqrt{\frac{\Pardim^2}{\numobs^{1- \newexponent}}}}_{ = 0} \; = \;
0,
\end{align*}
where the final step uses the scaling
condition~\eqref{EqnRegularScaling} on $\poly$.



\section{Discussion}
\label{sec:discussion}

Resampling methods such as the jackknife are especially well-suited to
modern estimators for which only ``black-box'' access is available.
The main contribution of this paper is to provide a precise and
non-asymptotic answer to the question: for estimating functionals of
$Z$-estimators in high-dimensional regimes, does the jackknife correction
yield significant improvement over standard plug-in methods?  When
using $Z$-estimators in dimension $\usedim$, we show that a standard
plug-in approach for estimating a functional is
$\sqrt{\numobs}$-consistent and asymptotically normal whenever
$\ParDim^2 \lesssim \numobs$.  Beyond this boundary, we show that the
plug-in method typically fails to have this favorable property, but
that jackknife correction can rectify it in the regime $\ParDim^{3/2}
\lesssim \numobs$.  This separation between plug-in and jackknife
provides strong motivation for using jackknife correction---despite
its higher computational costs---in the latter regime of
$(\numobs,\usedim)$-pairs.

Our work opens up a number of directions for future research. From a
technical perspective, our theoretical results require the
$Z$-estimation problem to satisfy certain smoothness assumption (see
Assumption~\ref{ass:tail}). It is interesting to understand to what
extent the same results can be extended to certain types of non-smooth
problems (e.g., quantile estimation).  Moreover, we have focused
exclusively on the standard jackknife, but we suspect that our
analysis techniques might be useful in studying other resampling
procedures, among them the delete-d jackknife~\cite{shao1989general},
delta-a-group jackknife~\cite{kott2001delete}, as well as the
bootstrap method and its variants~\cite{efron1992bootstrap}.  The
current analysis has exhibited a sharp distinction between standard
plug-in and standard jackknife in the high-dimensional setting; it
would be interesting to see if these other procedures also exhibit
such sharp distinctions.

Last, our current analysis of $Z$-estimators has focused on data that
is generated in an i.i.d. manner.  However, $Z$-estimators are also
deployed to datasets that exhibit sequential dependency, including in
time series data, treatment assignment, and reinforcement learning
(e.g.,~\cite{hadad2021confidence,DuaWanWai24}).  Extending these
debiasing procedures to accommodate sequentially collected data poses
another important area for future research.


\printbibliography

\appendix

\section{Proofs: Lemmas from~\Cref{ThmPlugin}}
\label{AppThmPlugin}

\noindent In this appendix, we collect together the proofs of lemmas
used in the proof of~\Cref{ThmPlugin}.

\subsection{Proof of~\Cref{LemUstatBound}}
\label{SecProofLemUstatBound}

\noindent In this appendix, we prove the tail bound on the
$U$-statistic component stated in~\Cref{LemUstatBound}. We start the
proof by introducing a known result:
\begin{lemma}[Lemma 4 in
the paper~\cite{su2023estimated}]
\label{lemma:u-stats-concentration}
Consider an i.i.d. zero mean sequence $\{(X_i, Y_i) \}_{i=1}^\numobs$,
and suppose there exist scalars $v, \sigma > 0$ and $\alpha \in [1,2]$
such that
\begin{align}
\label{eq:ustats-lemma-condition}  
  \lammax(\Exs [X X^{\mytrans}]) , \lammax(\Exs [YY^{\mytrans}]) \leq
  v^2, \quad \mbox{and} \quad \vecnorm{\vecnorm{X}{2}}{\psi_\alpha},
  \vecnorm{\vecnorm{Y}{2}}{\psi_\alpha} \leq \sigma \sqrt{\ParDim}.
\end{align}
Then for any $\delta \in (0,1)$, we have
\begin{align*}
   \abss{\inprod{ \frac{1}{\numobs} \sum_{i = 1}^\numobs \State_i}{
       \frac{1}{\numobs} \sum_{i = 1}^\numobs \Outcome_i} -
     \frac{1}{\numobs} \Exs [\inprod{X}{Y}]} & \leq \frac{c v^2
     \sqrt{\ParDim}}{\numobs} \log (1 / \delta) + \frac{c_\alpha
     \sigma^2 \ParDim}{\numobs^{3/2}} \log^{1/ 2 + 4 / \alpha}
   (\numobs / \delta), \quad \mbox{and} \\
 \abss{ \frac{1}{\numobs^2} \sum_{i \neq j} \inprod{\State_i}{
     \Outcome_i} } & \leq \frac{c v^2 \sqrt{\ParDim}}{\numobs} \log (1
 / \delta) + \frac{c_\alpha \sigma^2 \ParDim}{\numobs^{3/2}} \log^{1/
   2 + 4 / \alpha} (\numobs / \delta)
\end{align*}
with probability at least $1 - \delta$.
\end{lemma}

Recall the definition of $\sioneave,\sitwoave,\Ustat$ in
equation~\eqref{eq:def_phi}---\eqref{EqnDefnUstat}.  Introduce the
shorthand notation $\sione = \sionebase(\Sam_i)$ and $\sitwo =
\sitwobase(\Sam_i)$\footnote{We use these shorthand notations
throughout the appendices.}, so that we have
\begin{align*}
\sioneave = \frac{1}{\numobs} \sum_{i=1}^\numobs \sione, \quad
\mbox{and} \quad \sitwoave = \frac{1}{\numobs} \sum_{i=1}^\numobs
\sitwo.
\end{align*}
Moreover, define the uncentered $U$-statistic $\UstatTil \defn
\sqrt{\numobs} \inprod{\sioneave + \sitwoave^{\mytrans} \Mmat/2 }{
  \sitwoave}$, so that \mbox{$\Ustat = \UstatTil - \Exs[\UstatTil]$}
by definition.  Observe that $\UstatTil/\sqrt{\numobs}$ is defined as
an inner product between two empirical averages, and that by
definition of $\sione$ and $\sitwo$, we have $\Exs[\sioneave] = 0$ and
$\Exs[\sitwoave] = 0$.  Our next step is to compute the second moments
and Orlicz norms of the random vectors $\sione$ and $\sitwo$.
Defining the scalars
\begin{subequations}
\begin{align}
\label{eqn:def_of_v_sigma}
v_1^2 \defn \lambda_{\max} ( \Exs[\sione \sione^{\mytrans}] ) \quad
\mbox{and} \quad v_2^2 \defn \lambda_{\max} ( \Exs[\sitwo
  \sitwo^{\mytrans}] ), \quad \mbox{as well as} \\
\label{EqnDefnSigtil}
\sigtil_1 \defn \frac{\vecnorm{ \vecnorm{\sione}{2}
  }{\psi_1}}{\sqrt{\ParDim}} \quad \mbox{and} \quad \sigtil_2 \defn
\frac{\vecnorm{ \vecnorm{\sitwo}{2} }{\psi_1}}{\sqrt{\ParDim}},
\end{align}
\end{subequations}
we claim that
\begin{align}
\label{EqnAuxBounds}  
v_1 \leq c \frac{\BGradTar \OrParb}{\strongconvex}, \quad \sigtil_1\le
c \frac{\BGradTar \OrParb}{\strongconvex} \log (\ParDim), \quad v_2
\leq c \frac{\OrPara}{\strongconvex},\quad \sigtil_2 \leq c
\frac{\OrPara}{\strongconvex} \log (\ParDim).
\end{align}
We return to prove these claims momentarily.

Note that we have the operator norm bound $\opnorm{\Mmat} \leq
\PlainLip \, (1 + \frac{\sigma}{\strongconvex})$, and let $u \in
\real^\ParDim$ be an arbitrary unit norm vector.  Making use of the
conditions~\ref{ass:smooth_function}, ~\ref{ass:tail}
\mbox{and~\ref{ass:convergence},} each zero-mean random variable
$\inprod{u}{\sione + \Mmat \sitwo}$ has variance at most
\begin{multline*}
\Exs \big[ \inprod{u}{\sione + \Mmat \sitwo }^2 \big] \leq c
\big(v_1^2 + \opnorm{\Mmat}^2 v_2^2 \big) \leq c ( \frac{\BGradTar
  \OrParb}{\strongconvex})^2 + c(\BGradTartwo +
\frac{\BGradTar\OrParc}{\strongconvex})^2 \OrPara^2/\strongconvex^2
\revdefn (v')^2,
\end{multline*}
where we have made use of the bounds~\eqref{EqnAuxBounds}.  Similarly,
we have the Orlicz norm bound
\begin{align*}
\vecnorm{\vecnorm{\sione + \Mmat \sitwo }{2}}{\psi_1} / \sqrt{\ParDim}
\leq \sigtil_1 + \opnorm{\Mmat} \sigtil_2 \leq c ( \frac{L
  \sigma}{\strongconvex}) + c(L + \frac{L \sigma}{\strongconvex})
\frac{\sigma}{\strongconvex} \log \ParDim \revdefn \sigma'.
\end{align*}

We can now apply~\Cref{lemma:u-stats-concentration} to the pair
$(\sioneave + \Mmat \sitwoave) / v'$ and $\sitwoave / v_2$.  By doing
so, we are guaranteed to have
\begin{align*}
 \frac{1}{v' v_2} \Big| \inprod{\sioneave + \Mmat \sitwoave}{
   \sitwoave} - \Exs [\UstatTil/ \sqrt{\numobs}] \Big| \leq \frac{c
   \sqrt{\ParDim}}{\numobs} \log (1 / \delta) + \frac{c \ParDim
 }{\numobs^{3/2}} \{ \log^2{\ParDim} \}\{ \log^{9/2}(\numobs / \delta)
 \},
\end{align*}
with probability at least $1 - \delta$.  Thus, we have shown that
\begin{align*}
\abss{\Ustat} & \leq \polyshort\left( \frac{c \sqrt{\ParDim}}{
  \sqrt\numobs} \log (1 / \delta) + \frac{c \ParDim }{\numobs} \{
\log^2{\ParDim} \}\{ \log^{9/2}(\numobs / \delta) \}\right) \; \leq \;
\polyshort \sqrt{ \frac{\ParDim}{\numobs}} \log^{5}(\numobs/\delta)
\end{align*}
which completes the proof of~\Cref{LemUstatBound}. \\


\noindent It remains to prove our auxiliary claim.
\paragraph{Proof of the bounds~\eqref{EqnAuxBounds}:}
  
For any $\usedim$-dimensional random vector $Z$, Pisier's
inequality~\cite{pisier1983some} (see also \cite[Exercise
  2.18(a)]{wainwright2019high}) implies that
\begin{align}
\label{eq:from-directional-orlicz-to-norm-orlicz}
\vecnorm{ \vecnorm{Z}{2} }{\psi_1} \leq \sqrt{\usedim} \cdot \vecnorm{
  \vecnorm{Z}{\infty} }{\psi_{1}} \leq c_\alpha \sqrt{\usedim} \, \log
\usedim\, \max_{j \in [\usedim]}
\vecnorm{\inprod{\coordinate_j}{Z}}{\psi_1},
\end{align}
for a constant $c > 0$. From the definition~\eqref{eqn:def_of_v_sigma}
of $v_1$, we have
\begin{align*}
v_1^2 = \max_{\|u\|_2 = 1} \Exs[(\inprod{\sione}{u})^2] & =
\max_{\|u\|_2 = 1} \Exs[( \nabla\TargetFun(\TruePar)^{\mytrans} u-
  \IntSec_\TruePar^{\mytrans} \nabla_\Par \ZFun(\Sam_i, \TruePar)
  u)^2] \\
& \leq \max_{\|u\|_2 = 1} \Exs[(\IntSec_\TruePar^{\mytrans}
  \nabla_\Par \ZFun(\Sam_i, \TruePar) u)^2] \; \leq \; \Big(
\frac{\BGradTar}{\strongconvex} \Big)^2 \OrParb^2.
\end{align*}
From condition~\ref{ass:tail}, the random variable
$\inprod{u}{\sione}$ has sub-exponential tails for any unit-norm
vector $u$, whence $\vecnorm{ u^{\mytrans} \sione }{\psi_1} \leq c
\frac{\BGradTar}{\strongconvex} \OrParb$ for some universal constant
$c > 0$. Therefore, applying the Pisier
bound~\eqref{eq:from-directional-orlicz-to-norm-orlicz} yields
$\sigtil_1 \leq c \frac{\BGradTar}{\strongconvex} \OrParb \log
\ParDim$.

Via a similar argument, we obtain
\begin{align*}
v_2^2 = \max_{\|u\|_2 = 1} \Exs[(\sitwo^{\mytrans} u)^2] =
\max_{\|u\|_2 = 1} \Exs[(\ZFun(\Sam_i,
  \TruePar)^{\mytrans}\GradFun_\TruePar^{-1} u)^2] \leq
\OrPara^2/\strongconvex^2.
\end{align*}
Since $\inprod{u}{\sitwo}$ is a sub-exponential variable, it follows
that $\vecnorm{\inprod{u}{\sitwo}}{\psi_1} \leq c
\OrPara/\strongconvex$. Finally, invoking the Pisier
inequality~\eqref{eq:from-directional-orlicz-to-norm-orlicz} gives the
bound $\sigtil_2 \leq c \frac{\OrPara \log (\ParDim)}{\strongconvex}$.

\subsection{Proof of~\Cref{LemLinErrorBound}}
\label{SecProofLemLinErrorBound}

In this section, we prove the upper bound~\eqref{EqnLinErrorBound} on
the first-order approximation error $\LinError_\numobs$.  Introduce
the shorthand $\DelHat \defn \EstPar_\numobs - \TruePar$.  We first
observe that
\begin{align*}
\frac{1}{\Numobs}\sum_{i=1}^\Numobs \ZFun(\Sam_i, \TruePar) & = -
\frac{1}{\Numobs}\sum_{i=1}^\Numobs \Big[\ZFun(\Sam_i,
  \EstPar_\numobs) - \ZFun(\Sam_i, \TruePar) \Big] \; = \;
\frac{1}{\Numobs}\sum_{i=1}^\Numobs \Big[\ZFun(\Sam_i, \TruePar +
  \DelHat) - \ZFun(\Sam_i, \TruePar) \Big],
\end{align*}
where the first equality follows from the defining
condition~\eqref{eq:z_estimate_1} for the the $Z$-estimate
$\EstPar_\numobs$.  We now perform a Taylor series expansion, thereby
obtaining
\begin{align*}
\frac{1}{\Numobs}\sum_{i=1}^\Numobs \ZFun(\Sam_i, \TruePar) & = -
\left \{ \frac{1}{\Numobs}\sum_{i=1}^\Numobs \nabla_\Par \ZFun(\Sam_i,
\TruePar) \DelHat + \frac{1}{\Numobs}\sum_{i=1}^\Numobs \int_0^1
(1-s){ \nabla^2_\Par \ZFun_j(\Sam_i, \Par(s))[\DelHat, \DelHat]}
\textup{d}s \right \}.
\end{align*}
Combining with the definition~\eqref{EqnDefnLinError} of $\Errfirst$,
we find that
\begin{align*}
\enorm{\Errfirst} & = \Bigg\|\GradFun_{\TruePar}^{-1} \left[
  \frac{1}{\Numobs} \sum_{i=1}^\Numobs \nabla_\Par \ZFun(\Sam_i,
  \TruePar) + \GradFun_{\TruePar} \right] \DelHat + \frac{1}{\Numobs}
\sum_{i=1}^\Numobs \GradFun_{\TruePar}^{-1} \int_0^1 (1-s) \big[
  \nabla^2_\Par \ZFun(\Sam_i, \Par(s))[\DelHat, \DelHat] \big]
\textup{d}s \Bigg\|_2 \\
& \leq \opnorm{\GradFun_{\TruePar}^{-1}} \Big \{\Term_1 + \Term_2
\Big\}
\end{align*}
where
\begin{align*}
\Term_1 \defn \opnorm{ \frac{1}{\Numobs}\sum_{i=1}^\Numobs \nabla_\Par
  \ZFun(\Sam_i, \TruePar) + \GradFun_{\TruePar} } \enorm{\DelHat},
\quad \mbox{and} \quad \Term_2 \defn \sup_{\Par \in \ParSpace}
\opnorm{ \frac{1}{\Numobs}\sum_{i=1}^\Numobs \nabla^2_\Par
  \ZFun(\Sam_i, \Par)} \enorm{\DelHat}^2
\end{align*}
Note that $\opnorm{\GradFun_{\TruePar}^{-1}} \leq 1/\strongconvex$ by
assumption.  Thus, it remains to bound $\Term_1$ and $\Term_2$.
Recalling the shorthand $\dimdelta = \ParDim + \log(1/\delta)$,
condition~\ref{ass:convergence} ensures that $\enorm{\DelHat} \leq
\polyshort \sqrt{ \frac{\dimdelta \log(\numobs)}{\Numobs}}$ with
probability at least $1 - \delta$. \Cref{lm:op_norm_zfun} implies that
\begin{align*}
\opnorm{ \frac{1}{\Numobs}\sum_{i=1}^\Numobs \nabla_\Par \ZFun(\Sam_i,
  \TruePar) + \GradFun_{\TruePar} } \leq \polyshort \Big( \sqrt{
  \frac{\Ccerr}{\Numobs}} + \frac{\dimdelta \log \Numobs}{\Numobs}
\Big)
\end{align*}
with probability at least $1 - \delta$.  Combining the two pieces and
using the sample size condition, we find that
\begin{subequations}
\begin{align}
\label{EqnMalaOne}
\Term_1 & \leq \polyshort \sqrt{ \frac{\dimdelta
    \log(\numobs)}{\Numobs}} \Big( \sqrt{ \frac{\Ccerr}{\Numobs}} +
\frac{\dimdelta \log \Numobs}{\Numobs} \Big) \leq \polyshort
\frac{\dimdelta \log\numobs}{\Numobs}
\end{align}
with probability at least $1-\delta.$ Similarly,
\Cref{lm:sup_op_norm_zfun} implies that $\sup_{\Par \in \ParSpace}
\opnorm{ \frac{1}{\Numobs}\sum_{i=1}^\Numobs \nabla^2_\Par
  \ZFun(\Sam_i, \Par)} \leq c \OrParc\Big(1 + \sqrt{
  \frac{\dimdelta}{\Numobs}} + \frac{(\dimdelta \, (1+\log
  \Numobs))^{\InvOrConc}}{\Numobs} \Big)$.  Again, combining with the
bound on $\enorm{\DelHat}$ and the sample size condition, we find that
\begin{align}
  \label{EqnMalaTwo}
\Term_2 & \leq \polyshort { \frac{\dimdelta \log(\numobs)}{\Numobs}}
\Big(1 + \sqrt{ \frac{\dimdelta}{\Numobs}} + \frac{(\dimdelta \,
  (1+\log \Numobs))^{\InvOrConc}}{\Numobs} \Big) \leq \polyshort {
  \frac{\dimdelta \log(\numobs)}{\Numobs}}.
\end{align}
\end{subequations}
Finally, combining the two bounds~\eqref{EqnMalaOne}
and~\eqref{EqnMalaTwo}, we arrive at the bound
\begin{align*}
\enorm{\Errfirst} & \leq \polyshort \frac{\dimdelta
  \log\numobs}{\Numobs}
\end{align*}
with probability at least $1-\delta.$


\subsection{Proof of~\Cref{LemRemainTau}}
\label{sec:proof_LemRemainTau}

We make use of the shorthand $\DelHat \defn \EstPar_\numobs -
\TruePar$.  By the Lipschitz continuity of $ \nabla^2\TargetFun$, we
can write
\begin{align*}
|\RemainTau| & = \Big| \DelHat^{\mytrans} \Big[\int_0^1(1-s)[
    \nabla^2\TargetFun(\TruePar + s \DelHat) -
    \nabla^2\TargetFun(\TruePar)]\textup{d}s\Big] \DelHat \Big| \\
& \leq \enorm{\DelHat}^2 \Big[ \int_0^1 (1-s) \opnorm{ \nabla^2
    \TargetFun(\TruePar + s \DelHat) - \nabla^2 \TargetFun(\TruePar)}
  \textup{d}s \Big ] \; \leq \; \LipTarc \enorm{\DelHat}^3.
\end{align*}
Thus, using condition~\ref{ass:convergence} to bound
$\enorm{\DelHat}$, we find that
\begin{align*}
|\RemainTau| \leq \LipTarc (\polyshort \sqrt{ \log \Numobs} \sqrt{
  \frac{\myusedim}{\Numobs}})^{3} \leq \polyshort (
\frac{\myusedim}{\Numobs})^{3/2} \log^2 \numobs,
\end{align*}
as claimed in equation~\eqref{EqnRemainTauBound}. \\

Turning to the claimed bound on $\RemainZ$, we write
\begin{align}
|\RemainZ| & = \Biggr | \int_{0}^1 (1-s) \DelHat^{\mytrans} \Big\{
\frac{1}{\Numobs} \sum_{i=1}^\Numobs[{ \nabla^2_\Par \ZFun(\Sam_i,
    \TruePar + s \DelHat) [\IntSec_\TruePar]- \nabla^2_\Par
    \ZFun(\Sam_i, \TruePar)}[\IntSec_\TruePar]] \Big\} \DelHat
\textup{d}s \Biggr | \notag \\
& \leq \enorm{\DelHat}^2 \int_{0}^1 (1-s) \opnorm{ \frac{1}{\Numobs}
  \sum_{i=1}^\Numobs[{ \nabla^2_\Par \ZFun(\Sam_i, \TruePar + s
      \DelHat) [\IntSec_\TruePar]- \nabla^2_\Par \ZFun(\Sam_i,
      \TruePar)}[\IntSec_\TruePar]]} \textup{d}s \notag \\
& \leq \enorm{\DelHat}^3 \int_{0}^1 (1-s)s
\sup_{\Par\in \ParSpace}\opnorm{ \frac{1}{\Numobs}\sum_{i=1}^\Numobs[{
      \nabla^3_\Par \ZFun(\Sam_i, \Par)[\IntSec_\TruePar]}}
  \textup{d}s \notag \\
\label{eq:term_5_bound_eq1}  
& \leq c \frac{\OrPard\BGradTar}{\strongconvex} \sqrt{ \log
  (\Numobs/\delta)}\Big(1 + \sqrt{ \frac{\dimdelta}{\Numobs}} +
\frac{(\dimdelta \log {\Numobs})^{\InvOrCond}}{\Numobs}\Big) \;
\enorm{\DelHat}^3,
\end{align}
with probability at least $1 - \delta$, where the fourth line
uses~\Cref{lm:sup_op_norm_zfun} and the fact that
$\enorm{\IntSec_\TruePar} \leq
\opnorm{\GradFun_{\TruePar}^{-1}}\enorm{ \nabla \TargetFun(\TruePar)}
\leq \BGradTar/\strongconvex$.

We now impose the sample size
requirement~\eqref{eq:linear_approx_Z_est_sample_ahead}, and make use
condition~\ref{ass:convergence} to bound $\enorm{\DelHat}$.  In this
way, we find that
\begin{align*}
\abss{\RemainZ} & \leq c \frac{\OrPard\BGradTar}{\strongconvex} \sqrt{
  \log (\Numobs/\delta)}\Big(1 + \sqrt{ \frac{\dimdelta}{\Numobs}} +
\frac{(\dimdelta \log {\Numobs})^{3/2}}{\Numobs}\Big) (\polyshort
\sqrt{ \log \Numobs} \sqrt{ \frac{\myusedim}{\Numobs}})^{3}\\ &\le
\polyshort\sqrt{\log(\numobs/\delta)}(\sqrt{\log(\numobs)} \sqrt{
  \frac{\myusedim}{\Numobs}})^{3} \leq
\polyshort\Big(\frac{\dimdelta}{\numobs}\Big)^{3/2}\log^{2}
(\Numobs/\delta),
\end{align*}
with probability at least $1 - \delta$, as claimed in
equation~\eqref{EqnRemainZBound}.


\subsection{Proof of~\Cref{LemPieces}}
\label{SecProofLemPieces}

Introduce the shorthand $\Err_\numobs \defn \EstPar_\numobs -
\TruePar$.  Performing a Taylor series expansion of the function
$\TargetFun$ between $\EstPar_\numobs$ and $\TruePar$, we find that
\begin{subequations}
\begin{multline}
\label{eq:targetfun_taylor_1}  
\TargetFun(\EstPar_\numobs) = \TargetFun(\TruePar) + \inprod{ \nabla
  \TargetFun(\TruePar)}{\Err_\numobs} + \frac{1}{2}
\Err_\numobs^{\mytrans} \Big[\nabla^2 \TargetFun(\TruePar) \Big]
\Err_\numobs \\
+ \Err_\numobs^{\mytrans} \left \{ \int_0^1 (1-s) \Big[ \nabla^2
  \TargetFun(\TruePar + s (\EstPar_\numobs - \TruePar)) -
  \nabla^2\TargetFun(\TruePar) \Big] \textup{d}s \right \}
\Err_\numobs.
\end{multline}
Defining the function $F(\theta) \defn \frac{1}{\numobs}
\sum_{i=1}^\Numobs \inprod{\IntSec_\TruePar}{\ZFun(\Sam_i, \theta)}$,
a second Taylor series expansion yields
\begin{multline}
\label{eq:gradfun_taylor_1}  
F(\thetahat_\numobs) = F(\thetastar) + \frac{1}{\Numobs}
\sum_{i=1}^\Numobs \IntSec_\TruePar^{\mytrans} \nabla_\Par
\ZFun(\Sam_i, \TruePar) \Err_\numobs + \frac{1}{2\Numobs}
\sum_{i=1}^\Numobs \Err_\numobs^{\mytrans}
\vtinprod{\IntSec_\TruePar}{ \nabla^2_\Par \ZFun(\Sam_i, \TruePar)}
\Err_\numobs \\
+ \Err_\numobs^{\mytrans} \left \{ \frac{1}{\Numobs}
\sum_{i=1}^\Numobs \Big[ \int_{0}^1 (1-s)
  [\vtinprod{\IntSec_\TruePar}{ \nabla^2_\Par \ZFun \big(\Sam_i,
      \TruePar + s \Err_\numobs \big)- \nabla^2_\Par \ZFun(\Sam_i,
      \TruePar)}] \textup{d}s \Big] \right \} \Err_\numobs.
\end{multline}
\end{subequations}
Now observe that $F(\thetahat_\numobs) = 0$ by definition of the
$Z$-estimator.  Using this fact, and combining
equations~\eqref{eq:targetfun_taylor_1}
and~\eqref{eq:gradfun_taylor_1}, we find that
$\TargetFun(\EstPar_\numobs) - \TargetFun(\TruePar)=
\frac{\noiseW_{\numobs}}{ \sqrt{\numobs}} + \Term_2 + \Term_3 +
\RemainTau + \RemainZ$, where the stochastic noise term
$\noiseW_\numobs$ was defined in equation~\eqref{EqnDefnStochNoise};
and
\begin{subequations}
  \begin{align}
\label{EqnTermTwo}    
\Term_2 & \defn \Big[\nabla \TargetFun(\TruePar)^{\mytrans} -
  \frac{1}{\Numobs}\sum_{i=1}^\Numobs \IntSec_\TruePar^{\mytrans}
  \nabla_\Par \ZFun(\Sam_i, \TruePar) \Big] \cdot \Err_\numobs, \\
\label{EqnTermThree}    
\Term_3 & \defn \frac{1}{2} \Err_\numobs^{\mytrans} \left \{
\nabla^2\TargetFun(\TruePar) - \frac{1}{\Numobs} \sum_{i=1}^\Numobs
\vtinprod{\IntSec_\TruePar}{ \nabla^2_\Par \ZFun(\Sam_i, \TruePar)}
\right \} \Err_\numobs, \\
\RemainTau & \defn \Err_\numobs^{\mytrans} \left \{ \int_0^1(1-s)[
  \nabla^2 \TargetFun \big(\TruePar + s \Err_\numobs \big) - \nabla^2
  \TargetFun(\TruePar)] \textup{d}s \right \} \Err_\numobs, \notag \\
\RemainZ & \defn - \Err_\numobs^{\mytrans} \Big \{ \frac{1}{\Numobs}
\sum_{i=1}^\Numobs \left[ \int_{0}^1 (1-s)
  [\vtinprod{\IntSec_\TruePar}{ \nabla^2_\Par \ZFun \big(\Sam_i,
      \TruePar + s \Err_\numobs \big)- \nabla^2_\Par \ZFun(\Sam_i,
      \TruePar)}] \textup{d}s \right] \Big \} \Err_\numobs. \notag
  \end{align}
\end{subequations}

In order to complete the proof, we need to show that $\Term_2 +
\Term_3 = \frac{\Ustat}{\sqrt{\numobs}} + \frac{\highbias}{\numobs} +
\RemainL $.  A little bit of algebra shows that
\begin{align*}
\Term_2 + \Term_3 = \underbrace{\inprod{\sioneave}{ \sitwoave} +
  \sitwoave^{\mytrans} \Mmat \sitwoave/2}_{\equiv
  \Ustat/\sqrt{\Numobs} + \highbias/\numobs} + \underbrace{
  \inprod{\sioneave}{\Errfirst} + \Errfirst^{\mytrans} \Mhat
  \Errfirst/2 + \sitwoave^{\mytrans} \Mhat \Errfirst +
  \sitwoave^{\mytrans} (\Mhat - \Mmat) \sitwoave/2}_{\revdefn
  \RemainL},
\end{align*}
which completes the proof.


\section{Proofs: Lemmas for~\Cref{ThmJackknife}}
\label{AppThmJackknife}

\noindent In this appendix, we collect together the proofs of lemmas
used in the proof of~\Cref{ThmJackknife}.  The bound on LOO
differences from~\Cref{LemJackknifeKey} is of independent interest and
has an involved proof, so that we devote a separate section to it
(\Cref{sec:proof_LemJackknifeKey}).


\subsection{Proof of~\Cref{LemDebiasRemainder}}
\label{AppLemDebiasRemainder}

\noindent Recall the definitions~\eqref{EqnPieceOne},
~\eqref{EqnDefnTermTil2}, and~\eqref{EqnDefnTermTil3} of
$\Debias_\numobs$, $\TermTil_2$ and $\TermTil_3$, respectively, as
well as the shorthand $\LooDiff = \LOO{\EstPar}{\numobs-1}{i} -
\EstPar_\numobs$.  The remainder term is given by $\TermTil_R =
\Debias_\numobs - (\TermTil_2 + \TermTil_3)$.

\subsubsection{Main argument}

Define the matrices
\begin{align*}
\QmatTau_{i} & \defn \Big[ \int_0^1(1-s)[ \nabla^2
    \TargetFun(\EstPar_\numobs + s \LooDiff) - \nabla^2
    \TargetFun(\EstPar_\numobs)] \textup{d}s \Big], \quad \mbox{and}
\\
\QmatZ_{i,j} & \defn \int_{0}^1 (1-s) [\vtinprod{\IntSec_\TruePar}{
    \nabla^2_\Par \ZFun \big(\Sam_j, \EstPar_\numobs + s \LooDiff
    \big)- \nabla^2_\Par \ZFun(\Sam_j, \EstPar_\numobs)}] \textup{d}s.
\end{align*}
Using this notation, we claim that the remainder can be split as
\begin{align}
\label{EqnDefnTermTilPair}
\TermTil_R & = \underbrace{ \sum_{i=1}^\Numobs (\LooDiff)^{\mytrans}
  \QmatTau_{i} \LooDiff}_{\revdefn \TermTil_4} + \underbrace{ -
  \frac{1}{\Numobs}\sum_{i=1}^\Numobs \sum_{j\neq i} \Big \{
  (\LooDiff)^{\mytrans} \QmatZ_{i,j} \LooDiff \Big \}}_{\revdefn
  \TermTil_5}.
\end{align}
We return to prove this auxiliary claim momentarily.

By the definition of $\STerm_4$, we have
\begin{align*}
|\STerm_4| & \leq \sum_{i=1}^\Numobs |(\LooDiff)^{\mytrans} \Big[
  \int_0^1 (1-s) [ \nabla^2\TargetFun(\EstPar_\numobs + s \LooDiff) -
    \nabla^2\TargetFun(\EstPar_\numobs)] \textup{d}s \Big] \LooDiff|\\
& \leq \sum_{i=1}^\Numobs \Big[\int_0^1(1-s) \opnorm{
    \nabla^2\TargetFun(\EstPar_\numobs + s \LooDiff) -
    \nabla^2\TargetFun(\EstPar_\numobs)} \textup{d}s \Big]
\enorm{\LooDiff}^2 \\
& \leq \sum_{i=1}^\Numobs \Big[\int_0^1 \LipTarc (1-s)s \textup{d}s
  \Big] \enorm{\LooDiff}^3 \; \leq \; \LipTarc \sum_{i=1}^\Numobs
\enorm{\LooDiff}^3.
\end{align*}
Similarly, we have
\begin{align*}
|\STerm_5| & \leq \frac{1}{\Numobs} \sum_{i=1}^\Numobs |
(\LooDiff)^{\mytrans} \int_{0}^1 (1-s)
[\vtinprod{\IntSec_\TruePar}{\sum_{j\neq i} \nabla^2_\Par \ZFun
    \big(\Sam_j, \EstPar_\numobs + s \LooDiff \big) - \sum_{j \neq i}
    \nabla^2_\Par \ZFun(\Sam_j, \EstPar_\numobs)}] \textup{d} s
\LooDiff | \\
& \leq \big \{ \int_{0}^1 (1-s) s \textup{d}s \big \} \; \;
\sum_{i=1}^\Numobs \sup_{\Par \in \ParSpace} \opnorm{
  \frac{1}{\Numobs-1} \sum_{j\neq i} \nabla^3_\Par \ZFun(\Sam_j, \Par)
       [\IntSec_\TruePar] } \; \enorm{\LooDiff}^3 \\
 & \leq c \OrPard \sqrt{ \log (\Numobs/\delta)} \Big(1 + \sqrt{
  \frac{\myusedim}{\Numobs}} + \frac{(\myusedim \log
  {\Numobs})^{3/2}}{\Numobs} \Big) \sum_{i=1}^\Numobs
\enorm{\LooDiff}^3 \\
& \leq c \OrPard \sqrt{ \log (\Numobs/\delta)} \sum_{i=1}^\Numobs
\enorm{\LooDiff}^3
\end{align*}
where the third line follows
from~\Cref{lm:sup_op_norm_zfun}. Substituting
equation~\eqref{eq:LemJackknifeKey_eq1} into the upper bounds for
$\STerm_4$ and $\STerm_5$ yields the claim.


\subsubsection{Proof of the decomposition~\eqref{EqnDefnTermTilPair}}

Performing a Taylor expansion for $\TargetFun$ around
$\EstPar_\numobs$ yields
\begin{multline}
\label{eq:targetfun_taylor_2}  
\TargetFun(\LOO{\EstPar}{\numobs-1}{i}) - \TargetFun(\EstPar_\numobs)
= \nabla \TargetFun(\EstPar_\numobs)^{\mytrans} \LooDiff + \frac{1}{2}
(\LooDiff)^{\mytrans} \nabla^2\TargetFun(\EstPar_\numobs) \LooDiff \\
+ (\LooDiff)^{\mytrans} \Big[\int_0^1(1-s) [ \nabla^2
    \TargetFun(\EstPar_\numobs + s \LooDiff) - \nabla^2
    \TargetFun(\EstPar_\numobs)] \textup{d}s \Big] \LooDiff
\end{multline}
Summing over $i$, we obtain
 \begin{align}
  &~~\sum_{i=1}^\Numobs [\TargetFun(\LOO{\EstPar}{\numobs-1}{i}) -
     \TargetFun(\EstPar_\numobs)] \notag \\
& = \nabla\TargetFun(\EstPar_\numobs)^{\mytrans} \sum_{i=1}^\Numobs
   \LooDiff
   + \frac{1}{2} \sum_{i=1}^\Numobs (\LooDiff)^{\mytrans}
   \nabla^2\TargetFun(\EstPar_\numobs) \LooDiff \notag \\
& \qquad + \sum_{i=1}^\Numobs (\LooDiff)^{\mytrans}
   \Big[\int_0^1(1-s)[ \nabla^2 \TargetFun(\EstPar_\numobs + s
       \LooDiff) - \nabla^2 \TargetFun(\EstPar_\numobs)]
     \textup{d}s\Big] \LooDiff
 \end{align}
By~\cref{EqnDefnLOO}, we have
\begin{align*}
 \ZFun(\Sam_i, \EstPar_\numobs)= \sum_{j\neq
   i}(\ZFun(\Sam_j,\LOO{\EstPar}{\numobs-1}{i}) -
 \ZFun(\Sam_j,\EstPar_\numobs)).
\end{align*}
By Taylor expansion, we have
\begin{align}
 \IntSec_\TruePar^{\mytrans} \ZFun(\Sam_i, \EstPar_\numobs) & =
 \sum_{j\neq i} \IntSec_\TruePar^{\mytrans} \nabla_\Par
 \ZFun(\Sam_j,\EstPar_\numobs) \LooDiff + \frac{1}{2} \sum_{j\neq i}
 (\LooDiff)^{\mytrans} \vtinprod{\IntSec_\TruePar}{ \nabla^2_\Par
   \ZFun(\Sam_j,\EstPar_\numobs)} \LooDiff \notag \\
\label{eq:gradfun_taylor_2} 
 & \quad + \sum_{j\neq i} (\LooDiff)^{\mytrans} \int_{0}^1
(1-s)[\vtinprod{\IntSec_\TruePar}{ \nabla^2_\Par
    \ZFun(\Sam_j,\EstPar_\numobs + s \LooDiff) - \nabla^2_\Par
    \ZFun(\Sam_j,\EstPar_\numobs)}] \textup{d}s \LooDiff.
\end{align}

Putting together the pieces and using equation~\eqref{eq:z_estimate_1}
completes the proof.


\subsection{Proof of bound~\eqref{EqnHardTwo} from~\Cref{LemHard}}
\label{AppLemHardTwo} 

In this appendix, we prove the bound~\eqref{EqnHardTwo} on the
difference $\TermTil_2 - \Term_2$.  We begin by writing $\STerm_2 =
\STerm_{21} + \STerm_{22}$, where
\begin{align*}
\STerm_{21} & \defn \sum_{i=1}^\Numobs \Big[ \frac{\Numobs-1}{\Numobs}
  \nabla \TargetFun(\TruePar)^{\mytrans} - \frac{1}{\Numobs} \sum_{j
    \neq i} \IntSec_\TruePar^{\mytrans} \nabla_\Par \ZFun(\Sam_j,
  \TruePar) \Big] (\LOO{\EstPar}{\numobs-1}{i} - \EstPar_\numobs), \\
\STerm_{22} & \defn \sum_{i=1}^\Numobs \Big \{
\frac{\Numobs-1}{\Numobs} \big[ \nabla \TargetFun(\EstPar_\numobs) -
  \nabla \TargetFun(\TruePar) \big]^{\mytrans} - \frac{1}{\Numobs}
\sum_{j \neq i} \IntSec_\TruePar^{\mytrans} \big[ \nabla_\Par
  \ZFun(\Sam_j, \EstPar_\numobs) - \nabla_\Par \ZFun(\Sam_j, \TruePar)
  \big] \Big \} (\LOO{\EstPar}{\numobs-1}{i} - \EstPar_\numobs).
\end{align*}
The central part of our proof involves showing that
\begin{subequations}
\begin{align}
\label{lm:compare_term_2_auxi_1}  
|\Term_2 - \STerm_{21}| & \leq \polyshort\Big( \log^4(\Numobs/\delta)
\frac{\dimdelta^{3/2}}{\Numobs^{3/2}} + \log^5(\Numobs/\delta) \frac{
  \sqrt{\ParDim}}{\Numobs} \Big) \\
\label{lm:compare_term_2_auxi_2}
|\STerm_{22}| & \leq \polyshort \log^{4.5}(\Numobs/\delta)
\frac{\dimdelta^{3/2}}{\Numobs^{3/2}},
\end{align}
\end{subequations}
with probability at least $1 - \delta$.  Observe that the
claim~\eqref{EqnHardTwo} follows as a consequence of these two
auxiliary claims. \\

\noindent We prove claims~\eqref{lm:compare_term_2_auxi_1}
and~\eqref{lm:compare_term_2_auxi_2} in
Appendices~\ref{sec:proof_lm:compare_term_2_auxi_1}
and~\ref{sec:proof_lm:compare_term_2_auxi_2}, respectively.

\subsubsection{Proof of~\cref{lm:compare_term_2_auxi_1}}
\label{sec:proof_lm:compare_term_2_auxi_1}
We have
\begin{align*}
 & \quad \Term_2 - \STerm_{21} \\ & = \Big[ \nabla
    \TargetFun(\TruePar)^{\mytrans} - \frac{1}{\Numobs}
    \sum_{i=1}^\Numobs \IntSec_\TruePar^{\mytrans} \nabla_\Par
    \ZFun(\Sam_j, \TruePar) \Big] (\EstPar_\numobs - \TruePar) - \Big[
    \nabla\TargetFun(\TruePar)^{\mytrans}- \frac{1}{\Numobs}
    \sum_{i=1}^\Numobs\IntSec_\TruePar^{\mytrans} \nabla_\Par
    \ZFun(\Sam_j, \TruePar) \Big]
  \sum_{i=1}^\Numobs(\LOO{\EstPar}{\numobs-1}{i} - \EstPar_\numobs) \\
  & \quad + \frac{1}{\Numobs} \sum_{i=1}^\Numobs[
    \nabla\TargetFun(\TruePar)^{\mytrans} -
    \IntSec_\TruePar^{\mytrans} \nabla_\Par \ZFun(\Sam_i,
    \TruePar)](\LOO{\EstPar}{\numobs-1}{i} - \EstPar_\numobs) \\
& = \sioneave^{\mytrans} (\Errfirst + \sitwoave -
  \sum_{i=1}^\Numobs(\LOO{\EstPar}{\numobs-1}{i} - \EstPar_\numobs)) +
  \frac{1}{\Numobs} \sum_{i=1}^\Numobs
  \sione^{\mytrans}(\LOO{\EstPar}{\numobs-1}{i} - \EstPar_\numobs)
  \\ & = \Res_{21} + \Res_{22} + \Res_{23},
\end{align*}
where
\begin{align*}
\Res_{21} = \sioneave^{\mytrans} [\Errfirst -
  \sum_{i=1}^\Numobs(\LOO{\EstPar}{\numobs-1}{i} - \EstPar_\numobs)],
\quad \Res_{22} = \frac{1}{\Numobs^2} \sum_{i\neq j} \sione^{\mytrans}
\sitwo, \quad \Res_{23} = \frac{1}{\Numobs} \sum_{i=1}^\Numobs
\sione^{\mytrans} (\LOO{\EstPar}{\numobs-1}{i} - \EstPar_\numobs +
\frac{1}{\Numobs} \sitwo).
\end{align*}
By~\Cref{lm:op_norm_zfun} and $\enorm{\IntSec_{\TruePar}} \leq
\BGradTar/\strongconvex$, with probability at least $1 - \delta$,
\begin{align}
\enorm{\sioneave} \leq \opnorm{ \frac{1}{\Numobs} \sum_{i=1}^\Numobs
  \nabla_\Par \ZFun(\Sam_i, \TruePar) - \E[ \nabla_\Par \ZFun(\Sam_i,
    \TruePar)]} \enorm{\IntSec_{\TruePar}} \leq \polyshort \sqrt{
  \frac{\dimdelta}{\Numobs}}.\label{eq:concen_graident_targetfun}
\end{align} 
Since $\rnorm{\sione^{\mytrans} \Direc }{\psi_1} \leq
c\OrParb\BGradTar/\strongconvex$ for any $\Direc \in
\Sphere{\ParDim-1}$, applying concentration of sub-exponential
variables leads to
\begin{align} \label{eq:concen_graident_targetfun2}
\Prob(\max_{i} \enorm{\sione} \geq \sqrt{\ParDim}t_1) \leq
\Prob(\max_{i,j}|\sione^{\mytrans} e_j|\geq t_1) \leq
2\Numobs\ParDim\exp\{- \frac{t_1\strongconvex}{\OrParb\BGradTar} \}.
\end{align}
Choose $t_1=c\OrParb\BGradTar \log
(\Numobs\ParDim/\delta)/\strongconvex$
in~\cref{eq:concen_graident_targetfun2}, and
applying~\cref{eq:concen_graident_targetfun},~\Cref{LemLinErrorBound,LemJackknifeKey},
we have with probability $1 - \delta$
\begin{align}
|\Res_{21} + \Res_{23}| & \leq \polyshort \sqrt{
  \frac{\dimdelta}{\Numobs}} \cdot\polyshort \frac{{\dimdelta}
  \log^4(\Numobs/\delta)}{\Numobs} \notag \\ & \quad + \polyshort
\sqrt{\dimdelta} \log (\Numobs/\delta) \cdot \polyshort
\log^3(\Numobs/\delta) \frac{\dimdelta}{\Numobs^{3/2}} \notag \\ &
\leq \polyshort \log^4(\Numobs/\delta)
\frac{\dimdelta^{3/2}}{\Numobs^{3/2}}.
\label{eq:bound_term_21_23}
\end{align}
Similar to the proof of Lemma~\ref{LemUstatBound}, for the
U-statistics $\Remainder_{22}$ we have
\begin{align}
\label{eq:bound_term_22}   
|\Res_{22}| & \leq \polyshort \frac{ \sqrt{\ParDim}
  \log^{5}(\Numobs/\delta)}{\Numobs}
\end{align}
with probability at least $1 -
\delta$. Combining~\cref{eq:bound_term_21_23,eq:bound_term_22}
completes the proof.


\subsubsection{Proof of~\Cref{lm:compare_term_2_auxi_2}}
\label{sec:proof_lm:compare_term_2_auxi_2}

Performing a Taylor expansion at $\TruePar$ and reorganizing the
terms, we find
\begin{align*}
\STerm_{22} & = \sum_{i=1}^\Numobs \Big[
  \frac{\Numobs-1}{\Numobs}(\EstPar_\numobs - \TruePar)^{\mytrans}
  \nabla^2\TargetFun(\TruePar)(\LOO{\EstPar}{\numobs-1}{i} -
  \EstPar_\numobs)- \frac{1}{\Numobs} \sum_{j \neq i} (\EstPar_\numobs
  - \TruePar)^{\mytrans} \vtinprod{\IntSec_\TruePar}{ \nabla^2_\Par
    \ZFun(\Sam_j, \TruePar)}(\LOO{\EstPar}{\numobs-1}{i} -
  \EstPar_\numobs) \Big] \\
& \quad + \sum_{i=1}^\Numobs \Big[
  \frac{\Numobs-1}{\Numobs}(\EstPar_\numobs - \TruePar)^{\mytrans}
  \Big[\int_{0}^1[ \nabla^2\TargetFun(\TruePar + t(\EstPar_\numobs -
      \TruePar)) - \nabla^2
      \TargetFun(\TruePar)]dt\Big](\LOO{\EstPar}{\numobs-1}{i} -
  \EstPar_\numobs) \Big] \\
& \quad - \frac{1}{\Numobs} \sum_{i=1}^\Numobs (\EstPar_\numobs -
\TruePar)^{\mytrans} \vtinprod{\IntSec_\TruePar}{\sum_{j \neq i}
  \int_{0}^1[ \nabla^2_\Par \ZFun(\Sam_j, \TruePar + t(\EstPar_\numobs
    - \TruePar)) - \nabla^2_\Par \ZFun(\Sam_j,
    \TruePar)]dt}(\LOO{\EstPar}{\numobs-1}{i} - \EstPar_\numobs) \\
& = \Res_{4} + \Res_{5} + \Res_{6} + \Res_7 + \Res_8,
\end{align*}
where
\begin{align*}
\Res_{4} & = (\EstPar_\numobs - \TruePar)^{\mytrans} \Mhat
\Big[\sum_{i=1}^\Numobs(\LOO{\EstPar}{\numobs-1}{i} - \EstPar_\numobs)
  \Big], \\
\Res_{5} & = - \frac{1}{\Numobs} \sum_{i=1}^\Numobs (\EstPar_\numobs -
\TruePar)^{\mytrans} \Big[ \nabla^2\TargetFun(\TruePar)-
  \vtinprod{\IntSec_\TruePar}{ \nabla^2_\Par \ZFun(\Sam_i, \TruePar)}
  \Big](\LOO{\EstPar}{\numobs-1}{i} - \EstPar_\numobs), \\
\Res_{6} & = \frac{\Numobs-1}{\Numobs} (\EstPar_\numobs -
\TruePar)^{\mytrans} [\int_{0}^1[ \nabla^2\TargetFun(\TruePar +
    t(\EstPar_\numobs - \TruePar)) - \nabla^2\TargetFun(\TruePar)]dt]
        [\sum_{i=1}^\Numobs(\LOO{\EstPar}{\numobs-1}{i} -
          \EstPar_\numobs)], \\
\Res_{7} & = - \frac{1}{\Numobs} (\EstPar_\numobs -
\TruePar)^{\mytrans} \vtinprod{\IntSec_\TruePar}{\sum_{j=1
  }^\Numobs\int_{0}^1[ \nabla^2_\Par \ZFun(\Sam_j, \TruePar +
    t(\EstPar_\numobs - \TruePar)) - \nabla^2_\Par \ZFun(\Sam_j,
    \TruePar)]dt}[\sum_{i=1}^\Numobs (\LOO{\EstPar}{\numobs-1}{i} -
  \EstPar_\numobs)], \\
\Res_{8} & = \frac{1}{\Numobs} \sum_{i=1}^\Numobs (\EstPar_\numobs -
\TruePar)^{\mytrans} \vtinprod{\IntSec_\TruePar}{\int_{0}^1[
    \nabla^2_\Par \ZFun(\Sam_i, \TruePar + t(\EstPar_\numobs -
    \TruePar)) - \nabla^2_\Par \ZFun(\Sam_i,
    \TruePar)]dt}(\LOO{\EstPar}{\numobs-1}{i} - \EstPar_\numobs).
\end{align*}

\paragraph{Bounds on $\Res_4, \Res_6, \Res_7$.}
By condition~\ref{ass:convergence} and
\Cref{LemLinErrorBound,LemJackknifeKey,lm:sup_op_norm_zfun}, we obtain
\begin{align*}
|\Res_4| & \leq \polyshort \sqrt{ \log \Numobs} \sqrt{
  \frac{\dimdelta}{\Numobs}} \cdot1\cdot \log^4(\Numobs/\delta)
\frac{\dimdelta}{\Numobs} \leq \polyshort \log^{4.5}(\Numobs/\delta)
\frac{\dimdelta^{3/2}}{\Numobs^{3/2}}, \\ |\Res_6| & \leq \polyshort
\sqrt{ \log \Numobs} \sqrt{ \frac{\dimdelta}{\Numobs}} \cdot1\cdot
\log^4(\Numobs/\delta) \frac{\dimdelta}{\Numobs} \leq \polyshort
\log^{4.5}(\Numobs/\delta) \frac{\dimdelta^{3/2}}{\Numobs^{3/2}},
\\ |\Res_7| & \leq \enorm{\EstPar_\numobs -
  \TruePar}^2\enorm{\sum_{i=1}^\Numobs(\LOO{\EstPar}{\numobs-1}{i} -
  \EstPar_\numobs)} \int_{0}^1 (1-s)s \sup_{\Par\in \ParSpace}
\opnorm{ \frac{1}{\Numobs} \sum_{i=1}^\Numobs[{ \nabla^3_\Par
      \ZFun(\Sam_i, \Par)[\IntSec_\TruePar]}} \textup{d}s \notag \\ &
  \leq \polyshort \sqrt{ \log (\Numobs/\delta)}
  \cdot\enorm{\EstPar_\numobs -
    \TruePar}^2\cdot\enorm{\sum_{i=1}^\Numobs(\LOO{\EstPar}{\numobs-1}{i}
    - \EstPar_\numobs)} \leq \polyshort \log^{5.5}(\Numobs/\delta)
  \frac{\dimdelta^2}{\Numobs^2} \\ & \leq \polyshort
  \log^{4}(\Numobs/\delta) \frac{\dimdelta^{3/2}}{\Numobs^{3/2}}.
\end{align*}
\paragraph{Bounds on $\Res_5, \Res_8$.}
For $\Res_5$ and $\Res_8$, we have
\begin{align*}
 &\quad|\Res_5 + \Res_8|\\ & = \Big| \frac{1}{\Numobs}
  \sum_{i=1}^\Numobs \Big[\vtinprod{\IntSec_\TruePar}{\int_{0}^1[
        \nabla^2_\Par \ZFun(\Sam_i, \TruePar + t(\EstPar_\numobs -
        \TruePar))]dt} - \nabla^2\TargetFun(\TruePar)
    \Big][\EstPar_\numobs - \TruePar, \LOO{\EstPar}{\numobs-1}{i} -
    \EstPar_\numobs] \Big|\\ & \leq \frac{1}{\Numobs}
  \sum_{i=1}^\Numobs \sup_{\Par\in \ParSpace} \Big| \nabla^2_\Par
  \ZFun(\Sam_i, \Par)[\IntSec_\TruePar, \EstPar_\numobs - \TruePar,
    \LOO{\EstPar}{\numobs-1}{i} - \EstPar_\numobs] \Big| +
  \enorm{\IntSec_\TruePar} \cdot\opnorm{ \nabla^2\TargetFun(\TruePar)}
  \cdot\enorm{\EstPar_\numobs - \TruePar}
  \cdot\enorm{\LOO{\EstPar}{\numobs-1}{i} - \EstPar_\numobs} \\ & \leq
  \frac{1}{\Numobs} \sum_{i=1}^\Numobs \sup_{\Par\in \ParSpace} \Big|
  \nabla^2_\Par \ZFun(\Sam_i, \Par)[\IntSec_\TruePar, \EstPar_\numobs
    - \TruePar, \LOO{\EstPar}{\numobs-1}{i} - \EstPar_\numobs] \Big| +
  \polyshort \log^{2.5}(\Numobs/\delta)
  \frac{\dimdelta}{\Numobs^{3/2}},
\end{align*}
where the last line follows from
Assumption~\ref{ass:convergence}~and~\Cref{LemJackknifeKey}.  For any
random variables $l_a^{(-i)},l_b^{(-i)}$ that are independent of
$\Sam_i$, we have~\cite[Exercise 2.18(a)]{wainwright2019high}
\begin{align*}
 &\quad\Prob(\sup_{\Par\in \ParSpace} \enorm{ \nabla^2_\Par
    \ZFun(\Sam_i, \Par)[l_a^{(-i)},l_b^{(-i)}]} \geq
  \sqrt{\ParDim}t_1\enorm{l_a^{(-i)}} \enorm{l_b^{(-i)}}) \\ & \leq
  \Prob(\sup_{\Par\in \ParSpace} \Big| \nabla^2_\Par \ZFun(\Sam_i,
  \Par) \Big[ \frac{l_a^{(-i)}}{\enorm{l_a^{(-i)}}},
    \frac{l_b^{(-i)}}{\enorm{l_b^{(-i)}}},e_k\Big] \Big|\geq t_1,
  \text{ for some }k\in[\ParDim]) \\ & \leq \ParDim\exp \big( - \big\{
  \frac{c t_1}{ \OrParc } \big\}^\OrConc \big).
\end{align*}
Choosing $t_1=c\OrParc \log^{\InvOrConc}(\Numobs\ParDim/\delta)$,
$l_a^{(-i)} = \IntSec_\TruePar$, $l_b^{(-i)} =
\LOO{\EstPar}{\numobs-1}{i} - \TruePar$, we have
\begin{align}
\sup_{\Par\in \ParSpace} \enorm{ \nabla^2_\Par \ZFun(\Sam_i,
  \Par)[\IntSec_\TruePar, \LOO{\EstPar}{\numobs-1}{i} - \TruePar]}
\leq c\OrParc \log^{\InvOrConc}(\Numobs\ParDim/\delta) \sqrt{\ParDim}
\enorm{\IntSec_\TruePar} \enorm{\LOO{\EstPar}{\numobs-1}{i} -
  \TruePar}\label{eq:r_5_and_r_8_pf_1}
\end{align}
for all $i\in[\Numobs]$ with probability at least $1 - \delta$.
Therefore, with probability at least $1 - \delta$, for all
$i\in[\Numobs]$
\begin{align*}
 &\quad\, \sup_{\Par\in \ParSpace} \Big| \nabla^2_\Par \ZFun(\Sam_i,
  \Par)[\IntSec_\TruePar, \EstPar_\numobs - \TruePar,
    \LOO{\EstPar}{\numobs-1}{i} - \EstPar_\numobs] \Big|\\ & \leq
  \sup_{\Par\in \ParSpace} \vecnorm{ \nabla^2_\Par \ZFun(\Sam_i,
    \Par)[\IntSec_\TruePar, \EstPar_\numobs - \TruePar]}{2}
  \vecnorm{\LOO{\EstPar}{\numobs-1}{i} - \EstPar_\numobs}{2} \\ &
  \overset{(i)}{\leq } \sup_{\Par\in \ParSpace} \opnorm{
    \nabla_\Par^2\ZFun(\Sam_i, \Par)[\IntSec_\TruePar]}
  \enorm{\EstPar_\numobs - \LOO{\EstPar}{\numobs-1}{i}}
  \vecnorm{\LOO{\EstPar}{\numobs-1}{i} - \EstPar_\numobs}{2} +
  c\OrParc \log^{\InvOrConc}(\Numobs/\delta) \sqrt{\ParDim}
  \enorm{\IntSec_\TruePar} \enorm{\LOO{\EstPar}{\numobs-1}{i} -
    \TruePar} \vecnorm{\LOO{\EstPar}{\numobs-1}{i} -
    \EstPar_\numobs}{2} \\ &\overset{(ii)}{\leq } \polyshort\Big[ \log
    ^{5.5}(\Numobs/\delta) \frac{\dimdelta^{2}}{\Numobs^{2}} + \log
    ^{4.5}(\Numobs/\delta) \frac{\dimdelta^{3/2}}{\Numobs^{3/2}} \Big]
  \leq \polyshort \log^{4.5}(\Numobs/\delta)
  \frac{\dimdelta^{3/2}}{\Numobs^{3/2}},
\end{align*}
where equation (i) uses equation~\eqref{eq:r_5_and_r_8_pf_1}; equation
(ii) follows from
Assumption~\ref{ass:convergence},~\Cref{LemJackknifeKey} and the fact
that
\begin{align}
\sup_{\Par\in \ParSpace} \opnorm{ \nabla_\Par^2\ZFun(\Sam_i,
  \Par)[\IntSec_\TruePar]} & \leq \sup_{\Par\in \ParSpace} \fronorm{
  \nabla_\Par^2\ZFun(\Sam_i, \Par)[\IntSec_\TruePar]} \leq
\sqrt{\sum_{i,j=1}^\ParDim \sup_{\Par\in \ParSpace}|
  \nabla_\Par^2\ZFun(\Sam_i, \Par)[\IntSec_\TruePar,e_i, e_j]|^2}
\notag \\ & \leq \polyshort\ParDim
\log^{\InvOrConc}(2\Numobs\ParDim^2/\delta) \leq \polyshort\ParDim
\log^{3/2}(\Numobs/\delta). \label{eq:hessian_times_vector_op_norm_concen}
\end{align}
Here the last line uses Assumption~\ref{ass:tail} and the tail bound
for Orlicz-norm bounded random variable~\cite[Exercise
  2.18(a)]{wainwright2019high}. Therefore, we obtain the bound
\begin{align*}
|\Res_5 + \Res_8| & \leq \polyshort \log^{4.5}(\Numobs/\delta)
\frac{\dimdelta^{3/2}}{\Numobs^{3/2}}.
\end{align*}
Finally, combining the bounds for the term $\{ \Res_i \}_{i=4}^8$
yields~\cref{lm:compare_term_2_auxi_2}.


\subsection{Proof of bound~\eqref{EqnHardThree} from~\Cref{LemHard}}
\label{AppLemHardThree}

In order to prove this bound, we first decompose $\STerm_3$ as the sum
$\STerm_3 = \STerm_{31} + \STerm_{32}$, where
\begin{align*}
\STerm_{31} & \defn \frac{1}{2}
\sum_{i=1}^\Numobs(\LOO{\EstPar}{\numobs-1}{i} -
\EstPar_\numobs)^{\mytrans} \Big[ \frac{\Numobs-1}{\Numobs}
  \nabla^2\TargetFun(\TruePar) - \frac{1}{\Numobs} \sum_{j \neq i}
  \vtinprod{\IntSec_\TruePar}{ \nabla^2_\Par \ZFun(\Sam_j, \TruePar)}
  \Big] (\LOO{\EstPar}{\numobs-1}{i} - \EstPar_\numobs), \quad
\mbox{and} \\
\STerm_{32} & \defn \frac{1}{2}
\sum_{i=1}^\Numobs(\LOO{\EstPar}{\numobs-1}{i} -
\EstPar_\numobs)^{\mytrans} \Big[ \frac{\Numobs-1}{\Numobs} ( \nabla^2
  \TargetFun(\EstPar_\numobs) - \nabla^2\TargetFun(\TruePar)) -
  \frac{1}{\Numobs} \sum_{j \neq i} \vtinprod{\IntSec_\TruePar}{
    \nabla^2_\Par \ZFun(\Sam_j, \EstPar_\numobs) - \nabla^2_\Par
    \ZFun(\Sam_j, \TruePar)} \Big] (\LOO{\EstPar}{\numobs-1}{i} -
\EstPar_\numobs).
\end{align*}
We claim that, with probability least $1 - \delta$, the following
bounds hold:
\begin{subequations}
\begin{align}
\label{lm:compare_term_3_auxi_1}   
|\Term_3 - \STerm_{31}| & \leq \polyshort \log^5(\Numobs/\delta)
\frac{\dimdelta^{3/2}}{\Numobs^{3/2}}, \quad \mbox{and} \\
\label{lm:compare_term_3_auxi_2}
|\STerm_{32}| & \leq \polyshort \log^{5}(\Numobs/\delta)
\frac{\dimdelta^{3/2}}{\Numobs^{3/2}}.
\end{align} 
\end{subequations}
We prove equation
~\eqref{lm:compare_term_3_auxi_1}~and~\eqref{lm:compare_term_3_auxi_2}
in~\Cref{sec:proof_lm:compare_term_3_auxi_1}
and~\Cref{sec:proof_lm:compare_term_3_auxi_2},
respectively. Combining~\cref{lm:compare_term_3_auxi_1}
and~\cref{lm:compare_term_3_auxi_2} yields~\cref{EqnHardThree}.


\subsubsection{Proof of~\cref{lm:compare_term_3_auxi_1}}
\label{sec:proof_lm:compare_term_3_auxi_1}

Define $\Hesdif_i \defn \nabla^2\TargetFun(\TruePar) -
\vtinprod{\IntSec_\TruePar}{ \nabla^2_\Par \ZFun(\Sam_i, \TruePar)}$.
By the definition of $\Term_{3}$ and $\STerm_{31}$, we have
\begin{align*}
&\Term_{3} - \STerm_{31} \\ & = \frac{1}{2} (\EstPar_\numobs -
  \TruePar)^{\mytrans} \Mhat (\EstPar_\numobs - \TruePar) -
  \frac{1}{2} \sum_{i=1}^\Numobs(\LOO{\EstPar}{\numobs-1}{i} -
  \EstPar_\numobs)^{\mytrans} \Mhat (\LOO{\EstPar}{\numobs-1}{i} -
  \EstPar_\numobs) + \frac{1}{2\Numobs}
  \sum_{i=1}^\Numobs(\LOO{\EstPar}{\numobs-1}{i} -
  \EstPar_\numobs)^{\mytrans} \Hesdif_i (\LOO{\EstPar}{\numobs-1}{i} -
  \EstPar_\numobs) \\
 & = \frac{1}{2} [\Errfirst + \frac{1}{\Numobs} \sum_{i=1}^\Numobs
    \sitwo]^{\mytrans} \Mhat [\Errfirst + \frac{1}{\Numobs}
    \sum_{i=1}^\Numobs \sitwo] - \frac{1}{2}
  \sum_{i=1}^\Numobs(\LOO{\EstPar}{\numobs-1}{i} -
  \EstPar_\numobs)^{\mytrans} \Mhat (\LOO{\EstPar}{\numobs-1}{i} -
  \EstPar_\numobs) \\
 &\quad + \frac{1}{2\Numobs}
  \sum_{i=1}^\Numobs(\LOO{\EstPar}{\numobs-1}{i} -
  \EstPar_\numobs)^{\mytrans} \Hesdif_i (\LOO{\EstPar}{\numobs-1}{i} -
  \EstPar_\numobs) \\
& = \underbrace{ \frac{1}{\Numobs} \sum_{i=1}^\Numobs
    \sitwo^{\mytrans} \Mhat \Errfirst + \frac{1}{2}
    \Errfirst^{\mytrans} \Mhat\Errfirst}_{\Res_{31}} + \underbrace{
    \frac{1}{2\Numobs^2} \sum_{i\neq j} \sitwo^{\mytrans} \Mhat
    \sitwobase_j}_{\Res_{32}} \\
 &\quad + \underbrace{ \frac{1}{2} \sum_{i=1}^\Numobs\Big[
      \frac{\sitwo^{\mytrans}}{\Numobs} \Mhat \frac{\sitwo}{\Numobs} -
      (\LOO{\EstPar}{\numobs-1}{i} - \EstPar_\numobs)^{\mytrans} \Mhat
      (\LOO{\EstPar}{\numobs-1}{i} - \EstPar_\numobs)
      \Big]}_{\Res_{33}} + \underbrace{ \frac{1}{2\Numobs}
    \sum_{i=1}^\Numobs(\LOO{\EstPar}{\numobs-1}{i} -
    \EstPar_\numobs)^{\mytrans} \Hesdif_i (\LOO{\EstPar}{\numobs-1}{i}
    - \EstPar_\numobs)}_{\Res_{34}}.
\end{align*}
It remains to control the residuals $\{ \Res_{3i} \}_{i=1}^4$.


\paragraph{Bound on $\Res_{31}$.}

By Assumption~\ref{ass:convergence}~and~\Cref{LemLinErrorBound}, with
probability at least $1 - \delta$, we have
\begin{align}
\label{eqn:concentration-of-si2}
\enorm{ \frac{1}{\Numobs} \sum_{i=1}^\Numobs \sitwo} \leq
\enorm{\EstPar_\numobs - \TruePar} + \enorm{\Errfirst} \leq \polyshort
\log \Numobs \sqrt{ \frac{\dimdelta}{\Numobs}}.
\end{align}
From~\cref{EqnMBound}, we have $\opnorm{\Mhat} \leq \polyshort$ with
probability at least $1 - \delta$.  Therefore, we have established
that
\begin{align*}
|\Res_{31}| & \leq \enorm{ \frac{1}{\Numobs} \sum_{i=1}^\Numobs
  \sitwo} \opnorm{\Mhat} \enorm{\Errfirst} + \frac{1}{2}
\enorm{\Errfirst} \opnorm{\Mhat} \enorm{\Errfirst} \\ & \leq
\polyshort( \log^{2.5} \Numobs \frac{\dimdelta^{3/2}}{\Numobs^{3/2}} +
\log^{2} \Numobs \frac{\dimdelta^{2}}{\Numobs^{2}}) \leq \polyshort
\log^{2.5} \Numobs \frac{\dimdelta^{3/2}}{\Numobs^{3/2}}.
\end{align*}

\paragraph{Bound on $\Res_{33}$.} By the definitin of $\Res_{33}$,
with probability at least $1 - \delta$, we have
\begin{align*}
|\Res_{33}| & = \Big| \frac{1}{2} \sum_{i=1}^\Numobs\Big[(
  \frac{\sitwo}{\Numobs}-[\LOO{\EstPar}{\numobs-1}{i} -
    \EstPar_\numobs])^{\mytrans} \Mhat( \frac{\sitwo}{\Numobs} +
       [\LOO{\EstPar}{\numobs-1}{i} - \EstPar_\numobs]) \Big]
\Big|\\ & \leq \sum_{i=1}^\Numobs\enorm{
  \frac{\sitwo}{\Numobs}-[\LOO{\EstPar}{\numobs-1}{i} -
    \EstPar_\numobs]} \opnorm{\Mhat} \enorm{ \frac{\sitwo}{\Numobs} +
  [\LOO{\EstPar}{\numobs-1}{i} - \EstPar_\numobs]} \\ & \leq
\Numobs\polyshort \frac{{\dimdelta}
  \log^3(\Numobs/\delta)}{\Numobs^{3/2}} \cdot \frac{ \sqrt{\dimdelta}
  \log^2(\Numobs/\delta)}{\Numobs} \\ & \leq \polyshort \log
^{5}(\Numobs/\delta) \frac{\dimdelta^{3/2}}{\Numobs^{3/2}}.
\end{align*}
\paragraph{Bound on $\Res_{34}$.}
Using~\cref{eq:hessian_times_vector_op_norm_concen} from the proof of
bound~\eqref{EqnHardTwo}, we establish
\begin{align}
\opnorm{\Hesdif_i} & \leq \opnorm{ \nabla^2\TargetFun(\TruePar)} +
\opnorm{\vtinprod{\IntSec_\TruePar}{ \nabla^2_\Par \ZFun(\Sam_i,
    \TruePar)}} \leq \opnorm{ \nabla^2\TargetFun(\TruePar)} +
\sup_{\Par\in \ParSpace} \opnorm{\vtinprod{\IntSec_\TruePar}{
    \nabla^2_\Par \ZFun(\Sam_i, \TruePar)}} \notag \\ & \leq
\polyshort + \ParDim\polyshort
\log^{\InvOrConc}(2\Numobs\ParDim^2/\delta) \leq \polyshort\ParDim
\log^{3/2}(\Numobs/\delta)
\label{eq:claim:bound_hesdif}
\end{align}
for all $i\in[\Numobs]$ with probability at least $1 -
\delta$. Combining~\cref{eq:claim:bound_hesdif}
with~\Cref{LemJackknifeKey}, we obtain
\begin{align*}
|\Res_{34}| & \leq \frac{1}{2\Numobs} \sum_{i=1}^\Numobs
\enorm{\LOO{\EstPar}{\numobs-1}{i} - \EstPar_\numobs}^2
\opnorm{\Hesdif_i} \leq \polyshort \log^{5.5}(\Numobs/\delta)
\frac{\dimdelta^2}{\Numobs^2}.
\end{align*}
\paragraph{Bound on $\Res_{32}$.}
For $\Res_{32}$, we claim that with probability at least $1 - \delta$,
we have
\begin{align}
|\Res_{32}| \leq \polyshort
\Big(\frac{\sqrt\ParDim}{\numobs}+\frac{\ParDim^{3/2}}{\Numobs^{3/2}}
\Big) \log^{5}(\Numobs/\delta).\label{eq:claim:bound_res_32}
\end{align}
Combining the bounds for $\Res_{3i}(i=1,2,3,4)$
yields~\Cref{lm:compare_term_3_auxi_1}.  \\

\myunder{Proof of Claim~\eqref{eq:claim:bound_res_32}:} We write
$\Res_{32} = \Res_{35} + \Res_{36}+\Res_{37} $, where
\begin{align*}
\Res_{35} \defn \frac{1}{2\Numobs^2} \sum_{i\neq j} \sitwo^{\mytrans}
\Mmat \sitwobase_j, ~~ \Res_{36}\defn \frac{1}{2\Numobs^2} \sum_{i, j
  =1}^\numobs \sitwo^{\mytrans}(\Mhat - \Hesdif) \sitwobase_j,~~
\Res_{37}\defn -\frac{1}{2\Numobs^2} \sum_{i =1}^\numobs
\sitwo^{\mytrans}(\Mhat - \Hesdif) \sitwobase_i.
\end{align*}
Similar to Lemma~\ref{LemUstatBound}, with probability at least $1 -
\delta$, we have
\begin{align*}
|\Res_{35}| & \leq \polyshort \frac{ \sqrt{\ParDim}}{\Numobs}
\log^{5}(\Numobs/\delta).
\end{align*}
Moreover, since $\Res_{36}=\RemainLthree$ in
equation~\eqref{EqnDefnRemainL}, we have by
eqaution~\eqref{EqnRemainLBound} that
\begin{align*}
|\Res_{35}| & \leq \polyshort \Big(\frac{\dimdelta}{\numobs}
\Big)^{3/2} \log^2(\numobs/\delta)
\end{align*}
with probability at least $1 - \delta$.  Finally,
\begin{align*}
|\Res_{37}| & \leq \frac{1}{2\Numobs^2} \sum_{i =1}^\numobs
\opnorm{\Mhat - \Hesdif} \,
\vecnorm{\sitwo}{2}^2\\ &\overset{(i)}{\leq}
\frac{\polyshort}{\numobs}
\sup_{i\in[\numobs]}\vecnorm{\ZFun(\Sam_i,\TruePar)}{2}^2 \cdot
\opnorm{\Mhat - \Hesdif} \\ &\overset{(ii)}{\leq}
\polyshort\frac{\ParDim}{\numobs} \log(\numobs/\delta)\cdot \Big(
\sqrt{ \frac{\Ccerr\log(\numobs/\delta)}{\Numobs}} + \frac{\dimdelta
  \log^{5/2} (\Numobs/\delta)}{\Numobs} \Big)\leq \polyshort
\Big(\frac{\dimdelta}{\numobs} \Big)^{3/2} \log^2(\numobs/\delta)
\end{align*}
where step~(i) uses the definition of $\sitwo$ in
equation~\eqref{EqnSidefn}; step~(ii) follows from
equation~\eqref{EqnMBound} and the property of sub-Exponential vectors
that
\begin{align*}
   \vecnorm{\ZFun(\Sam_i,\TruePar)}{2}^2
   =\sum_{j=1}^\ParDim\sup_{i\in[\numobs]}|\inprod{\ZFun(\Sam_i,\TruePar)}{e_j}|^2
   \leq \polyshort\ParDim \log(\numobs\ParDim/\delta) \leq
   \polyshort\ParDim \log(\numobs/\delta)
\end{align*}
with probability at least $1-\delta.$ Putting the bounds on
$\Res_{35}, \Res_{36}, \Res_{37}$ together
yields~\cref{eq:claim:bound_res_32}.


\subsubsection{Proof of~\Cref{lm:compare_term_3_auxi_2}}
\label{sec:proof_lm:compare_term_3_auxi_2}

By
Assumption~\ref{ass:convergence},~\Cref{lm:sup_op_norm_zfun,LemJackknifeKey},
with probability at least $1 - \delta$, we have
\begin{align*}
|\STerm_{32}| & \leq \LipTarc \sum_{i=1}^\Numobs
\enorm{\LOO{\EstPar}{\numobs-1}{i} - \EstPar_\numobs}^2
\enorm{\EstPar_\numobs - \TruePar} +
\sum_{i=1}^\Numobs\enorm{\LOO{\EstPar}{\numobs-1}{i} -
  \EstPar_\numobs}^2 \opnorm{ \frac{1}{\Numobs} \sum_{j \neq i}
  \vtinprod{\IntSec_\TruePar}{ \nabla^2_\Par \ZFun(\Sam_j,
    \EstPar_\numobs) - \nabla^2_\Par \ZFun(\Sam_j, \TruePar)}} \\
 & \leq \LipTarc\sum_{i=1}^\Numobs \enorm{\LOO{\EstPar}{\numobs-1}{i}
  - \EstPar_\numobs}^2\enorm{\EstPar_\numobs - \TruePar} +
\sum_{i=1}^\Numobs\enorm{\LOO{\EstPar}{\numobs-1}{i} -
  \EstPar_\numobs}^2\enorm{\EstPar_\numobs - \TruePar}
\sup_{\Par\in \ParSpace} \opnorm{ \frac{1}{\Numobs} \sum_{j \neq i} {
    \nabla^3_\Par \ZFun(\Sam_j, \Par)[\IntSec_\TruePar]}} \\
& \leq \polyshort \frac{\dimdelta^{3/2}}{\Numobs^{3/2}}
\log^{5}(\Numobs/\delta)
\end{align*}
for some constant $\polyshort = \poly > 0$.


\section{LOO difference control: Proof of ~\Cref{LemJackknifeKey}}
\label{sec:proof_LemJackknifeKey}

In this appendix, we prove~\Cref{LemJackknifeKey} that bounds various
functions of LOO differences $\LooDiff = \LOO{\EstPar}{\numobs-1}{i} -
\EstPar_\numobs$. We prove the bounds~\eqref{eq:LemJackknifeKey_eq1},~\eqref{eq:LemJackknifeKey_eq3},~\eqref{eq:LemJackknifeKey_eq4} in Section~\ref{sec:pf_eq:LemJackknifeKey_eq1},~\ref{sec:pf_eq:LemJackknifeKey_eq3},~\ref{sec:pf_eq:LemJackknifeKey_eq4}, respectively.

\subsection{Proof of bound~\eqref{eq:LemJackknifeKey_eq1}}
\label{sec:pf_eq:LemJackknifeKey_eq1}

Consider the inequalities
\begin{align}
  \label{eq:LemJackknifeKey}
\enorm{\LooDiff} \; \stackrel{(a)}{\leq} \; \frac{c}{\strongconvex
  \Numobs} \enorm{\ZFun(\Sam_i, \EstPar_\numobs)} \quad \mbox{and}
\quad \enorm{\ZFun(\Sam_i, \EstPar_\numobs)} \stackrel{(b)}{\leq}
\polyshort { \sqrt{\dimdelta} \log^2(\Numobs/\delta)}.
\end{align}
By the union bound, it suffices that each bound holds simultaneouly
for all $i \in [\numobs]$ with probability at least $1-\delta/2$.
Sections~\ref{sec:pf_eq:LemJackknifeKey_eq11}~and~\ref{sec:pf_eq:LemJackknifeKey_eq12},
respectively, are devoted the proofs of these two claims.


\subsubsection{Proof of bound~\eqref{eq:LemJackknifeKey}(a)}
\label{sec:pf_eq:LemJackknifeKey_eq11}

By~\cref{eq:z_estimate_1} and the Newton-Leibniz formula, we have
\begin{align*}
\ZFun(\Sam_i, \EstPar_\numobs )= \sum_{j \neq i}(\ZFun(\Sam_j,
\LOO{\EstPar}{\numobs-1}{i}) - \ZFun(\Sam_j, \EstPar_\numobs))=
\sum_{j \neq i} \Big[ \int_0^1 \nabla _{\Par} \ZFun \big(\Sam_j,
  \EstPar_\numobs + t \LooDiff \big ) dt \Big] \LooDiff,
\end{align*}
and thus the relation
\begin{align}
\label{eq:thetai-theta_1}
\LooDiff = \frac{1}{\Numobs - 1} \Biggr \{ \frac{1}{\Numobs - 1}
\sum_{j \neq i} \int_0^1 \nabla _{\Par} \ZFun(\Sam_j, \EstPar_\numobs
+ t \LooDiff ) dt \Biggr \}^{-1} \ZFun(\Sam_i, \EstPar_\numobs).
\end{align}
Consequently, it suffices to lower bound the minimum singular value;
in particular, we claim that there is a constant $c > 0$
such that
\begin{align}
\label{eq:pf_jackknifekey_111}  
\sigma_{\min} \Big(\frac{1}{\Numobs - 1} \sum_{j \neq i} \int_0^1
\nabla _{\Par} \ZFun(\Sam_j, \EstPar_\numobs + t \LooDiff ) dt \Big)
\geq c \strongconvex
\end{align}
for all $i\in[\numobs]$ with probability at least $1 -
\tfrac{\delta}{2}$.  Introducing the shorthand $\DelHat \defn
\EstPar_\numobs - \TruePar$, we have
\begin{multline*}
\opnorm{\sum_{j \neq i} \int_0^1 [ \nabla_\Par \ZFun(\Sam_j,
    \EstPar_\numobs + t \LooDiff ) - \nabla_{\Par} \ZFun(\Sam_j,
    \TruePar)] dt} \leq \sup_{\Par \in \ParSpace} \opnorm{\sum_{j \neq
    i} \nabla^2_{\Par} \ZFun(\Sam_j, \Par)} \; \big \{
\enorm{\LooDiff} + \enorm{\DelHat} \big \},
\end{multline*}
Now define the matrices
\begin{align}
\label{eq:pf_jackknifekey_112}   
  \Bmat_{\Numobs,i} & \defn \GradFun_{\TruePar} - \frac{1}{\Numobs -
    1} \sum_{j \neq i} \int_0^1 \nabla_{\Par} \ZFun(\Sam_j,
  \EstPar_\numobs + t(\LOO{\EstPar}{\numobs-1}{i} - \EstPar_\numobs))
  dt \qquad \mbox{for $i = 1, \ldots, \numobs$.}
\end{align}
Applying~\Cref{lm:sup_op_norm_zfun,lm:op_norm_zfun} and the triangle
inequality yields, for each $i \in [\numobs]$, the upper bound
\begin{align}
\label{eq:est_ani_thetai-theta}
\opnorm{\Bmat_{\Numobs,i}} & \leq c\OrParc\Big(1 + \sqrt{
  \frac{\dimdelta}{\Numobs}} + \frac{(\dimdelta \log
  {\Numobs})^{\InvOrConc}}{\Numobs}
\Big)(\enorm{\LOO{\EstPar}{\numobs-1}{i} - \TruePar} +
\enorm{\EstPar_\numobs - \TruePar}) + c\OrParb\Big( \sqrt{
  \frac{\dimdelta}{\Numobs}} + \frac{(\dimdelta \log
  {\Numobs})}{\Numobs} \Big) \notag \\ & \overset{(i)}{\leq}
c\OrParc(\enorm{\LOO{\EstPar}{\numobs-1}{i} - \TruePar} +
\enorm{\EstPar_\numobs - \TruePar}) + c\OrParb \sqrt{
  \frac{\dimdelta}{\Numobs}} \notag \\ & \leq \polyshort \sqrt{ \log
  \Numobs} \sqrt{ \frac{\ParDim + \log (\numobs/
    \delta)}{\Numobs}}\overset{(ii)}{\leq} \strongconvex/2
\end{align}
valid with probability at least $1 - \delta/2$.  Here in
step~(i)~and~(ii), we made use of the assumed lower bound on the
sample size~\eqref{eq:linear_approx_Z_est_sample_ahead}. Combining
equation~\eqref{eq:pf_jackknifekey_112}~and~\eqref{eq:est_ani_thetai-theta}
with Assumption~\ref{ass:convergence} and applying triangle inequality
yields the lower bound~\eqref{eq:pf_jackknifekey_111}.


\subsubsection{Proof of bound~\eqref{eq:LemJackknifeKey}(b)}
\label{sec:pf_eq:LemJackknifeKey_eq12}

We introduce the shorthand $\LinApp^{(-i)} \defn \Numobs^{-1} \sum_{j
  \neq i} \sitwowithout_j$, and write $\EstPar_\numobs - \TruePar=
\LinApp^{(-i)} + \Errfirst + \Numobs^{-1} \sitwo$.  We also recall the
definitions of $\Errfirst$ and $\sitwowithout_j$ from
equations~\eqref{EqnDefnLinError} and~\eqref{eq:def_psi},
respectively, and observe that $\LinApp^{(-i)}$ is
independent of $\Sam_i$.

With these observations in place, performing a Taylor expansion yields
\begin{align}
\label{eq:lm_thetai-theta_lipschitz}
\enorm{\ZFun(\Sam_i, \EstPar_\numobs)} & \leq \enorm{\ZFun(\Sam_i,
  \TruePar)} + \enorm{ \nabla _\Par\ZFun(\Sam_i,
  \TruePar)(\EstPar_\numobs - \TruePar)}+ \sup_{\Par\in \ParSpace}
\sqrt{\sum_{j=1}^\ParDim \big[\nabla ^2_\Par\ZFun(\Sam_i, \Par)[e_j,
      \EstPar_\numobs - \TruePar, \EstPar_\numobs - \TruePar]\big]^2}
\notag \\
&\leq
\LooBound_{a1} + \LooBound_{a2}+\LooBound_{a3},
 \end{align}
where
\begin{align*}
 \LooBound_{a1}
 &\defn
 \enorm{\ZFun(\Sam_i, \TruePar)},\\
 \LooBound_{a2}
 &\defn
 \enorm{ \nabla
  _\Par\ZFun(\Sam_i, \TruePar) \LinApp^{(-i)}} +
\sup_{\Par\in \ParSpace} \opnorm{ \nabla _\Par\ZFun(\Sam_i,
  \Par)}(\enorm{ \frac{1}{\Numobs} \sitwo} + \enorm{\Errfirst}),\\
  \LooBound_{a3}
  &\defn
  2\sup_{\Par\in \ParSpace} \opnorm{ \nabla
  ^2_\Par\ZFun(\Sam_i, \Par)}(\enorm{ \frac{1}{\Numobs} \sitwo} +
\enorm{\Errfirst})^2+2 \sup_{\Par\in \ParSpace}
\sqrt{\sum_{j=1}^\ParDim \big[\nabla ^2_\Par\ZFun(\Sam_i, \Par)[e_j,
      \LinApp^{(-i)}, \LinApp^{(-i)}]\big]^2}
\end{align*} are the terms corresponding to the zeroth to second order gradients of $\ZFun$. 
 
It remains to provide bounds for each term  and combine them. \\

\myunder{Intermediate bounds: }
we first collect some intermediate results used to derive bounds for $\LooBound_{a1},\LooBound_{a2},\LooBound_{a3}$.
Letting $t_1=c\OrPara \log (\ParDim\Numobs /\delta)$ and using
Assumption~\ref{ass:tail}, we have
\begin{align}
\Prob(\LooBound_{a1} \geq \sqrt{\ParDim}t_1) & \leq
\Prob(|\ZFun_j(\Sam_i, \TruePar)|\geq t_1,~\text{for some
}j\in[\ParDim]) \leq 2\exp\{- \frac{t_1}{c\OrPara} + \log \ParDim\}
\leq \frac{\delta}{\Numobs} \label{eq:union_gaussian_eq1_add}.
\end{align}

For any $i$, by~\cref{eqn:concentration-of-si2} and a union bound, we have 
\begin{align} \label{eq:bound-on-l-minus-i}
\enorm{\LinApp^{(-i)}} \leq \polyshort \log \numobs \sqrt{
  \frac{\dimdeltatil}{n}}
\end{align} with probability at least $1-\delta/\numobs$, where $\dimdeltatil\defn  \ParDim+\log(\numobs/\delta)$.


Let   $t_2=c\OrParb \log (\ParDim\Numobs /\delta)$ and $t_3=c\OrParc  \log^{\InvOrConc} (\ParDim\Numobs /\delta) $. Since $\LinApp^{(-i)}$ is independent of $\Sam_i$, we have  
\begin{subequations}
\begin{align}
&\Prob(\enorm{ \nabla _\Par\ZFun(\Sam_i, \TruePar) \LinApp^{(-i)}} \geq
  \sqrt{\ParDim}t_2\enorm{\LinApp^{(-i)}}) \notag \\ &\qquad\leq
  \Prob(| \nabla _\Par\ZFun(\Sam_i, \TruePar)[e_j, \LinApp^{(-i)}]| \geq
  t_2\enorm{\LinApp^{(-i)}},~\text{for some }j\in[\ParDim]) \leq
  2\exp\{- \frac{ct_2}{\OrParb} + \log \ParDim\} \le
  \frac{\delta}{\Numobs} .\label{eq:union_gaussian_eq2} \\ &\Prob\Bigg(
  \sup_{\Par\in \ParSpace} \sqrt{\sum_{j=1}^\ParDim
    \big[\nabla ^2_\Par\ZFun(\Sam_i,
      \Par)[e_j, \LinApp^{(-i)}, \LinApp^{(-i)}]\big]^2} \geq
  \sqrt{\ParDim}t_3\enorm{\LinApp^{(-i)}}^2\Bigg) \notag \\ &\qquad\leq
  \Prob( \sup_{\Par\in \ParSpace}| \nabla ^2_\Par\ZFun(\Sam_i,
  \TruePar)[e_j, \LinApp^{(-i)}, \LinApp^{(-i)}]| \geq
  t_3\enorm{\LinApp^{(-i)}}^2,~\text{for some }j\in[\ParDim]) \leq
  2\exp\{- \frac{ct^\OrConc_3}{\OrParc^\OrConc} + \log \ParDim\}
  \notag \\ &\qquad\leq \frac{\delta}{\Numobs}
  .\label{eq:union_gaussian_eq2_add}
\end{align}

\end{subequations}
Moreover, from~\Cref{lm:sup_op_norm_zfun} with $\numobs=1$, we have
\begin{subequations}
\begin{align}
\sup_{\Par\in \ParSpace} \opnorm{ \nabla _\Par\ZFun(\Sam_i, \Par)} & \leq
\polyshort\dimdeltatil\log^2(\numobs/\delta), \label{eq:union_gaussian_eq3}
\\ \sup_{\Par\in \ParSpace} \opnorm{ \nabla ^2_\Par\ZFun(\Sam_i, \Par)} &
\leq \polyshort\dimdeltatil^{\InvOrConc}
, \label{eq:union_gaussian_eq3_add}
\end{align} 
for all $i\in[\numobs]$ with probability at least $1-\delta.$ 
\end{subequations}
Lastly, by Assumption~\ref{ass:tail} and the concentration properties of sub-Exponential variables, we have 
\begin{align}
  \enorm{\sitwobase_i} \leq \sqrt\ParDim \log(\numobs/\delta) \label{eq:single_psi_bound}
\end{align}
for all $i\in[\numobs]$ with probability at least $1-\delta.$\\

\myunder{Combining the bounds:}
\begin{subequations}  
 we proceed to combine the intermediate bounds to obtain bounds for $\LooBound_{a1},\LooBound_{a2},\LooBound_{a3}$. 
From equation~\eqref{eq:union_gaussian_eq1_add}, we have 
\begin{align}
  \LooBound_{a1}\leq c\OrPara \sqrt{\ParDim} \log (\Numobs\ParDim/\delta)
    \label{eq:loobound_a1_bound}
\end{align}
with probability at least $1-\delta$.
From~\cref{eq:single_psi_bound,eq:bound-on-l-minus-i,eq:union_gaussian_eq2,eq:union_gaussian_eq3} and Lemma~\ref{LemLinErrorBound}, we have
\begin{align}
  \LooBound_{a2}
  &\leq 
  c\OrParb
   \sqrt{\ParDim} \log (\ParDim\Numobs /\delta) \enorm{\LinApp^{(-i)}}
   + \polyshort\log^2(\numobs/\delta) \dimdeltatil (\enorm{
     \frac{1}{\Numobs} \sitwo} + \enorm{\Errfirst})
     \notag\\
     &\leq 
       \polyshort \log^2(\numobs/\delta)\sqrt{\frac{\ParDim\dimdeltatil}{\numobs}}
   + \polyshort\log^{3.5}(\numobs/\delta) \sqrt{\dimdelta}\frac{\dimdeltatil}{\numobs}
    \label{eq:loobound_a2_bound}
\end{align} with probability at least $1-\delta$.
Lastly, from~\cref{eq:single_psi_bound,eq:bound-on-l-minus-i,eq:union_gaussian_eq2_add,eq:union_gaussian_eq3_add} and Lemma~\ref{LemLinErrorBound},   we have
\begin{align}
  \LooBound_{a3}
  &\leq 
 \polyshort\dimdelta^{\InvOrConc}(\enorm{ \frac{1}{\Numobs} \sitwo}
   + \enorm{\Errfirst})^2 + c\OrParc\sqrt{\ParDim}
   \log^{\InvOrConc}(\numobs/\delta) \enorm{\LinApp^{(-i)}}^2 
     \notag\\
     &\leq 
       \polyshort \log^2(\numobs/\delta)\sqrt{\frac{\dimdelta^{3.5}}{\numobs^2}}
   + \polyshort\log^3(\numobs/\delta) \frac{\dimdeltatil\dimdelta}{\numobs}
    \label{eq:loobound_a3_bound}
\end{align}
\end{subequations}
with probability at least $1-\delta$.  Putting the
bounds~\eqref{eq:loobound_a1_bound},~\eqref{eq:loobound_a2_bound},~\eqref{eq:loobound_a3_bound}
together and using the sample size
condition~\eqref{eq:linear_approx_Z_est_sample_ahead} yields
bound~\eqref{eq:LemJackknifeKey}(b).


\subsection{Proof of the bound~\eqref{eq:LemJackknifeKey_eq3}}
\label{sec:pf_eq:LemJackknifeKey_eq3}

We prove~\cref{eq:LemJackknifeKey_eq3} by bounding
\begin{align}
& \Sma_0 \defn \enorm{\LooDiff + \frac{1}{\Numobs-1}
    \GradFun_{\TruePar}^{-1} \ZFun(\Sam_i, \EstPar_\numobs)}, \quad
  \Sma_1\defn\enorm{ \frac{1}{\Numobs-1} \GradFun_{\TruePar}^{-1}
    \ZFun(\Sam_i, \TruePar)- \frac{1}{\Numobs-1}
    \GradFun_{\TruePar}^{-1} \ZFun(\Sam_i,
    \EstPar_\numobs)} \label{eq:thetai-theta-tri_2}
\end{align}
in equation~\eqref{eq:theta_i-theta_prec_1} and
\eqref{eq:theta_i-theta_prec3}, respectively. Equation~\eqref{eq:LemJackknifeKey_eq3} follows from the bounds and the triangle inequality.\\
 
\myunder{Bound on $\Sma_0$:} 
by~\cref{eq:est_ani_thetai-theta}~and~bound~\eqref{eq:LemJackknifeKey}(b) from the proof of bound~\eqref{eq:LemJackknifeKey_eq1} and
Woodbury's identity, we have
\begin{align}
\Sma_0 & = \frac{1}{\Numobs-1} \enorm{(\E[ \nabla _{\Par}
    \ZFun(\Sam_i, \TruePar)] +
  \Bmat_{\Numobs,i})^{-1}\Bmat_{\Numobs,i} \GradFun_{\TruePar}^{-1}
  \ZFun(\Sam_i, \EstPar_\numobs)} \notag \\
\label{eq:theta_i-theta_prec_1}
& \leq \polyshort\log^2(\Numobs/\delta) \log \Numobs
\frac{\dimdelta}{\Numobs^{3/2}}
\end{align}
with probability at least $1-\delta,$
where $\Bmat_{\Numobs,i}$ is defined in~\eqref{eq:pf_jackknifekey_112}.
\\

\myunder{Bound on $\Sma_1$:}  performing a Taylor expansion, we have
\begin{align}
& \qquad \enorm{\ZFun(\Sam_i, \EstPar_\numobs) - \ZFun(\Sam_i,
    \TruePar)} \notag \\
& \leq \sup_{\Par\in \ParSpace} \enorm{ \nabla _\Par\ZFun(\Sam_i,
    \Par)(\EstPar_\numobs - \LOO{\EstPar}{\numobs-1}{i})} + \enorm{
    \nabla _\Par\ZFun(\Sam_i, \TruePar)(\LOO{\EstPar}{\numobs-1}{i} -
    \TruePar)} \notag \\
& \qquad + \sup_{\Par\in \ParSpace} \sqrt{\sum_{j=1}^\ParDim
    \big[\nabla ^2_\Par\ZFun(\Sam_i,
      \Par)[e_j, \LOO{\EstPar}{\numobs-1}{i} -
        \TruePar, \LOO{\EstPar}{\numobs-1}{i}- \TruePar]\big]^2}
  \notag \\
& \leq \sup_{\Par\in \ParSpace} \opnorm{ \nabla _\Par\ZFun(\Sam_i,
    \Par)} \enorm{\EstPar_\numobs - \LOO{\EstPar}{\numobs-1}{i}} +
  c\OrParb \log (\ParDim\Numobs/\delta)
  \sqrt{\ParDim} \enorm{\LOO{\EstPar}{\numobs-1}{i} -
    \TruePar}+c\OrParc\log^{\InvOrConc}(\numobs/\delta) \sqrt{\ParDim} \enorm{\LOO{\EstPar}{\numobs-1}{i}
    - \TruePar}^2 
\label{eq:theta_i-theta_prec2}  
\end{align}
  for all $i\in[\numobs]$ which probability at least $1 - \delta/2$,
where the second inequality follows
from~\cref{eq:union_gaussian_eq2,eq:union_gaussian_eq2_add} with
$\LinApp^{(-i)}$ replaced by $\LOO{\EstPar}{\numobs-1}{i} - \TruePar$.

Combining this with equation~\eqref{eq:union_gaussian_eq3},~\eqref{eq:LemJackknifeKey_eq1} and
Assumption~\ref{ass:convergence}, we further obtain
\begin{align*}
    & \qquad \enorm{\ZFun(\Sam_i, \EstPar_\numobs) - \ZFun(\Sam_i,
    \TruePar)} \\
& \leq \polyshort \Big( \frac{\dimdelta^{1/2} \dimdeltatil
   \log^4(\Numobs/\delta) }{\Numobs} +
\log(\Numobs/\delta) \sqrt{\log\numobs}{ \frac{\sqrt{\ParDim}
    \sqrt{\dimdeltatil}}{ \sqrt\Numobs}} + \log^{\InvOrConc}
(\numobs/\delta)
\log\numobs\frac{\dimdeltatil\sqrt{\ParDim}}{{\numobs}} \Big) \\
& \leq \polyshort \log^2(\Numobs/\delta) \frac{\dimdelta}{ \sqrt\Numobs}.
\end{align*}
Therefore,
\begin{align}
\abss{\Sma_1} \leq \polyshort \log^2(\Numobs/\delta)
\frac{\dimdelta}{\Numobs^{3/2}} \label{eq:theta_i-theta_prec3}.
\end{align} 
  for all $i\in[\numobs]$ which probability at least $1 - \delta/2$.

\subsection{Proof of bound~\eqref{eq:LemJackknifeKey_eq4} from Lemma~\ref{LemJackknifeKey}}\label{sec:pf_eq:LemJackknifeKey_eq4}

Define the empirical estimate \mbox{$\EmpGZF \defn- {\Numobs^{-1}}
  \sum_{i=1}^\Numobs \nabla_\Par \ZFun(\Sam_i, \EstPar_\numobs)$} of
$\GradFun_{\TruePar}$, along with the auxiliary quantities
\begin{align}
\label{eq:def_js_1}  
\EmpGZFa & \defn - \frac{1}{\Numobs} \sum_{j=1}^\Numobs \int_0^1
\nabla_{\Par} \ZFun(\Sam_j, \EstPar_\numobs + t \LooDiff) dt, \quad
\mbox{and} \; \EmpGZFb \defn -\frac{1}{\Numobs} \sum_{j=1}^\Numobs
\nabla_\Par \ZFun(\Sam_j, \LOO{\EstPar}{\numobs-1}{i}).
\end{align}


\subsubsection{Main argument}

From equation~\eqref{eq:z_estimate_1}, we have the equivalence
\begin{equation*}
\ZFun(\Sam_i, \LOO{\EstPar}{\numobs-1}{i}) =
\sum_{j=1}^\Numobs(\ZFun(\Sam_j, \LOO{\EstPar}{\numobs-1}{i}) -
\ZFun(\Sam_j, \EstPar_\numobs))=-(\Numobs\GradFunminusi)
(\LOO{\EstPar}{\numobs-1}{i} - \EstPar_\numobs),
\end{equation*}
from which we can write 
\begin{align}\LOO{\EstPar}{\numobs-1}{i} - \EstPar_\numobs
= - \GradFunminusi^{-1} \ZFun(\Sam_i,
\LOO{\EstPar}{\numobs-1}{i})/\Numobs. \label{eq:thetai-theta_1_repr2}
\end{align}
Our analysis is based on the
upper bound
\begin{align}
\Sma & \leq \frac{1}{\Numobs} \enorm{\sum_{i=1}^\Numobs [ -
    \GradFunminusi ^{-1} + \GradFun_{\Numobs}^{-1}] \ZFun(\Sam_i,
  \LOO{\EstPar}{\numobs-1}{i})} + \frac{1}{\Numobs}
\enorm{\sum_{i=1}^\Numobs \GradFun_{\Numobs}^{-1} \ZFun(\Sam_i,
  \LOO{\EstPar}{\numobs-1}{i})} \notag \\
\label{eq:thetai-theta-eq4-pf1}  
& \leq \underbrace{\frac{1}{\Numobs} \sum_{i=1}^\Numobs
  \opnorm{\GradFunminusi^{-1}} \cdot \enorm{\Big( - \GradFunminusi +
    \GradFun_{\Numobs} \Big) \GradFun_{\Numobs}^{-1} \ZFun(\Sam_i,
    \LOO{\EstPar}{\numobs-1}{i})}}_{\Sma_2} +
\underbrace{\frac{1}{\Numobs} \enorm{\sum_{i=1}^\Numobs
    \GradFun_{\Numobs}^{-1} \ZFun(\Sam_i,
    \LOO{\EstPar}{\numobs-1}{i})}}_{\Sma_3}.
\end{align} 
In order to complete the proof of~\cref{eq:LemJackknifeKey_eq4}, it
suffices to show that
\begin{align}
\label{eq:sum_thetai-theta}
\Sma_2 \stackrel{(i)}{\leq} \polyshort \log^4(\Numobs/\delta)
\frac{\dimdelta}{\Numobs}, \quad \mbox{and} \quad
\Sma_3 \stackrel{(ii)}{\leq} \polyshort \log^4(\Numobs/\delta)
\frac{\dimdelta}{\Numobs}
\end{align}
with probability at least $1-\delta$.


\subsubsection{Proof of claim~\eqref{eq:sum_thetai-theta}(i)}
In the following, we first derive upper bounds for
$\opnorm{\GradFunminusi ^{-1}}, \opnorm{ \GradFun_{\Numobs}^{-1}},
\opnorm{ - \GradFunminusi + \GradFun_{\Numobs}}, \enorm{\ZFun(\Sam_i,
  \LOO{\EstPar}{\numobs-1}{i})} $. The
claim~\eqref{eq:sum_thetai-theta}(i) then follows directly from
combining these bounds.  \\

\myunder{Bounds on individual terms.}  We start with bounding the
difference $\opnorm{ - \GradFunminusi + \GradFun_{\Numobs}}$.
By~\Cref{lm:sup_op_norm_zfun}, equation~\eqref{eq:LemJackknifeKey_eq1}
and the sample size
condition~\cref{eq:linear_approx_Z_est_sample_ahead}, we have
\begin{align}
\opnorm{ - \GradFunminusi + \GradFun_{\Numobs}} & \leq
\sup_{\Par\in \ParSpace} \opnorm{ \frac{1}{\Numobs} \sum_{j=1}^\Numobs
  \nabla ^2_{\Par} \ZFun(\Sam_j, \Par)}
\enorm{\LOO{\EstPar}{\numobs-1}{i} - \EstPar_\numobs} \notag \\
\label{eq:thetai-theta-eq4}
& \stackrel{(i)}{\leq} \polyshort \frac{ \sqrt{\dimdelta}
  \log^2(\Numobs/\delta)}{\Numobs} \stackrel{(ii)}{\leq}
\strongconvex/4
\end{align}
for all $i \in [\Numobs]$ with probability at least $1 - \delta$.

Next, we prove a upper bound for $\opnorm{ \GradFun_{\Numobs}^{-1}}$
(or equivalently, a lower bound for
$\sigma_{\min}(\GradFun_\Numobs)$). A similar upper bound for
$\opnorm{\GradFunminusi^{-1}}$ can be established through the triangle
inequality.  Namely, with probability at least $1 - \delta$
\begin{align}
\sigma_{\min}(\GradFun_\Numobs) &\geq \strongconvex -
\opnorm{\GradFun_\Numobs - \GradFun_\TruePar} \notag \\
& \geq \strongconvex - \opnorm{\GradFun_\Numobs + \frac{1}{\Numobs}
  \sum_{i=1}^\Numobs \nabla _\Par\ZFun(\Sam_i, \TruePar)} - \opnorm{
  \frac{1}{\Numobs} \sum_{i=1}^\Numobs \nabla _\Par\ZFun(\Sam_i,
  \TruePar) + \GradFun_\TruePar} \notag \\
& \geq \strongconvex- \sup_{\Par\in \ParSpace} \opnorm{
  \frac{1}{\Numobs} \sum_{i=1}^\Numobs \nabla ^2_{\Par} \ZFun(\Sam_j,
  \Par)} \enorm{\EstPar_\numobs - \TruePar} - \opnorm{
  \frac{1}{\Numobs} \sum_{i=1}^\Numobs \nabla_\Par\ZFun(\Sam_i,
  \TruePar) + \GradFun_\TruePar} \notag \\
& \stackrel{(i)}{\geq} \strongconvex- c\OrParc\Big(1 + \sqrt{
  \frac{\dimdelta}{\Numobs}} + \frac{(\dimdelta \log
  {\Numobs})^{\InvOrConc}}{\Numobs} \Big) \polyshort \sqrt{ \log
  \Numobs} \sqrt{ \frac{\dimdelta}{\Numobs}}-c\OrParb\Big( \sqrt{
  \frac{\dimdelta\log(\numobs/\delta)}{\Numobs}} + \frac{\dimdelta
  \log^2(\numobs/\delta)}{\Numobs} \Big) \notag \\
\label{eq:thetai-theta-eq4-pf3}
& \stackrel{(ii)}{\geq} \strongconvex-\polyshort \sqrt{ \log
  (\Numobs/\delta)} \sqrt{ \frac{\dimdelta}{\Numobs}}
\geq\strongconvex/2,
\end{align}
where step (i) follows from
Assumption~\ref{ass:convergence},~\Cref{lm:sup_op_norm_zfun,lm:op_norm_zfun}; step (ii)
 uses the sample size assumptions
in~\cref{eq:linear_approx_Z_est_sample_ahead}. 

Lastly, with
probability at least $1 - \delta$ for all
$i\in[\Numobs]$
\begin{align}
\enorm{\ZFun(\Sam_i, \LOO{\EstPar}{\numobs-1}{i})}
& \leq
\enorm{\ZFun(\Sam_i, \EstPar_\numobs)} + \enorm{\ZFun(\Sam_i,
  \LOO{\EstPar}{\numobs-1}{i}) - \ZFun(\Sam_i, \EstPar_\numobs)}
\notag \\ 
& \leq 
\enorm{\ZFun(\Sam_i, \EstPar_\numobs)} +
\sup_{\Par\in \ParSpace} \enorm{\nabla  \ZFun(\Sam_i, \Par)} \enorm{\LOO{\EstPar}{\numobs-1}{i}
  - \EstPar_\numobs} \notag \\ 
  & \stackrel{(i)}{\leq} 
  \polyshort\Big(
\log^2(\Numobs/\delta) \sqrt{\dimdelta} + \frac{ \sqrt{\dimdelta}
  \log^2(\Numobs/\delta)}{\Numobs} \log(\numobs/\delta) \dimdeltatil\Big)
\notag \\
  & \stackrel{(ii)}{\leq} 
\polyshort \log^2(\Numobs/\delta) \sqrt{\dimdelta},
\label{eq:hattheta-thetai_zfun}
\end{align}
where step (i) follows
from bound~\eqref{eq:LemJackknifeKey_eq1},  equation~\eqref{eq:LemJackknifeKey}(b)~and~\eqref{eq:union_gaussian_eq3} in the proof of bound~\eqref{eq:LemJackknifeKey_eq1};
step (ii) uses the sample size
condition~\eqref{eq:linear_approx_Z_est_sample_ahead}.\\

\myunder{Combining the bounds:} Substituting
equations~\eqref{eq:thetai-theta-eq4-pf3},~\eqref{eq:thetai-theta-eq4}(i),~\eqref{eq:thetai-theta-eq4}(ii),~and~\eqref{eq:hattheta-thetai_zfun} into~\cref{eq:thetai-theta-eq4-pf1} yields
the bound
\begin{align*}
 \Sma_2 & \leq \frac{1}{\Numobs} \sum_{i=1}^\Numobs
 \frac{1}{\strongconvex^2} \opnorm{ - \GradFunminusi +
   \GradFun_{\Numobs}} \enorm{\ZFun(\Sam_i,
   \LOO{\EstPar}{\numobs-1}{i})} \leq \polyshort
 \log^4(\Numobs/\delta) \frac{{\dimdelta}}{\Numobs},
\end{align*}
valid with probability at least $1 - \delta$.


\subsubsection{Proof of claim~\eqref{eq:sum_thetai-theta}(ii)}
\label{sec:pf_eq:sum_thetai-theta_claim2}

By the lower bound on $\sigma_{\min}( \GradFun_{\Numobs})$
in~\cref{eq:thetai-theta-eq4-pf3} from the proof of
claim~\eqref{eq:sum_thetai-theta}(i), it remains to show
\begin{align}
 \frac{1}{\Numobs} \enorm{\sum_{i=1}^\Numobs \ZFun(\Sam_i,
   \LOO{\EstPar}{\numobs-1}{i})} \leq \polyshort
 \log^4(\Numobs/\delta)
 \frac{\dimdelta}{\Numobs} \label{eq:sum_thetai-theta_claim2_pf_1}.
\end{align}
We establish this claim via the following two-step argument.

\textbf{Step 1.\quad} First, we approximate the LHS
of~\cref{eq:sum_thetai-theta_claim2_pf_1} by
$\frac{1}{\Numobs^2}\enorm{\sum_{i=1}^\Numobs \nabla _\Par\ZFun(\Sam_i,
  \EstPar_\numobs) \GradFun_\Numobs^{-1} \ZFun(\Sam_i,
  \EstPar_\numobs)}$, and show that
\begin{align}
 \enorm{ \frac{1}{\Numobs} \sum_{i=1}^\Numobs \ZFun(\Sam_i,
    \LOO{\EstPar}{\numobs-1}{i})- \frac{1}{\Numobs^2}
    \sum_{i=1}^\Numobs \nabla _\Par\ZFun(\Sam_i, \EstPar_\numobs)
    \GradFun_\Numobs^{-1} \ZFun(\Sam_i, \EstPar_\numobs)} 
    \leq
     \polyshort \log^4(\Numobs/\delta) \frac{\dimdelta}{\Numobs}\label{eq:sum_thetai-theta_claim2_pf_00}
\end{align}
 with probability at least $1-\delta/2.$

\textbf{Step 2.\quad} Second, we prove that
\begin{align}
\label{eq:sum_thetai-theta_claim2_pf_3}
 \frac{1}{\Numobs} \vecnorm{\sum_{i=1}^\Numobs
   \nabla _\Par\ZFun(\Sam_i,
   \EstPar_\numobs) \GradFun_\Numobs^{-1} \ZFun(\Sam_i,
   \EstPar_\numobs)}{2} \leq 
   \polyshort
 \log(\Numobs/\delta){\ParDim}
\end{align} 
 with probability at least $1-\delta/2.$ We defer the proof of Step 2 to Section~\ref{sec:pf_step2_in_sma_1}.

Combining the two steps yields claim.\\

\myunder{Proof of Step 1.}
\begin{align}
  \quad \enorm{ \frac{1}{\Numobs} \sum_{i=1}^\Numobs \ZFun(\Sam_i,
    \LOO{\EstPar}{\numobs-1}{i})- \frac{1}{\Numobs^2}
    \sum_{i=1}^\Numobs \nabla _\Par\ZFun(\Sam_i, \EstPar_\numobs)
    \GradFun_\Numobs^{-1} \ZFun(\Sam_i, \EstPar_\numobs)} \notag
& \leq \Sma_{31}+\Sma_{32}+\Sma_{33},
\end{align}
where 
\begin{align*}
\Sma_{31}&\defn \frac{1}{\Numobs} \enorm{\sum_{i=1}^\Numobs\big[
    \ZFun(\Sam_i, \LOO{\EstPar}{\numobs-1}{i}) - \ZFun(\Sam_i,
    \EstPar_\numobs)- \nabla _\Par\ZFun(\Sam_i,
    \EstPar_\numobs)(\LOO{\EstPar}{\numobs-1}{i} - \EstPar_\numobs)
    \big]}, \\ \Sma_{32}&\defn \frac{1}{\Numobs}
\enorm{\sum_{i=1}^\Numobs \nabla _\Par\ZFun(\Sam_i, \EstPar_\numobs)
  \big(\LOO{\EstPar}{\numobs-1}{i} - \EstPar_\numobs +
  \frac{\GradFun_\Numobs^{-1} \ZFun(\Sam_i,
    \LOO{\EstPar}{\numobs-1}{i} \big)}{\Numobs})}, \\ \Sma_{33}&\defn
\frac{1}{\Numobs^2} \enorm{\sum_{i=1}^\Numobs \nabla
  _\Par\ZFun(\Sam_i, \EstPar_\numobs)
  \GradFun_\Numobs^{-1}(\ZFun(\Sam_i, \LOO{\EstPar}{\numobs-1}{i}) -
  \ZFun(\Sam_i, \EstPar_\numobs))}.
  \end{align*}
  We establish bounds on  $\Sma_{31},\Sma_{32},\Sma_{33}$, respectively.
  From Lemma~\ref{lm:sup_op_norm_zfun} and the bound on $\enorm{\LOO{\EstPar}{\numobs-1}{i} - \EstPar_\numobs}$
in~\cref{eq:LemJackknifeKey_eq1}, we have
  \begin{align*}
  \Sma_{31}&\leq
    \frac{1}{\Numobs} \sum_{i=1}^\Numobs\sup_{\Par\in \ParSpace} \opnorm{
    \nabla ^2_\Par\ZFun(\Sam_i,
    \Par)} \enorm{\LOO{\EstPar}{\numobs-1}{i} - \EstPar_\numobs}^2 \\
    &\leq 
    \polyshort \log^4(\Numobs/\delta) \frac{\dimdelta\dimdeltatil^{3/2}}{\Numobs^2}
  \end{align*}with probability at least $1-\delta.$
   From~\cref{eq:thetai-theta_1_repr2,eq:union_gaussian_eq3,eq:hattheta-thetai_zfun,eq:thetai-theta-eq4}, we have
  \begin{align*}
  \Sma_{32}&\leq
  \frac{1}{\Numobs^2} \sum_{i=1}^\Numobs\sup_{\Par\in \ParSpace} \opnorm{
    \nabla _\Par\ZFun(\Sam_i, \Par)} \opnorm{ - \GradFunminusi
    +
    \GradFun_\Numobs^{-1}} \opnorm{\ZFun(\Sam_i,
    \LOO{\EstPar}{\numobs-1}{i})} \\
    &\leq
    \polyshort\log^6(\Numobs/\delta) 
  \frac{\dimdelta\dimdeltatil}{\Numobs^2}
  \end{align*}
  with probability at least $1-\delta.$ Lastly, using Lemma~\ref{lm:sup_op_norm_zfun},~equation~\eqref{eq:LemJackknifeKey_eq1}~and~\eqref{eq:thetai-theta-eq4}, we obtain
 \begin{align*}   
     \Sma_{33}&\leq
  \frac{1}{\Numobs^2} \sum_{i=1}^\Numobs
  \sup_{\Par\in \ParSpace} \opnorm{ \nabla _\Par\ZFun(\Sam_i,
    \Par)} \opnorm{\GradFun_\Numobs^{-1}} \sup_{\Par\in \ParSpace} \opnorm{
    \nabla \ZFun(\Sam_i, \Par)} \enorm{\LOO{\EstPar}{\numobs-1}{i} -
    \EstPar_\numobs} \\
    & \leq
\polyshort  \log^6(\Numobs/\delta)  \frac{
    \sqrt{\dimdelta} \dimdeltatil^2}{\Numobs^2}\end{align*}
    with probability at least $1-\delta.$
   Putting pieces together and using the sample size condition~\eqref{eq:linear_approx_Z_est_sample_ahead} yields bound~\eqref{eq:sum_thetai-theta_claim2_pf_00}.
%


\subsubsection{Proof of Step 2 in Section~\ref{sec:pf_eq:sum_thetai-theta_claim2}} \label{sec:pf_step2_in_sma_1}
Adopt the shorthand notations $\check{\GZInv}_i$, $\GZ_i$ and $\ZF_i$ for the
quantities $ \nabla _\Par\ZFun(\Sam_i,
\EstPar_\numobs) \GradFun_\Numobs^{-1}$, $ \nabla _\Par\ZFun(\Sam_i,
\EstPar_\numobs)$, and $\ZFun(\Sam_i, \EstPar_\numobs)$,
respectively. Let $s_j(\Mmat)$ denote the $j^{th}$-largest
singular value of a matrix $\Mmat$. At a high-level, the proof consists of two parts: (a) an algebraic
argument to reduce the problem to a simplified form; and (b) provide an upper bound for  the simplified problem. Specifically, in the simplification step, we show that
\begin{subequations}
\begin{align}
    \frac{1}{\Numobs} \vecnorm{\sum_{i=1}^\Numobs
    \nabla _\Par\ZFun(\Sam_i,
    \EstPar_\numobs) \GradFun_\Numobs^{-1} \ZFun(\Sam_i,
    \EstPar_\numobs)}{2} 
  &\leq
  \frac{\sqrt{\ParDim}}{\strongconvex\numobs}\cdot \sqrt{\sum_{k=1}^{\ParDim} \opnorm{\sum_{i=1}^{\Numobs} \ZF_{ik} \GZ_i}^2}\label{eq:step_2_in_sma_1_eq1}
\end{align}
and to bound the RHS of equation~\eqref{eq:step_2_in_sma_1_eq1}, we prove
that with probability at least $1 - \delta$
\begin{align}
 \label{eq:sum_thetai-theta_claim2_pf_8}  
|\ZF_{ik}| \leq \polyshort \log(\Numobs/\delta),~~~\text{for all } i
\in [\Numobs], k \in [\ParDim].
\end{align}
\end{subequations}
Step 2 then follows from combining these two results and a
$\ZF_{ik}$-weighted version of~\eqref{eq:bounded_grad_zfun_lm}
in~\Cref{lm:sup_op_norm_zfun} (which can be established following the
proof of~\eqref{eq:bounded_grad_zfun_lm} in
Lemma~\ref{lm:sup_op_norm_zfun} but with a slight modification in the
bound~\eqref{EqnTailEps}), more precisely,
\begin{align}
\label{eq:sum_thetai-theta_claim2_pf_9}.  
\opnorm{ \frac{1}{\Numobs} \sum_{i=1}^{\Numobs} \ZF_{ik} \GZ_i}^2 \leq
\polyshort \log^2(\Numobs/\delta) ,~~~ \mbox{for all $k \in
  [\ParDim]$}
\end{align}

\myunder{Proof of part (a):}
\begin{align}
& \frac{1}{\Numobs} \vecnorm{\sum_{i=1}^\Numobs \nabla_\Par
    \ZFun(\Sam_i, \EstPar_\numobs) \GradFun_\Numobs^{-1} \ZFun(\Sam_i,
    \EstPar_\numobs)}{2} = \sup_{\Direc \in \Sphere{\ParDim-1}}
  \frac{1}{\Numobs} \sum_{i=1}^{\Numobs} \Direc^{\mytrans}
  \check{\GZInv}_i\ZF_i = \sup_{\Direc \in \Sphere{\ParDim-1}}
  \tr\Big( \frac{1}{\Numobs} \sum_{i=1}^{\Numobs}
  \ZF_i\Direc^{\mytrans} \GZ_i\GradFun_\Numobs^{-1} \Big) \notag
  \\ &\quad \leq \sup_{\Direc \in \Sphere{\ParDim-1}}
  \opnorm{\GradFun_\Numobs^{-1}} \sum_{j=1}^\ParDim s_j\Big(
  \frac{1}{\Numobs} \sum_{i=1}^{\Numobs} \ZF_i\Direc^{\mytrans}
  \GZ_i\Big) \leq \sqrt{\ParDim} \opnorm{\GradFun_\Numobs^{-1}}
  \sup_{\Direc \in \Sphere{\ParDim-1}} \fronorm{ \frac{1}{\Numobs}
    \sum_{i=1}^{\Numobs} \ZF_i\Direc^{\mytrans} \GZ_i} \notag
  \\ &\quad \leq \frac{ \sqrt{\ParDim}}{\strongconvex} \sup_{\Direc
    \in \Sphere{\ParDim-1}} \fronorm{ \frac{1}{\Numobs}
    \sum_{i=1}^{\Numobs} \ZF_i\Direc^{\mytrans} \GZ_i}
   \label{eq:sum_thetai-theta_claim2_pf_4},
   \end{align}
   where the first inequality follows from Von Neumann's trace
   inequality~\cite{mirsky1975trace}, the second inequality uses
   $\fronorm{\cdot} = \sqrt{ \sumn s_j(\cdot)^2}$ and Cauchy-Schwartz
   inequality, and the last inequality
   uses~\cref{eq:thetai-theta-eq4-pf3} from previous proof.  Moreover,
   a little algebra shows
   \begin{align}
  &\qquad \sup_{\Direc \in \Sphere{\ParDim-1}} \fronorm{
       \frac{1}{\Numobs} \sum_{i=1}^{\Numobs} \ZF_i\Direc^{\mytrans}
       \GZ_i}^2 = \sup_{\Direc \in \Sphere{\ParDim-1}}
     \frac{1}{\Numobs^2} \sum_{i=1}^{\Numobs} \sum_{j=1}^{\Numobs}
     \tr(\ZF_i\Direc^{\mytrans} \GZ_i\GZ_j\Direc\ZF_j^{\mytrans})
     \notag \\ & = \sup_{\Direc \in \Sphere{\ParDim-1}}
     \frac{1}{\Numobs^2} \sum_{i=1}^{\Numobs} \sum_{j=1}^{\Numobs}
     \Direc^{\mytrans}(\ZF_i^{\mytrans} \ZF_j) \GZ_i\GZ_j\Direc =
     \sup_{\Direc \in \Sphere{\ParDim-1}} \Direc^{\mytrans} \Big[
       \frac{1}{\Numobs^2} \sum_{i=1}^{\Numobs}
       \sum_{j=1}^{\Numobs}(\ZF_i^{\mytrans} \ZF_j) \GZ_i\GZ_j\Big]
     \Direc \notag \\
\label{eq:sum_thetai-theta_claim2_pf_5}     
 & = \frac{1}{\Numobs^2} \opnorm{\sum_{i=1}^{\Numobs}
  \sum_{j=1}^{\Numobs} (\ZF_i^{\mytrans} \ZF_j) \GZ_i\GZ_j} =
\frac{1}{\Numobs^2} \opnorm{\sum_{i=1}^{\Numobs}
  \sum_{j=1}^{\Numobs} \sum_{k=1}^{\ParDim} \ZF_{ik} \ZF_{jk} \GZ_i
  \GZ_j}.
\end{align}
Also, note that
\begin{align}
\label{eq:sum_thetai-theta_claim2_pf_7}  
\opnorm{\sum_{i=1}^{\Numobs} \sum_{j=1}^{\Numobs} \sum_{k=1}^{\ParDim} \ZF_{ik}
  \ZF_{jk} \GZ_i \GZ_j} & \leq
\sum_{k=1}^{\ParDim} \opnorm{\sum_{i=1}^{\Numobs} \sum_{j=1}^{\Numobs} \ZF_{ik}
  \ZF_{jk} \GZ_i\GZ_j} \leq
\sum_{k=1}^{\ParDim} \opnorm{\sum_{i=1}^{\Numobs} \ZF_{ik} \GZ_i}^2.
\end{align}Putting the pieces together completes the proof.\\

\myunder{Proof of part (b):}
a Taylor series expansion yields
\begin{align}
|\ZF_{ik}| 
& = 
|[e_k^{\mytrans} \ZFun(\Sam_i, \EstPar_\numobs)-e_k^{\mytrans}
\ZFun(\Sam_i, \LOO{\EstPar}{\numobs-1}{i})] + [e_k^{\mytrans} \ZFun(\Sam_i,
\LOO{\EstPar}{\numobs-1}{i})-e_k^{\mytrans} \ZFun(\Sam_i, \TruePar)] +
e_k^{\mytrans} \ZFun(\Sam_i, \TruePar)| \notag \\
& \revdefn \Sma_{31}+\Sma_{32}+\Sma_{33}+\Sma_{34}\label{eq:sum_thetai-theta_claim2_pf_10},
\end{align}
where
\begin{align*}
    \Sma_{31}
    &\defn
    \sup_{\Par \in \ParSpace} \enorm{e_k^{\mytrans} \nabla
  _\Par\ZFun(\Sam_i, \Par)} \enorm{\LOO{\EstPar}{\numobs-1}{i} -
  \EstPar_\numobs},~~~~~~\,\,\,\,
    \Sma_{32}
    \defn
   |\nabla_\Par \ZFun(\Sam_i, \TruePar)[e_k,
  \LOO{\EstPar}{\numobs-1}{i}- \TruePar]|,\\
      \Sma_{33}
    &\defn
   |\sup_{\Par\in \ParSpace}|\nabla ^2_\Par\ZFun(\Sam_i,
\TruePar)[e_k, \LOO{\EstPar}{\numobs-1}{i}- \TruePar,
  \LOO{\EstPar}{\numobs-1}{i}- \TruePar]|,~~~~~~~~~~~~~~
    \Sma_{34}
    \defn
  |e_k^{\mytrans} \ZFun(\Sam_i,
\TruePar)|.
\end{align*}

Since $\LOO{\EstPar}{\numobs-1}{i} - \TruePar$ is independent of $\Sam_i$, we have
from Assumption~\ref{ass:tail} that
   \begin{subequations}
   \begin{align}
 &\Prob(\Sma_{31} \geq t_1 \sqrt{\ParDim} \enorm{\LOO{\EstPar}{\numobs-1}{i} -
  \EstPar_\numobs}) \leq
     \Prob(| \nabla _\Par\ZFun(\Sam_i, \Par)[e_k,e_j]|\geq
     t_1 ,~\text{for some~}j\in[\ParDim]) \leq 2\ParDim\exp\Big( -
     \frac{ct_1}{\OrParb} \Big), \label{eq:sum_thetai-theta_claim2_pf_105}
     \\
      &\Prob( \Sma_{32}\geq
     t_1\enorm{\LOO{\EstPar}{\numobs-1}{i} - \TruePar}) \leq
     2\exp\Big( -
     \frac{ct_1}{\OrParb} \Big), \label{eq:sum_thetai-theta_claim2_pf_11}
     \\ 
     &
     \Prob(  \Sma_{33}\geq
     t_2\enorm{\LOO{\EstPar}{\numobs-1}{i} - \TruePar}^2) \leq
     2\exp\Big( -
     \frac{ct_2^{\OrConc}}{\OrParc^\OrConc} \Big), \label{eq:sum_thetai-theta_claim2_pf_115}
     \\ 
     &
     \Prob( \Sma_{34}\geq t_3) \leq
     2\exp\Big( -
     \frac{ct_2}{\OrPara} \Big) \label{eq:sum_thetai-theta_claim2_pf_12}.
   \end{align}    
   \end{subequations}
Letting $t_1 = c \OrParb \log (\Numobs/\delta)$, $t_2 = c \OrParc
\log^{\InvOrConc} (\Numobs/\delta)$ and $t_3 = c \OrPara \log
(\Numobs/\delta)$ in the inequalities, and combining them with
Assumption~\ref{ass:convergence} and a union bound over
$i\in[\Numobs],k\in[\ParDim]$, we obtain
\begin{align*}
|\ZF_{ik}|
 &\leq
 \Sma_{31}+ \Sma_{32}+ \Sma_{33}+ \Sma_{34}\\
 \\
 & \leq \polyshort \left(
  \log (\Numobs/\delta) \frac{ \sqrt{\dimdelta\ParDim}
    \log^2(\Numobs/\delta)}{\Numobs} + \log (\Numobs/\delta) \sqrt{
    \frac{ \dimdeltatil \log \Numobs}{\Numobs}} + \log^{\InvOrConc}
  (\Numobs/\delta) \log\numobs\frac{\dimdeltatil}{\numobs} + \log
  (\Numobs/\delta) \right) \\ 
  & \leq \polyshort \log(\Numobs/\delta)
   \end{align*}
   for all $i\in[\numobs],k\in[\ParDim]$ with probability at least $1-\delta,$
where the last line follows from the sample size assumption
in~\eqref{eq:linear_approx_Z_est_sample_ahead}.


\section{Jackknife variance estimation: Proof of~\Cref{PropJackVarEst}} 
\label{SecProofPropJackVarEst}

We now prove our guarantee on the consistency of jackknife variance
estimation.  The argument is based on decomposing the difference
$\numobs \JackVar - \VarTrue$ into the sum of four terms.  These terms
depend on various approximate versions of the matrix $\Jmat$: in
particular, we define \mbox{$\EmpGZF \defn- {\Numobs^{-1}}
  \sum_{i=1}^\Numobs \nabla_\Par \ZFun(\Sam_i,
  \EstPar_\numobs)$,}\footnote{Note that the definitions are
consistent with those in~\eqref{eq:def_js_1} from
Section~\ref{sec:pf_eq:LemJackknifeKey_eq4}.} as well as
\begin{align*}
\EmpGZFa & \defn - \Numobs^{-1} \sum_{j=1}^\Numobs \int_0^1
\nabla_{\Par} \ZFun(\Sam_j, \EstPar_\numobs + t \LooDiff) dt, \quad
\EmpGZFb\defn - {\Numobs^{-1}} \sum_{j=1}^\Numobs
\nabla_\Par\ZFun(\Sam_j, \LOO{\EstPar}{\numobs-1}{i}), \\
\EmpGZFc & \defn - {\Numobs^{-1}} \sum_{j \neq i}
\nabla_\Par\ZFun(\Sam_j, \LOO{\EstPar}{\numobs-1}{i}), \quad
\mbox{and} \quad \EmpGZFcall\defn - {\Numobs^{-1}} \sum_{j=1 }^\Numobs
\nabla_\Par\ZFun(\Sam_j, \LOO{\EstPar}{\numobs-1}{i}).
\end{align*}

With this notation, we have the decomposition $\numobs \JackVar -
\TrueVar^2 = \sum_{j=1}^4 \vterm_4$, where
\begin{subequations}
\begin{align}
\label{eq:jac_var_pf_result1} 
\vterm_1& \defn |\Numobs \JackVar^2 - (\Numobs-1) \sum_{i=1}^\Numobs[
  \nabla
  \TargetFun(\EstPar_\numobs)^{\mytrans}(\LOO{\EstPar}{\numobs-1}{i} -
  \EstPar_\numobs)]^2|, \\
 \label{eq:jac_var_pf_result2}
 \vterm_2&
\defn
(\Numobs-1) \Big| \sumn\Big\{[ \nabla
   \TargetFun(\EstPar_\numobs)^{\mytrans}(\LOO{\EstPar}{\numobs-1}{i} -
   \EstPar_\numobs)]^2-[ \nabla
   \TargetFun(\LOO{\EstPar}{\numobs-1}{i})^{\mytrans}
   \frac{\EmpGZFc^{-1} \ZFun(\Sam_i,
     \LOO{\EstPar}{\numobs-1}{i})}{\Numobs}]^2\Big\} \Big|,  \\
 \label{eq:jac_var_pf_result3}
 \vterm_3&
\defn
\frac{\Numobs-1}{\Numobs^2}| \sumn \Big\{[
  \nabla\TargetFun(\LOO{\EstPar}{\numobs-1}{i})^{\mytrans}{\EmpGZFc^{-1} \ZFun(\Sam_i,
    \LOO{\EstPar}{\numobs-1}{i})}]^2-[
  \nabla\TargetFun(\TruePar)^{\mytrans}{\GradFun_\TruePar^{-1} \ZFun(\Sam_i,
    \TruePar)}]^2\Big\}| , \\
\vterm_4&
\defn
\frac{\Numobs-1}{\Numobs}| \frac{1}{\Numobs} \sum_{i=1}^\Numobs[
  \nabla\TargetFun(\TruePar)^{\mytrans}{\GradFun_\TruePar^{-1} \ZFun(\Sam_i,
    \TruePar)}]^2 - \E[
  \nabla\TargetFun(\TruePar)^{\mytrans}{\GradFun_\TruePar^{-1} \ZFun(\Sam_i,
    \TruePar)}]^2|
    \label{eq:jac_var_pf_result4}.
\end{align}
\end{subequations}
We prove the claim~\eqref{EqnJackVarEst} by showing that all of the
following inequalities hold with probability at least $1 -\delta$:
\begin{subequations}
\begin{align}
\label{EqnJackVarOneTwo}  
  |\vterm_1| \stackrel{(i)}{\leq} \polyshort \log^6(\Numobs/\delta)
\frac{\dimdelta^{3/2}}{\Numobs},
\quad |\vterm_2|
    \stackrel{(ii)}{\leq} \polyshort \log^5(\Numobs/\delta) \frac{\dimdelta}{\Numobs}, \\
\label{EqnJackVarThreeFour}        
      |\vterm_3| \stackrel{(iii)}{\leq} \polyshort
      \log^3(\Numobs/\delta) \sqrt{ \frac{\dimdelta}{\Numobs}}, \quad
      \mbox{and} \quad |\vterm_4| \stackrel{(iv)}{\leq} \polyshort
      \log^3(\Numobs/\delta) \frac{1}{ \sqrt{\Numobs}},
\end{align}
\end{subequations}
for some constants $\polyshort= \poly$.

Combining these bounds and and applying the triangle inequality, we
find that
\begin{align*}
|\Numobs \JackVar^2 - \VarTrue^2| = |\Numobs \JackVar^2 - \E[ \nabla
  \TargetFun(\TruePar)^{\mytrans} {\GradFun_\TruePar^{-1} \ZFun(\Sam_i,
    \TruePar)} ]^2 | \leq \polyshort \log^6(\Numobs/\delta)
\frac{\dimdelta^{3/2}}{\Numobs}.
\end{align*}
The stated claim of Proposition~\ref{PropJackVarEst} follows since the RHS of this
inequality goes to zero by the sample size assumption. \\

\noindent Thus, the remainder of our proof is devoted to establishing
the four claimed bounds.


\paragraph{Proof of inequality~\eqref{EqnJackVarOneTwo}(i):}

Performing a Taylor expansion around $\EstPar_\numobs$, we obtain
\begin{align*}
\Numobs \JackVar^2 & = (\Numobs-1) \sum_{i=1}^\Numobs
\Big[\TargetFun(\LOO{\EstPar}{\numobs-1}{i}) -
  \TargetFun(\EstPar_\numobs)-
  \frac{\sum_{j=1}^\Numobs[\TargetFun(\EstPar_\numobs^{(-j)}) -
      \TargetFun(\EstPar_\numobs)]}{\Numobs} \Big]^2 \\
&=
(\Numobs-1) \sum_{i=1}^\Numobs\Big[ \vtermtil^{(i)}_{12}+\vtermtil^{(i)}_{12}
\Big]^2,
\end{align*}
where 
\begin{align*}
\vtermtil^{(i)}_{11} & \defn \nabla
\TargetFun(\EstPar_\numobs)^{\mytrans}(\LOO{\EstPar}{\numobs-1}{i} -
\EstPar_\numobs)\\
\vtermtil^{(i)}_{12} & \defn - \frac{
  \nabla\TargetFun(\EstPar_\numobs)^{\mytrans}
  \sum_{j=1}^\Numobs(\EstPar_\numobs^{(-j)} -
  \EstPar_\numobs)}{\Numobs} +
\nabla^2\tau(\IntPar^{(-j)})[\LOO{\EstPar}{\numobs-1}{i} -
  \EstPar_\numobs, \LOO{\EstPar}{\numobs-1}{i} - \EstPar_\numobs] \\
&~~~~~~ -\frac{1}{\Numobs} \sum_{j=1}^\Numobs
\nabla^2\tau(\IntPar^{(-i)})[\EstPar_\numobs^{(-j)} - \EstPar_\numobs,
  \EstPar_\numobs^{(-j)} - \EstPar_\numobs] .
\end{align*}
and
$\IntPar^{(-i)} = \lambda_i\EstPar_\numobs + (1 - \lambda_i)
\LOO{\EstPar}{\numobs-1}{i}$ for some $\lambda_i\in[0,1]$.
Let $\vterm^{(i)}_{11}=|\vtermtil^{(i)}_{11}|$ and $\vterm^{(i)}_{12}=|\vtermtil^{(i)}_{12}|$. By some basic algebra, we have
\begin{align*}
\vterm_1 & \leq 2 (\Numobs-1) \sum_{i=1}^\numobs (\vterm^{(i)}_{11} +
\vterm^{(i)}_{12}) \vterm^{(i)}_{12}.
\end{align*}

Since by~\Cref{LemJackknifeKey} and the Lipschitz continuity assumptions on
$\nabla\TargetFun, \nabla^2\TargetFun$, we have
\begin{align*}
{\vterm^{(i)}_{11} } & \leq \polyshort \frac{
  \sqrt{\dimdelta}}{\Numobs} \log^2(\Numobs/\delta) , \text{ and} \\
{\vterm^{(i)}_{12} } & \leq \polyshort \log^4(\Numobs/\delta)
\frac{\dimdelta}{\Numobs^2}
\end{align*}
for all $i \in [\numobs]$ with probability at least $1-\delta$. Putting pieces together 
completes the proof of inequality~\eqref{EqnJackVarOneTwo}(i).


\paragraph{Proof of inequality~\eqref{EqnJackVarOneTwo}(ii):}

Recall equation~\eqref{eq:thetai-theta_1_repr2} that
\begin{align*}
\LOO{\EstPar}{\numobs-1}{i} - \EstPar_\numobs
= - \GradFunminusi^{-1}
\ZFun(\Sam_i, \LOO{\EstPar}{\numobs-1}{i})/\Numobs.
\end{align*}
Therefore, we have
\begin{align}
\vterm_2 & = (\Numobs-1) \Big| \sumn\Big\{[ \nabla
  \TargetFun(\EstPar_\numobs)^{\mytrans}(\LOO{\EstPar}{\numobs-1}{i} -
  \EstPar_\numobs)]^2 - [ \nabla
  \TargetFun(\LOO{\EstPar}{\numobs-1}{i})^{\mytrans}
  \frac{\EmpGZFc^{-1} \ZFun(\Sam_i,
    \LOO{\EstPar}{\numobs-1}{i})}{\Numobs}]^2 \Big\} \Big| \notag \\ &
= \frac{\Numobs-1}{\Numobs^2}| \sumn \Big\{
[\nabla\TargetFun(\EstPar_\numobs)^{\mytrans}{\EmpGZFa^{-1}
    \ZFun(\Sam_i, \LOO{\EstPar}{\numobs-1}{i})}]^2 -
[\nabla\TargetFun(\LOO{\EstPar}{\numobs-1}{i})^{\mytrans}{\EmpGZFc^{-1}
    \ZFun(\Sam_i, \LOO{\EstPar}{\numobs-1}{i})}]^2\Big\}| \notag \\
\label{eq:jac_var_pf_eq_3}
& \leq 2 \sup_{i \in [\Numobs]} \vterm_{21}^{(i)} (\vterm_{21}^{(i)} +
\vterm_{22}^{(i)}),
 \end{align}
 where
 \begin{align*}
     \vterm_{21}^{(i)} & \defn|
     \nabla\TargetFun(\EstPar_\numobs)^{\mytrans}{\EmpGZFa^{-1}
       \ZFun(\Sam_i, \LOO{\EstPar}{\numobs-1}{i})} - \nabla
     \TargetFun(\LOO{\EstPar}{\numobs-1}{i})^{\mytrans}{\EmpGZFc^{-1}
       \ZFun(\Sam_i,
       \LOO{\EstPar}{\numobs-1}{i})}|,\\
\vterm_{22}^{(i)} & \defn |\nabla \TargetFun(\LOO{\EstPar}{\numobs
  -1}{i})^{\mytrans}{\EmpGZFc^{-1} \ZFun(\Sam_i,
  \LOO{\EstPar}{\numobs-1}{i})}|.
 \end{align*}
Let us collect some results used to bound $ \vterm_{21}^{(i)}$ and
$\vterm_{22}^{(i)}$. We claim that with probability at least $1 -
\delta$, for all $i \in [\Numobs]$, we have
\begin{subequations}
\begin{align}
\label{eq:jac_var_pf_claim_0}     
\enorm{\TargetFun(\EstPar_\numobs)},~\enorm{\TargetFun(\LOO{\EstPar}{\numobs-1}{i})},
~\enorm{ \nabla\TargetFun(\EstPar_\numobs)},~ &\enorm{
  \nabla\TargetFun(\LOO{\EstPar}{\numobs-1}{i})},~\opnorm{\EmpGZFa^{-1}},~\opnorm{\EmpGZFc^{-1}},~\opnorm{\EmpGZFcall^{-1}}
\leq \polyshort \\
\label{eq:jac_var_pf_claim_1}
\enorm{\ZFun(\Sam_i, \LOO{\EstPar}{\numobs-1}{i})} & \leq \polyshort
\log^2(\Numobs/\delta) \sqrt{\dimdelta} \\
|\nabla
\TargetFun(\LOO{\EstPar}{\numobs-1}{i})^{\mytrans}{\EmpGZFc^{-1}
  \ZFun(\Sam_i, \LOO{\EstPar}{\numobs-1}{i})}| & \leq
\polyshort\enorm{
  \nabla\TargetFun(\LOO{\EstPar}{\numobs-1}{i})^{\mytrans}
  \EmpGZFc^{-1}} \log (\Numobs/\delta) \notag \\
\label{eq:jac_var_pf_claim_2}
& \leq \polyshort \log (\Numobs/\delta) \\
|\nabla\TargetFun(\LOO{\EstPar}{\numobs-1}{i})^{\mytrans}{\EmpGZFc^{-1}
  \nabla\ZFun(\Sam_i, \LOO{\EstPar}{\numobs-1}{i})}| & \leq \polyshort
\enorm{ \nabla \TargetFun(\LOO{\EstPar}{\numobs-1}{i})^{\mytrans}
  \EmpGZFc^{-1}} \sqrt{\ParDim} \log (\Numobs/\delta) \notag \\
\label{eq:jac_var_pf_claim_25} 
& \leq \polyshort \sqrt{\ParDim} \log (\Numobs/\delta), \\
\enorm{\nabla \TargetFun(\EstPar_\numobs) - \nabla
  \TargetFun(\LOO{\EstPar}{\numobs-1}{i})},~~\enorm{\EmpGZFa^{-1} -
  \EmpGZFcall^{-1}} & \leq \polyshort \enorm{\EstPar_\numobs -
  \LOO{\EstPar}{\numobs-1}{i}} \notag \\
\label{eq:jac_var_pf_claim_3}
& \leq \polyshort \frac{\sqrt{\dimdelta}}{\Numobs}
\log^2(\Numobs/\delta).
\end{align}
\end{subequations}
We return the proof of these claims momentarily.

Given these bounds, we have from claim~\eqref{eq:jac_var_pf_claim_2}
that
\begin{align}
\vterm_{22}^{(i)}\leq\polyshort\sqrt{\ParDim}\log(\numobs/\delta)
\label{eq:jac_var_pf_eq_40}
\end{align}
for all $i \in [\numobs]$ with probability at least $1 - \delta.$
Moreover, by the triangle inequality
\begin{align}
\vterm_{21}^{(i)} & \leq \enorm{ \nabla \TargetFun(\EstPar_\numobs) -
  \nabla \TargetFun(\LOO{\EstPar}{\numobs-1}{i})}
\opnorm{\EmpGZFa^{-1}} \enorm{\ZFun(\Sam_i,
  \LOO{\EstPar}{\numobs-1}{i})} + \enorm{ \nabla
  \TargetFun(\LOO{\EstPar}{\numobs-1}{i})(\EmpGZFa^{-1} -
  \EmpGZFcall^{-1})} \enorm{\ZFun(\Sam_i,
  \LOO{\EstPar}{\numobs-1}{i})} \notag \\
& \quad + \enorm{ \nabla
  \TargetFun(\LOO{\EstPar}{\numobs-1}{i})(\EmpGZFcall^{-1} -
  \EmpGZFc^{-1})} \enorm{\ZFun(\Sam_i, \LOO{\EstPar}{\numobs-1}{i})}
\notag \\
& \stackrel{(i)}{\leq} \enorm{ \nabla \TargetFun(\EstPar_\numobs) -
  \nabla \TargetFun(\LOO{\EstPar}{\numobs-1}{i})}
\opnorm{\EmpGZFa^{-1}} \enorm{\ZFun(\Sam_i,
  \LOO{\EstPar}{\numobs-1}{i})} + \enorm{ \nabla
  \TargetFun(\LOO{\EstPar}{\numobs-1}{i})} \opnorm{\EmpGZFa^{-1} -
  \EmpGZFcall^{-1}} \enorm{\ZFun(\Sam_i, \LOO{\EstPar}{\numobs-1}{i})}
\notag \\
& \quad + \frac{1}{\Numobs} \enorm{ \nabla
  \TargetFun(\LOO{\EstPar}{\numobs-1}{i}) \EmpGZFc^{-1}{ \nabla
    \ZFun(\Sam_i, \LOO{\EstPar}{\numobs-1}{i})}}
\opnorm{\EmpGZFcall^{-1}} \enorm{\ZFun(\Sam_i,
  \LOO{\EstPar}{\numobs-1}{i})} \notag \\
\label{eq:jac_var_pf_eq_4}
& \stackrel{(ii)}{\leq} \polyshort \log^4(\Numobs/\delta)
\frac{\dimdelta}{\Numobs}
\end{align}
for all $i\in[\numobs]$ with probability at least $1-\delta$, where
step (i) uses the submultiplicity of the matrix operator norm and
Woodbury's matrix identity; step (ii) follows from the
claims~\eqref{eq:jac_var_pf_claim_0},~\eqref{eq:jac_var_pf_claim_1},~\eqref{eq:jac_var_pf_claim_25},~\eqref{eq:jac_var_pf_claim_3}
that bound each individual term.

Substituting
equation~\eqref{eq:jac_var_pf_eq_40}~and~\eqref{eq:jac_var_pf_eq_4}
into equation~\eqref{eq:jac_var_pf_eq_3} yields
inequality~\eqref{EqnJackVarOneTwo}(ii).  \\

\myunder{Proof of the claims:}\, claim~\eqref{eq:jac_var_pf_claim_0}
follows from the Lipschitz continuity of $\TargetFun,
\nabla\TargetFun$, Assumption~\ref{ass:tail} and an argument similar
to equation~\eqref{eq:thetai-theta-eq4-pf3} from the proof of
Lemma~\ref{LemJackknifeKey}.  Claim~\eqref{eq:jac_var_pf_claim_1} is
derived in~\cref{eq:hattheta-thetai_zfun} from the same proof. The
first lines of
claims~\eqref{eq:jac_var_pf_claim_2}~and~\eqref{eq:jac_var_pf_claim_25}
use Assumption~\ref{ass:tail} and the fact that $ \nabla
\TargetFun(\LOO{\EstPar}{\numobs-1}{i}), \EmpGZFc,
\LOO{\EstPar}{\numobs-1}{i}$ are independent of $\Sam_i$; the second
lines of
claims~\eqref{eq:jac_var_pf_claim_2}~and~\eqref{eq:jac_var_pf_claim_25}
follow from claim~\eqref{eq:jac_var_pf_claim_0}.  The first line of
claim~\eqref{eq:jac_var_pf_claim_3} follows from the Lipschitz
continuity of $ \nabla \TargetFun$,
claim~\eqref{eq:jac_var_pf_claim_1}, Woodbury's matrix identity and an
argument similar to~\eqref{eq:thetai-theta-eq4} in the proof of
Lemma~\ref{LemJackknifeKey}; the second line of
claim~\eqref{eq:jac_var_pf_claim_3} follows from
equation~\eqref{eq:LemJackknifeKey_eq1}~in~\Cref{LemJackknifeKey}.


\paragraph{Proof of inequality~\eqref{EqnJackVarThreeFour}(iii):}

Applying the triangle inequality, we have the upper bound
\begin{align}
 \vterm_3 & = \frac{\Numobs-1}{\Numobs^2}|\sum_{i=1}^\Numobs\Big\{[
   \nabla
   \TargetFun(\LOO{\EstPar}{\numobs-1}{i})^{\mytrans}{\EmpGZFc^{-1}
     \ZFun(\Sam_i, \LOO{\EstPar}{\numobs-1}{i})}]^2-[ \nabla
   \TargetFun(\TruePar)^{\mytrans}{\GradFun_\TruePar^{-1}
     \ZFun(\Sam_i, \TruePar)}]^2\Big\}| \notag \\
\label{eq:jac_var_pf_eq_5} 
 & \leq 2 \sup_{i \in [\Numobs]} \vterm_{31}^{(i)} \; \big \{
\vterm_{31}^{(i)} + \vterm_{32}^{(i)} \big \},
 \end{align}
where \mbox{$\vterm_{32}^{(i)} \defn |\nabla
  \TargetFun(\TruePar)^{\mytrans}{\GradFun_\TruePar^{-1} \ZFun(\Sam_i,
    \TruePar)}|$,} and
\begin{align*}
\vterm_{31}^{(i)} & \defn |\nabla
\TargetFun(\LOO{\EstPar}{\numobs-1}{i})^{\mytrans}{\EmpGZFc^{-1}
  \ZFun(\Sam_i, \LOO{\EstPar}{\numobs-1}{i})} - \nabla
\TargetFun(\TruePar)^{\mytrans}{\GradFun_\TruePar^{-1} \ZFun(\Sam_i,
  \TruePar)}|
 \end{align*}
Likewise, we claim the following holds with probability at least $1 -
\delta$ for all $i\in[\Numobs]$,
\begin{subequations}
\begin{align}
\opnorm{\GradFun_\TruePar^{-1} - \EmpGZFc^{-1}}
& \leq 
\polyshort \sqrt{ \frac{{\dimdelta}}{\Numobs}} \log (\Numobs/\delta) \label{eq:jac_var_pf_claim_4}, \\
\enorm{ \nabla  \TargetFun(\TruePar) -  \nabla  \TargetFun(\LOO{\EstPar}{\numobs-1}{i})}
  &\leq 
  \polyshort \sqrt{ \frac{{\dimdelta}}{\Numobs}} \log (\Numobs/\delta)
\label{eq:jac_var_pf_claim_45} \\
\vterm_{32}^{(i)}
 & \leq 
 \polyshort \log (\Numobs/\delta) \label{eq:jac_var_pf_claim_5}.
\end{align}
\end{subequations}
We return to the proof of these claims later.

It remains to provide a bound on $\vterm_{31}^{(i)}$.  By the triangle
inequality, we have the upper bound \mbox{$\vterm_{31}^{(i)} \leq
  \vterm_{311}^{(i)} + \vterm_{312}^{(i)} + \vterm_{313}^{(i)}$,}
where
 \begin{align*}
  \vterm_{311}^{(i)}
  &\defn 
| \nabla  \TargetFun(\LOO{\EstPar}{\numobs-1}{i})^{\mytrans}{\EmpGZFc^{-1} \ZFun(\Sam_i, \LOO{\EstPar}{\numobs-1}{i})}
-
 \nabla  \TargetFun(\TruePar)^{\mytrans}{\EmpGZFc^{-1} \ZFun(\Sam_i, \LOO{\EstPar}{\numobs-1}{i})}|, \notag\\
  \vterm_{312}^{(i)}
  &\defn
| \nabla  \TargetFun(\TruePar)^{\mytrans}{\EmpGZFc^{-1} \ZFun(\Sam_i, \LOO{\EstPar}{\numobs-1}{i})} -  \nabla  \TargetFun(\TruePar)^{\mytrans}{\GradFun_\TruePar^{-1} \ZFun(\Sam_i, \LOO{\EstPar}{\numobs-1}{i})}| \notag, \\
 \vterm_{313}^{(i)}
  &\defn
| \nabla  \TargetFun(\TruePar)^{\mytrans}{\GradFun_\TruePar^{-1} \ZFun(\Sam_i, \LOO{\EstPar}{\numobs-1}{i})} -  \nabla  \TargetFun(\TruePar)^{\mytrans}{\GradFun_\TruePar^{-1} \ZFun(\Sam_i, \TruePar)}|.
\end{align*}
Using the fact that $\nabla \TargetFun(\LOO{\EstPar}{\numobs-1}{i}),
\nabla \TargetFun(\TruePar),\EmpGZFc^{-1},\LOO{\EstPar}{\numobs-1}{i}$
are independent of $\Sam_i$ and Assumption~\ref{ass:tail}, we have
\begin{align*}
\vterm_{311}^{(i)} & \leq \; \polyshort \enorm{ \nabla
  \TargetFun(\LOO{\EstPar}{\numobs-1}{i}) - \nabla
  \TargetFun(\TruePar)} \opnorm{\EmpGZFc^{-1}} \log (\Numobs/\delta)
\; \stackrel{(i)}{\leq} \; \polyshort \log^2(\Numobs/\delta) \sqrt{
  \frac{\dimdelta}{{\Numobs}}}
\end{align*}
with probability at least $1-\delta.$ Here step~(i) uses
equations~\eqref{eq:jac_var_pf_claim_0},
~\eqref{eq:jac_var_pf_claim_45}, and
\mbox{Assumption~\ref{ass:convergence}.} For $\vterm_{312}^{(i)}$,
similarly
\begin{align*}
 \vterm_{312}^{(i)} & {\leq} \polyshort \enorm{ \nabla
   \TargetFun(\TruePar)} \opnorm{\EmpGZFc^{-1} -
   \GradFun_\TruePar^{-1}} \log (\Numobs/\delta) \;
 \stackrel{(ii)}{\leq} \; \polyshort \log^2(\Numobs/\delta) \sqrt{
   \frac{\dimdelta}{{\Numobs}}}
\end{align*}
with probability at least $1-\delta,$
where step~(ii) uses equation~\eqref{eq:jac_var_pf_claim_4}. For $\vterm_{313}^{(i)}$, performing a Taylor expansion, we find
\begin{align*}
    \vterm_{313}^{(i)}
  &
{\leq} 
 |\nabla  \ZFun(\Sam_i, \TruePar)[\GradFun_\TruePar^{-1} \nabla  \TargetFun(\TruePar), \LOO{\EstPar}{\numobs-1}{i} - \TruePar]|
 +  \sup_{\Par\in \ParSpace}
\big| \nabla ^2\ZFun(\Sam_i, \Par)[\GradFun_\TruePar^{-1} \nabla  \TargetFun(\TruePar), \LOO{\EstPar}{\numobs-1}{i} - \TruePar, \LOO{\EstPar}{\numobs-1}{i}-\TruePar]\big|
\\
 & \stackrel{(iii)}{\leq}
 \polyshort \log^2(\Numobs/\delta) \sqrt{ \frac{\dimdelta}{{\Numobs}}} 
\end{align*}
with probability at least $1-\delta,$
where step~(iii) follows from  Assumption~\ref{ass:tail},~\ref{ass:convergence}, and the fact that  $\GradFun_\TruePar^{-1}
\nabla  \TargetFun(\TruePar), \LOO{\EstPar}{\numobs-1}{i} - \TruePar$ are
independent of $\Sam_i$. Putting pieces together, we obtain
\begin{align}
     \vterm_{31}^{(i)}
  &\leq \polyshort \log^2(\Numobs/\delta) \sqrt{ \frac{\dimdelta}{{\Numobs}}}, \label{eq:jac_var_pf_eq_7}
\end{align}
for all $i$ with probability at least $1-\delta$.

Putting the bounds on $\vterm^{(i)}_{31}$ and $\vterm^{(i)}_{32}$ together yields inequality~\eqref{EqnJackVarThreeFour}(iii).\\

\myunder{Proof of the claims:}\, claim~\eqref{eq:jac_var_pf_claim_4}
follows from Woodbury's matrix identity,
claim~\eqref{eq:jac_var_pf_claim_0} in the previous proof and the
bound
\begin{align*}
\opnorm{\EmpGZFc - \GradFun_\TruePar} & \leq  \opnorm{ \frac{ \nabla  \ZFun(\Sam_i, \LOO{\EstPar}{\numobs-1}{i})}{\Numobs}} + 
\opnorm{\EmpGZFb - \EmpGZF} + \opnorm{\EmpGZF - \GradFun_\TruePar} \\
 & \leq 
 \frac{\sup_{\Par\in \ParSpace} \opnorm{ \nabla  \ZFun(\Sam_i, \Par)}}{\Numobs} +  \frac{1}{\Numobs} \sup_{\Par\in \ParSpace} \opnorm{\sum_{i=1}^\Numobs \nabla ^2\ZFun(\Sam_i, \Par)} \enorm{\LOO{\EstPar}{\numobs-1}{i} - \TruePar} + \opnorm{\EmpGZF - \GradFun_\TruePar} \\
 & \leq 
\polyshort\Big( \log^2(\numobs/\delta) \frac{\dimdeltatil}{\Numobs} +  \log (\Numobs/\delta) \sqrt{ \frac{\dimdelta}{\Numobs}} +  \sqrt{ \frac{\dimdelta}{\Numobs}} \Big) \leq \polyshort \log (\Numobs/\delta) \sqrt{ \frac{\dimdelta}{\Numobs}},
\end{align*}
where the third line uses Assumption~\ref{ass:convergence},~\Cref{lm:sup_op_norm_zfun,lm:op_norm_zfun}.  Claim~\eqref{eq:jac_var_pf_claim_45} is due to the Lipschitz continuity of $ \nabla  \TargetFun$ and Assumption~\ref{ass:convergence}. Claim~\eqref{eq:jac_var_pf_claim_5} follows from Assumption~\ref{ass:tail}.


\paragraph{Proof of inequality~\eqref{EqnJackVarThreeFour}(iv):}

Introduce the shorthand notation $\OneDirScore(\Sam_i) \defn | \nabla
\TargetFun(\TruePar)^{\mytrans} \GradFun_\TruePar^{-1} \ZFun(\Sam_i,
\TruePar)|^2$. Then by Assumption~\ref{ass:tail} and properties of
Orlicz norm, we have $\sigma_{1/2} \defn \rnorm{\OneDirScore(\Sam_i) -
  \E[\OneDirScore(\Sam_i)]}{\psi_{1/2}} \leq \polyshort$.  Moreover,
\begin{align*}
v^2_{1/2} \defn \E[\OneDirScore(\Sam_i) - \E[\OneDirScore(\Sam_i)]]^2
\leq c \sigma^2_{1/2},
~~\text{and~~}
\E\Big[\Big|\frac{1}{\Numobs} \sum_{i=1}^\Numobs \OneDirScore(\Sam_i)
 - \E[\OneDirScore(\Sam_i)\Big|\Big]
 \leq 
 \frac{v_{1/2}}{\sqrt{\numobs}}.
\end{align*}
Invoking~\Cref{LemAdam} on the single-element function class
$\Fclass= \{\OneDirScore(\Sam_i) - \E[\OneDirScore(\Sam_i)] \}$, we
find
 \begin{align*}
 \Prob(\vterm_4\geq t+\frac{c v_{1/2}}{\sqrt{\numobs}})
 &\leq
 \Prob\Big[\big| \frac{1}{\Numobs} \sum_{i=1}^\Numobs \OneDirScore(\Sam_i)
 - \E[\OneDirScore(\Sam_i)] \big|\geq t+\frac{c v_{1/2}}{\sqrt{\numobs}}\Big]  \\
 &\leq 
\exp \big(- \frac{c\Numobs t^2}{ \sigma_{1/2}^2} \big) + 3 \exp \big( - \big\{ \frac{c\Numobs t}{ \sigma_{1/2}  \log^2\Numobs } \big\}^{1/2} \big).
\end{align*}
Letting $t=c\sigma_{1/2}( \sqrt{ \frac{ \log (1/\delta)}{\Numobs}} +  \frac{ \log^2(1/\delta) \log^2\Numobs}{\Numobs}) \leq  \polyshort \log^3(\Numobs/\delta) \frac{1}{ \sqrt{\Numobs}}$ yields inequality~\eqref{EqnJackVarThreeFour}(iv).


\section{A sufficient condition for~\ref{ass:convergence}}
\label{sec:suff_cond_of_conv}

In this appendix, we provide a sufficient condition for the
convergence condition~\ref{ass:convergence} to hold.

\myassumption{Con}{ass:concavity_app}{ The matrix $-\nabla \PopZfun(\Par)$
  is positive semidefinite for all $\Par \in \ParSpace$ and strictly
  positive definite for $\TruePar$---namely, we have $-\nabla
  \PopZfun(\TruePar) \succeq \strongconvex\IdMat_{\ParDim}$ for some
  constant $\strongconvex>0$.  }

 
 We now show that~\ref{ass:concavity} is a sufficient condition
 for~\ref{ass:convergence} when the sample size
 condition~\eqref{eq:linear_approx_Z_est_sample_ahead} is in force.

\begin{lemma}[A sufficient condition for~\ref{ass:convergence}]
\label{LemLinErrorBound_verify}
In addition to conditions~\ref{ass:tail}~and~\ref{ass:concavity},
suppose that the sample size satisfies $\Numobs/
  \log^{2}(\Numobs/\delta)  \geq \polyshort \myusedim$.  Then
\begin{align}
\label{eq:est_err_of_estpar_1_verify}  
\vecnorm{\EstPar_\numobs - \TruePar}{2} \leq c
\frac{\Longc}{\strongconvex} \sqrt{ \log \Numobs}
\sqrt{\frac{\myusedim}{\Numobs}}
\end{align}
with probability at least $1 - \delta$.
\end{lemma}
\noindent As a direct consequence, Assumption~\ref{ass:convergence} is
satisfied.


\subsection{Proof of~\Cref{LemLinErrorBound_verify}}
\label{proof:LemLinErrorBound_verify}

Our proof makes use of two auxiliary lemmas, which we begin by
stating.  We recall the notation $\PopZfun(\Par) = \Exs[\ZFun(\Sam,
  \Par)]$, and use $c, c'$ to denote universal positive constants.
\begin{lemma}
\label{lm:zfun_convexity}
Under conditions~\ref{ass:tail} and~\ref{ass:concavity}, for all $\Par
\in \ParSpace$, we have
\begin{align}
\label{eq:inequality-likelihood3}
\inprod{\PopZfun(\Par)}{\TruePar - \Par} & \geq
\begin{cases}
 \frac{\strongconvex}{2} \vecnorm{\Par - \TruePar}{2}^2, & \mbox{if
   $\vecnorm{\Par - \TruePar}{2} \leq
   \frac{\strongconvex}{2c\LipZFb}$, and} \\
\frac{\strongconvex^2}{4c\LipZFb} \vecnorm{\Par - \TruePar}{2}, &
\mbox{otherwise.}
\end{cases}
\end{align}
\end{lemma}
\noindent See~\Cref{sec:app-proof-zfun_convexity} for the proof.

\begin{lemma}
\label{lm:empirical_error}
Under condition~\ref{ass:tail}, there exists $\polyshort= \poly$ such
that given any $\delta \in (0, 1)$ and a sample size $\Numobs/
\log^{2} (\Numobs/\delta) \geq \polyshort \Ccerr$, we have
\begin{align}
\label{eq:bound-for-supem}
\SupEm \defn \sup_{\Par\in \ParSpace}
\vecnorm{\sum_{i=1}^\numobs\ZFun_\Numobs(\Sam_i, \Par)/\numobs -
  \PopZfun(\Par)}{2} \leq c \Longc \sqrt{ \log \Numobs} \sqrt{
  \frac{\ParDim + \log (1 / \delta)}{\Numobs}}
\end{align}
with probability $1 - \delta$.
\end{lemma}
\noindent See~\Cref{sec:app-proof-empirical_error} for the proof.\\

Using these auxiliary claims, we are equipped to prove the
claim~\eqref{eq:est_err_of_estpar_1_verify}
from~\Cref{LemLinErrorBound_verify}.  Introduce the shorthand
notation $\Par(s) \defn \TruePar + s(\EstPar_\numobs - \TruePar)$ and
$\bar\ZFun_\Numobs(\Par) \defn \sum_{i=1}^\Numobs\ZFun(\Sam_i,
\Par)/\Numobs$.  From the definition~\eqref{eq:z_estimate_1} of the
$Z$-estimator, we have $\bar\ZFun_\numobs(\EstPar_\numobs) = 0$,
whence
\begin{align}
\label{eq:linear_approx_Z_est1}
\inprod{\PopZfun(\EstPar_\numobs)}{\TruePar- \EstPar_\numobs} =
\inprod{\PopZfun(\EstPar_\numobs) -
  \bar\ZFun_\Numobs(\EstPar_\numobs)}{\TruePar - \EstPar_\numobs} \leq
\SupEm \cdot \vecnorm{\TruePar - \EstPar_\numobs }{2}.
\end{align}
Combining~\cref{eq:linear_approx_Z_est1} with~\Cref{lm:zfun_convexity}
yields the lower bound
\begin{align}
\label{eq:linear_approx_Z_est2}
\SupEm \cdot \vecnorm{\TruePar - \EstPar_\numobs }{2} & \ge
\begin{dcases}
  \frac{\strongconvex}{2} \vecnorm{\EstPar_\numobs - \TruePar}{2}^2, &
  \mbox{if $\vecnorm{\EstPar_\numobs - \TruePar}{2} \le
    \frac{\strongconvex}{2c\LipZFb}$, and} \\
 \frac{\strongconvex^2}{4c\LipZFb} \vecnorm{\EstPar_\numobs -
   \TruePar}{2}, & \mbox{otherwise.}
\end{dcases}
\end{align}
Now condition on the event that~\cref{eq:bound-for-supem} and the
sample size condition from~\Cref{LemLinErrorBound_verify} hold.  We
can then solve the fixed-point inequality for
$\vecnorm{\EstPar_\numobs - \TruePar}{2}$ to establish that
$\vecnorm{\EstPar_\numobs - \TruePar}{2} \leq c
\frac{\Longc}{\strongconvex} \sqrt{ \log \Numobs} \sqrt{ \frac{\ParDim
    + \log (1 / \delta)}{\Numobs}}$ with probability at least $1 -
\delta$.


\subsubsection{Proof of~\Cref{lm:zfun_convexity}} 
\label{sec:app-proof-zfun_convexity}

Define the function $\Par(s) \defn \TruePar + s(\Par - \TruePar)$ for
$s \in [0,1]$.  Performing Taylor series expansion of $\ZFun$ around
$\TruePar$ yields
\begin{align}
\inprod{\PopZfun(\Par)}{\TruePar - \Par} & = \inprod{\PopZfun(\Par) -
  \PopZfun(\TruePar)}{\TruePar - \Par} \notag \\
& = \binprod{\int_{0}^{1} \nabla  \PopZfun(\Par(s))(\Par -
  \TruePar) \textup{d}s}{\TruePar - \Par} \notag \\
\label{eq:zfun_convex_eq1}
& = (\Par - \TruePar)^{\mytrans} \nabla  \PopZfun(\TruePar)(\TruePar-\Par) +
(\Par - \TruePar)^{\mytrans} \Big( \int_{0}^{1} \nabla  \PopZfun(\Par(s))
\textup{d}s- \nabla  \PopZfun(\TruePar) \Big)(\Par - \TruePar).
\end{align}
By Assumption~\ref{ass:tail}~and~\ref{ass:concavity}, it follows that
for any $\Par$ satisfying $\vecnorm{\Par - \TruePar}{2} \le
\strongconvex/(2 c \LipZFb)$
\begin{subequations}
\begin{align}
\label{eq:inequality-likelihood2}
\binprod{\PopZfun(\Par)}{\TruePar - \Par} \ge
        {\strongconvex\vecnorm{\Par - \TruePar}{2}^2} - c\LipZFb
        \vecnorm{\Par - \TruePar}{2}^3 \geq
        \frac{\strongconvex}{2} \vecnorm{\Par - \TruePar}{2}^2.
\end{align}
When $\vecnorm{\Par - \TruePar}{2} \ge \strongconvex/(2c\LipZFb)$, let
$\Par = \TruePar + tv$ with $t = \vecnorm{\Par - \TruePar}{2}$ and
$v=(\Par - \TruePar)/\vecnorm{\Par - \TruePar}{2}$, we have
$\inprod{\PopZfun(\Par)}{ \frac{\TruePar - \Par}{\vecnorm{\TruePar -
      \Par}{2}}} = \inprod{\PopZfun(\TruePar + tv)}{-v}$. Taking the
derivative with respect to $t$, we find that \mbox{$\frac{d}{dt}
  \inprod{\PopZfun(\TruePar + tv)}{-v} = v^{\mytrans} \nabla  \PopZfun
  (\TruePar + tv)(-v) \geq 0$,} where the last inequality uses
Assumption~\ref{ass:concavity}. Therefore, we obtain the smallest
value of $\PopZfun(\TruePar + tv)^{\mytrans} (-v)$ when $t = \vecnorm{\Par -
  \TruePar}{2} = \strongconvex/(2c\LipZFb)$. Therefore, when
$\vecnorm{\Par - \TruePar}{2} \ge \strongconvex/(2c\LipZFb)$, by
equation~\eqref{eq:inequality-likelihood2}, we have
\begin{align}
\label{eq:beta-error-greater}
\inprod{\PopZfun(\Par)}{\TruePar - \Par} & \geq
\frac{\strongconvex^2}{ 4c\LipZFb} \vecnorm{\Par - \TruePar}{2}.
\end{align}
\end{subequations}
Combining equations~\eqref{eq:inequality-likelihood2}
and~\eqref{eq:beta-error-greater} concludes the proof of
Lemma~\ref{lm:zfun_convexity}.


\subsubsection{Proof of~\Cref{lm:empirical_error}}
\label{sec:app-proof-empirical_error}

Recall that $\SupEm \defn \sup_{\Par\in \ParSpace}
\vecnorm{\ZFun_\Numobs(\Par) - \PopZfun(\Par)}{2}$.
We have $\SupEm = \sup_{\Direc \in \Sphere{\ParDim-1}}
\SupEm(\Direc)$.
Let $\Coverset \defn \{\Direc_1, \ldots, \Direc_M\}$ be a $1/2$-covering of $\Sphere{\ParDim-1}$ in the Euclidean norm; such a covering exists with $M \leq 5^\ParDim$. Additionally, let $\ParSpace_0 \defn \{\Par^1, \ldots, \Par^N\}$ be an $\eps$-covering of $\ParSpace$ in the Euclidean norm for some $\eps > 0$, where $N \leq (1 + 2\ConParNorm/\eps)^\ParDim$~\cite[Example 5.8]{wainwright2019high}.
  
This proof consists of two steps: (a) a discretization argument to simplify the  problem to a finite maximum over a family of zero-mean random variables $Q_{\Direc,\Par}$ indexed by  $(\Direc,\Par)\in\Coverset\times\ParSpace_0$; (b) proving a tail bound on each $Q_{\Direc,\Par}$, and combined with a union bound on the cover set. More specifically, in the discretization step, we prove
\begin{subequations}
\begin{align}
    \SupEm 
    \leq 
    2 \max_{i \in [M]} \sup_{j \in [N]} \SupEm_{\Direc_i, \Par^j} 
    +
c \LipZFa \eps,\label{eq:sufficient_cond_concentration_1}
\end{align} with probability at least $1-\delta/2$, where $\SupEm_{\Direc, \Par} =   \frac{1}{\Numobs}
  \sum_{i=1}^\Numobs \Direc^{\mytrans} \ZFun(\Sam_i, \Par) - \E[\Direc^{\mytrans}
    \ZFun(\Sam, \Par)]$; and for each individual term $\SupEm_{\Direc,\Par}$, we show
\begin{align}
    \Prob[\SupEm_{\Direc_i, \Par^j} \geq t] \leq 2
\exp(-c\Numobs t^2/\OrPara^2) \label{eq:sufficient_cond_concentration_2}
\end{align}
for any $0 \leq t \leq c{\OrPara}$ for some
constant $c > 0$.
\end{subequations}
Combining these two results with the union bound, we obtain
\begin{align}
\SupEm \leq c\OrPara \sqrt{
  \frac{\ParDim \log (1 + 2\ConParNorm/\eps) + \log
    (1/\delta)}{\Numobs}} + c\OrParb\eps
\label{eq:disc_supem_5}
\end{align}
with probability at least $1 - \delta$. 

Finally, set $\eps=c
\sqrt{\dimdelta/\Numobs}$ for some constant $c > 0$. We can verify
that $t$ chosen in~\eqref{eq:sufficient_cond_concentration_2} is less than
$c\OrPara$. Therefore, we arrive at
\begin{align*}
\SupEm \leq c \Big[ \OrParb + \OrPara \sqrt{| \log (1 +
    \ConParNorm)|} + \OrPara \sqrt{ \log \Numobs } \Big] \sqrt{
  \frac{\dimdelta}{\Numobs}} \leq c \Longc \sqrt{ \log \Numobs} \sqrt{
  \frac{\dimdelta}{\Numobs}},
\end{align*} 
which completes the proof of~\Cref{lm:empirical_error}.\\

\myunder{Proof of the discretization bound~\eqref{eq:sufficient_cond_concentration_1}:}
define 
$ \SupEm(\Direc) \defn
  \sup_{\Par\in \ParSpace} \SupEm_{\Direc, \Par}.
$
 It then follows from a
standard discretization argument~\cite[Chap 6.]{wainwright2019high}
that
\begin{align}
\SupEm \leq  2\max_{\Direc_1, \ldots, \Direc_M} \SupEm(\Direc).\label{eq:disc_supem_15}
\end{align}
 Further discretizing $\SupEm(\Direc)$ yields
 \begin{align}
 \max_{\Direc_1, \ldots, \Direc_M} \SupEm(\Direc) 
 & \leq
 2\max_{i\in[M]} \big[\sup_{j\in[N]} \SupEm_{\Direc_i, \Par^j} +
   \sup_{\enorm{\Par_1 - \Par_2} \leq \eps }|\SupEm_{\Direc_i, \Par^1}
   - \SupEm_{\Direc_i, \Par^2}|\big] \notag \\
& 
\leq 2\max_{i\in[M]} \big[\sup_{j\in[N]} \SupEm_{\Direc_i, \Par^j}
   \notag \\
& 
\quad \quad + \sup_{\enorm{\Par_1 - \Par_2} \leq \eps }|
   \frac{1}{\Numobs} \sum_{i=1}^\Numobs \Direc^{\mytrans} (\ZFun(\Sam_i,
   \Par_1) - \ZFun(\Sam_i, \Par_2))| + \big|
   \E[\Direc^{\mytrans}(\ZFun(\Sam, \Par_1) - \ZFun(\Sam, \Par_2))] \big|
   \Big] \notag \\
\label{eq:disc_supem_2}    
& \leq 
2 \max_{i \in [M]} \sup_{j \in [N]} \SupEm_{\Direc_i, \Par^j} +
2 \sup_{\enorm{\Par_1 - \Par_2} \leq \eps } \enorm{ \frac{1}{\Numobs}
  \sum_{i=1}^\Numobs (\ZFun(\Sam_i, \Par_1) - \ZFun(\Sam_i, \Par_2))}
+ 2 \LipZFa \eps.
\end{align}
Moreover, by~\Cref{lm:sup_op_norm_zfun} and the sample size condition,
with probability at least $1 - \delta/2$
\begin{align}
\sup_{\enorm{\Par_1 - \Par_2} \leq \eps } \enorm{ \frac{1}{\Numobs}
  \sum_{i=1}^\Numobs (\ZFun(\Sam_i, \Par_1) - \ZFun(\Sam_i, \Par_2))}
& \leq \sup_{\enorm{\Par_1 - \Par_2} \leq \eps }
\sup_{\Par\in \ParSpace } \opnorm{ \frac{1}{\Numobs}
  \sum_{i=1}^\Numobs \nabla _\Par \ZFun(\Sam_i, \Par)} \enorm{\Par_1 -
  \Par_2} \notag \\
& \leq \Big[
c \OrParb+\polyshort\Big( \sqrt{\log(\numobs/\delta)} \sqrt{
  \frac{\dimdelta}{\Numobs}} + \frac{\log^2(\numobs/\delta) \dimdelta
}{\Numobs} \Big)\Big] 
\eps\notag\\
&\leq
c\OrParb\eps
\label{eq:disc_supem_4}
 \end{align}
for some constant $c$, $\polyshort= \poly>0.$ Putting pieces together yields the bound~\eqref{eq:sufficient_cond_concentration_1}.\\

\myunder{Proof of the tail bound~\eqref{eq:sufficient_cond_concentration_2}:}
Since $\rnorm{\Direc^{\mytrans} \ZFun(\Sam_i, \Par)}{\psi_1} \leq \OrPara$,
bound~\eqref{eq:sufficient_cond_concentration_2} follows directly from the concentration properties of sub-Exponential variables. 


\section{Verification of assumptions}
\label{SecVerify}

In this section, we collect the derivatives of the $Z$-functions in
linear regression, generalized linear models, inverse propensity score
weighting and weak instrumental variable examples introduced
in~\Cref{SecExamples} and verify the regularity conditions required
for our main results to hold.  Recall that for any $c>0$ and a random
variable $X$, we define the Orlicz norm
\begin{align*}
\rnorm{X}{\psi_{c}} \defn \inf\{ \lambda>0 \, \mid \,
\E[\psi_{c}(|X|/\lambda)] \leq 1 \}, ~~~~~~\text{where}~~~\psi_c(x)
\defn \exp(x^c)-1.
\end{align*}
For a random variable $Y \in \R$, we say it is sub-Gaussian with
parameter $\nu$ if $\rnorm{Y}{\psi_2} \leq \nu$. For a random vector $X
\in \R^\ParDim$, we say it is sub-Gaussian with parameter $\nu$ if it
is zero-mean and ${\inprod{\Direc}{X}}$ is sub-Gaussian with parameter
$\nu$ for any fixed unit vector $\Direc\in \sphere^{\ParDim-1}$.  See
the books~\cite{wainwright2019high,vershynin2018high} for some
standard properties of Orlicz norms.


\subsection{Linear regression}
\label{sec:linear_verify}

The linear regression problem from~\Cref{ExaLinearRegression} is
defined by the $Z$-function $\ZFun(\Sam_i, \Par) \defn \covar_i
\Zfunexam(Y_i- \inprod{\covar_i}{\Par})$. By some basic algebra, we
have
\begin{subequations}
\begin{align}
  \nabla _{\Par} \ZFun(\Sam_i, \Par) & = - \covar_i \covar_i^{\mytrans}
  \Zfunexam^\prime \big(Y_i- \inprod{\covar_i}{\Par}
  \big), \label{eq:linear_derivative_1} \\
  \nabla _{\theta}^2\ZFun(\Sam_i, \Par)[\Direc,v,w] & =
  \inprod{\covar_i}{\Direc} \inprod{\covar_i}{v} \inprod{\covar_i}{w}
  \Zfunexam^{\prime\prime}(Y_i- \inprod{\covar_i}{\Par}), \quad
  \mbox{and}
 \label{eq:linear_derivative_2} \\
\nabla _{\theta}^3 \ZFun(\Sam_i, \Par)[\Direc,v,w, \Direcd] & = -
\inprod{\covar_i}{\Direc} \inprod{\covar_i}{v} \inprod{\covar_i}{w}
\inprod{\covar_i}{\Direcd} \Zfunexam^{\prime\prime\prime} \big(Y_i-
\inprod{\covar_i}{\Par} \big) \label{eq:linear_derivative_3}
\end{align} 
\end{subequations}
for $\Direc, \Direcb, \Direcc, \Direcd\in \mathbb{S}^{\ParDim-1}$.  Next,
we introduce a set of assumptions on the linear regression problem
under which Assumption~\ref{ass:tail}~and~\ref{ass:concavity} are
satisfied.

\paragraph{Assumptions on linear regression.}
There exist constants $c_1,c_2,c_3,c_4, \ConParNorm>0$ such that
$|\Zfunexam(0)|,|\Zfunexam^\prime|,
|\Zfunexam^{\prime\prime}|,|\Zfunexam^{\prime\prime\prime}|$ are
bounded by $c_1$; $X$ is a sub-Gaussian vector with parameter $c_2$
and $Y$ is sub-Gaussian with parameter $c_3$;
$\Zfunexam^{\prime\prime} \geq c_4$ and
$\sigma_{\min}(\E[\covar_i\covar_i^{\mytrans}]) \geq c_4$; the parameter
space $\ParSpace$ is an open convex set with $\enorm{\Par} \leq
\ConParNorm$ for all $\Par\in \ParSpace$.

\paragraph{Verification of Assumption~\ref{ass:tail}.}
For any unit vectors $\Direc, \Direcb, \Direcc, \Direcd\in \Sphere{\ParDim-1}$,  we have
\begin{subequations}
\begin{align}
    |\inprod{\Direc}{\ZFun(\Sam_i, \Par)}| 
    & = 
    |\inprod{\Direc}{\covar_i}|\cdot|\Zfunexam(Y_i- \inprod{\covar_i}{\Par})|
    \leq 
     |\inprod{\Direc}{\covar_i}|\cdot(c_0+c_1|Y_i|+c_1\inprod{\covar_i}{\Par}) \label{eq:linear_verify_1}, \\
    | \nabla _{\Par} \ZFun(\Sam_i,
  \Par)[\Direc, \Direcb]|
  & = 
| \Zfunexam^\prime \big(Y_i- \inprod{\covar_i}{\Par} \big)|\cdot | \inprod{\covar_i}{\Direc}|\cdot|\inprod{\covar_i}{\Direcb}|\leq c_1| \inprod{\covar_i}{\Direc}|\cdot|\inprod{\covar_i}{\Direcb}|\label{eq:linear_verify_2}
 \end{align} and similarly
 \begin{align}
 \sup_{\Par\in \ParSpace}| \nabla ^2_{\Par} \ZFun(\Sam_i,
  \Par)[\Direc, \Direcb, \Direcc]|
  & = 
 \sup_{\Par\in \ParSpace}| \Zfunexam^{\prime\prime} \big(Y_i- \inprod{\covar_i}{\Par} \big)|\cdot | \inprod{\covar_i}{\Direc}|\cdot|\inprod{\covar_i}{\Direcb}|\cdot| \inprod{\covar_i}{\Direc}|\notag \\
 &\qquad
 \leq c_2| \inprod{\covar_i}{\Direc}|\cdot|\inprod{\covar_i}{\Direcb}|\cdot| \inprod{\covar_i}{\Direc}|\label{eq:linear_verify_3}, \\
  \sup_{\Par\in \ParSpace}| \nabla ^3_{\Par} \ZFun(\Sam_i,
  \Par)[\Direc, \Direcb, \Direcc, \Direcd]|
  & = 
 \sup_{\Par\in \ParSpace}| \Zfunexam^{\prime\prime\prime} \big(Y_i- \inprod{\covar_i}{\Par} \big)|\cdot | \inprod{\covar_i}{\Direc}|\cdot|\inprod{\covar_i}{\Direcb}|\cdot| \inprod{\covar_i}{\Direcc}|\cdot | \inprod{\covar_i}{\Direcd}|\notag \\
 &\qquad
 \leq c_3| \inprod{\covar_i}{\Direc}|\cdot|\inprod{\covar_i}{\Direcb}|\cdot| \inprod{\covar_i}{\Direcc}|\cdot | \inprod{\covar_i}{\Direcd}|\notag \\
  &\qquad =
 \Big[\frac{| \inprod{\covar_i}{\Direc}|}{c_2} \Big]\cdot \Big[c_3c_2|\inprod{\covar_i}{\Direcb}|\cdot| \inprod{\covar_i}{\Direcc}|\cdot | \inprod{\covar_i}{\Direcd}|\Big]\label{eq:linear_verify_4}.
\end{align}
\end{subequations}
Therefore, by the definition of sub-Gaussian random variables and
standard properties of Orlicz norms,we can verify that
condition~\ref{ass:tail} holds with some $\OrPara, \OrParb, \OrParc,
\OrPard > 0$ depending on $c_1, c_2, c_3, \ConParNorm$.

\paragraph{Verification of Assumption~\ref{ass:concavity}.}

From equation~\eqref{eq:linear_derivative_1} and the assumptions on
linear regression, we have
\begin{align*}
\nabla \PopZfun(\Par) = \E[\nabla \ZFun(\covar_i,
  \Par)] =-\E[\covar_i\covar_i^{\mytrans}
  \Zfunexam^\prime(Y_i-\inprod{\covar_i}{\Par})] \preceq
-c_4\E[\covar_i\covar_i^{\mytrans}] \preceq-c_4^2\IdMat_\ParDim
\end{align*}
for all $\Par\in \ParSpace$. Assumption~\ref{ass:concavity} is hence
satisfied.


\subsection{Generalized linear models}
\label{sec:glm_verify}

The generalized linear model problems from~\Cref{ExaGLM} have moment
equations with the $Z$-function $\ZFun(\Sam_i, \Par) = \SuffStat(Y_i)
\covar_i - \Zfunexam^\prime(\inprod{\covar_i}{\Par}) \covar_i$. The
associated derivatives are
\begin{align*}
&\nabla _{\theta} \ZFun(\Sam_i, \Par) 
= - \covar_i\covar_i^{\mytrans}
\Zfunexam^{\prime\prime}(\inprod{\covar_i}{\Par}), 
\\
& \nabla
_{\theta}^2\ZFun(\Sam_i, \Par)[\Direc,v,w] = -
\Zfunexam^{\prime\prime\prime}(\inprod{\covar_i}{\Par})
\inprod{\covar_i}{\Direc} \inprod{\covar_i}{v} \inprod{\covar_i}{w} ,
\\ 
& \nabla _{\theta}^3\ZFun(\Sam_i, \Par)[\Direc,v,w, \Direcd] = -
\Zfunexam^{\prime\prime\prime\prime}(\inprod{\covar_i}{\Par})
\inprod{\covar_i}{\Direc} \inprod{\covar_i}{v} \inprod{\covar_i}{w}
\inprod{\covar_i}{\Direcd}
\end{align*}
for $\Direc, \Direcb, \Direcc, \Direcd\in \mathbb{S}^{\ParDim-1}$.  We
show that Assumption~\ref{ass:tail}~and~\ref{ass:concavity} are
satisfied under the following set of assumptions.

\paragraph{Assumptions on generalized linear models.}

There exist some constants $c_1,c_2,c_3, \ConParNorm>0$ such that
$|\Zfunexam^\prime(0)|,
|\Zfunexam^{\prime\prime}|,|\Zfunexam^{\prime\prime\prime}|,|\Zfunexam^{\prime\prime\prime\prime}|$
are bounded by $c_1$; $X$ is a sub-Gaussian vector with parameter
$c_2$ and $\SuffStat(Y)$ is sub-Gaussian with parameter $c_3$; the
parameter space $\ParSpace$ is an open convex set with $\enorm{\Par}
\leq \ConParNorm$ for all $\Par\in \ParSpace$. Moreover, there exists
some constant $c_4>0$ such that $\Zfunexam^{\prime\prime} \geq 0$
and \begin{align*} \sigma_{\min}(\E[\indicator_{\Event}
    \covar_1\covar_1^{\mytrans}]) \geq c_4,
\end{align*}
where 
\begin{align*} \Event\defn \{\covar_1: |\inprod{\covar_1}{\TruePar}|\leq c_5 \} \text{ and }c_5\defn \inf_{c\geq0} \{\min\{\Zfunexam^{\prime\prime}(c), \Zfunexam^{\prime\prime}(-c) \} \leq c_4\}.
\end{align*}

\paragraph{Verification of Assumption~\ref{ass:tail}.}
We can verify this assumption using the derivatives of $\ZFun$ and
following the same arguments as for the linear regression problem in
Section~\ref{sec:linear_verify}. Thus we omit the proof here.
\paragraph{Verification of Assumption~\ref{ass:concavity}.}
Note that 
\begin{align*}
\PopZfun(\Par) & = \E[\nabla  \ZFun(\Sam_i, \Par)] =
-\E\Big[{\Zfunexam^{\prime\prime}(\inprod{\covar_i}{\Par})} \cdot
  \covar_i\covar_i^{\mytrans} \Big]\preceq \mathbf{0} \text{ for all
} \Par\in \ParSpace,~~~~~ \text{and} \\
\PopZfun(\TruePar) & = \E[\nabla  \ZFun(\Sam_i, \TruePar)] = - \E
\Big[{\Zfunexam^{\prime\prime}(\inprod{\covar_i}{\TruePar})} \cdot
  \covar_i\covar_i^{\mytrans} \Big] \preceq -\E \Big
    [\indicator_{\Event}{\Zfunexam^{\prime\prime}(\inprod{\covar_i}{\TruePar})} \cdot
      \covar_i\covar_i^{\mytrans} \Big]\\ &\leq -c_4 \E
    [\indicator_{\Event} \covar_1\covar_1^{\mytrans}] \preceq -c_4^2
    \IdMat_\ParDim.
\end{align*}

\paragraph{Logistic regression satisfies the boundedness assumption
  on $\Zfunexam$.}

For logistic regression, we have $\Zfunexam(t) = \log(1 + e^t)$, and
hence
\begin{align*}
\Zfunexam^\prime(x) = \frac{e^x}{1 + e^x},~~~
\Zfunexam^{\prime\prime}(x) = \frac{e^x}{(1 + e^x)^2},~~~
\Zfunexam^{\prime\prime\prime}(x) = \frac{e^x (1 - e^{x})}{(1 +
  e^x)^3}, \quad \Zfunexam^{\prime\prime\prime}(x) =
\frac{e^{3x}-4e^{2x}+e^x }{(1 + e^x)^4}.
\end{align*}
It can be readily verified that $|\Zfunexam^\prime(0)|,
|\Zfunexam^{\prime\prime}|,|\Zfunexam^{\prime\prime\prime}|,|\Zfunexam^{\prime
  \prime \prime \prime} |$ are all bounded by some universal constant
$c_1>0$ for all $x$, and $\Zfunexam^{\prime\prime} \geq 0$.


\subsection{Inverse propensity score weighting}\label{sec:ipw_verify}
Recall the definitions and notations in~\Cref{ExaIPW}. 
The inverse propensity score weighting (IPW) problem has the $Z$-function 
\begin{align*}
\ZFun(\Sam_i, \Par) =
\begin{bmatrix}
\covar_i \big(\act_i - \frac{\exp(\inprod{\covar_i}{\beta})}{1 +
  \exp(\inprod{\covar_i}{\beta})}\big)\\ \frac{\Action_i \Outcome_i}{
  \logistic{\State_i}{\beta}} - \frac{(1 - \Action_i) \Outcome_i}{1 -
  \logistic{\State_i}{\beta}} - \tau
\end{bmatrix}.
\end{align*}
We use $u_{1:\ParDim-1}$ to denote the first $\ParDim-1$ coordinates of
$u$ and write $\Par=(\beta,\tau)$. Adopt the shorthand notation $\smoothtrun$ for $\inprod{\covar_i}{\beta}$.  By some basic algebra, we have  
\begin{align*}
 & \nabla_{\theta} \ZFun(\Sam_i, \Par) =
\begin{bmatrix}
- \covar_i \frac{ e^{\smoothtrun}}{(1 + e^{\smoothtrun})^2}\covar_i^{\mytrans}  & 0_{\ParDim-1}\\
 - \Action_i \Outcome_i \covar_i^{\mytrans} e^{-\smoothtrun}
-(1 - \Action_i) \Outcome_i \covar_i^{\mytrans} e^{\smoothtrun}  & -1
\end{bmatrix}
\\
 & \nabla_{\theta}^2\ZFun(\Sam_i, \Par)[u,v,w] = -  \frac{ e^{\smoothtrun}(1 -  e^{\smoothtrun})}{(1 +  e^{\smoothtrun})^3}\inprod{\covar_i}{u_{1:\ParDim-1}}
\inprod{\covar_i}{v_{1:\ParDim-1}}\inprod{\covar_i}{w_{1:\ParDim-1}}
\\
 &\quad\quad   +  u_{\ParDim }\Big(
\Action_i \Outcome_i  e^{-\smoothtrun}
-(1 - \Action_i) \Outcome_i  e^{\smoothtrun}\Big) \inprod{\covar_i}{v_{1:\ParDim-1}} \inprod{\covar_i}{w_{1:\ParDim-1}} 
, \\
 & \nabla_{\theta}^3\ZFun(\Sam_i, \Par)[u,v,w,\Direcd]\\
&=-  \Big[
  \frac{e^{s}-4e^{2s}+e^{3s} }{(1 + e^x)^4}
\Big]\inprod{\covar_i}{u_{1:\ParDim-1}}
\inprod{\covar_i}{v_{1:\ParDim-1}}\inprod{\covar_i}{w_{1:\ParDim-1}}\inprod{\covar_i}{\Direcd_{1:\ParDim-1}}\\
 &~   +  u_{\ParDim }\Big(
 - \Action_i \Outcome_i  e^{-\smoothtrun}
-(1 - \Action_i) \Outcome_i   e^{\smoothtrun} \Big)\inprod{\covar_i}{v_{1:\ParDim-1}} \inprod{\covar_i}{w_{1:\ParDim-1}} \inprod{\covar_i}{\Direcd_{1:\ParDim-1}} 
\end{align*}
for $\Direc,v,w,\Direcd\in \mathbb{S}^{\ParDim-1}$.
We make the following set of assumptions on the IPW problem.

\paragraph{Assumptions on the inverse propensity weighting problem.}
There exist some constants  $c_1,c_2,\ConParNorm,c_3,c_4>0$ such that   $X$ is a sub-Gaussian vector with parameter $c_1$ and $Y$  is  sub-Gaussian with parameter $c_2$; the covariance $\E[\covar\covar^{\mytrans}]\succeq c_3\IdMat_{\ParDim-1}$ the parameter space $\betaspace$ is an open convex set with $\enorm{\beta}\leq\ConParNorm$ for all $\beta\in\betaspace$; a strong overlap assumption on the propensity score, namely, $\varphi(\inprod{\beta}{\covar} \big)\in[c_4,1-c_4]$ for all $X$ and $\beta\in\betaspace$.
It is important to note that by the boundedness of propensity score and the sub-Gaussianity of $Y_i$, we have
\begin{align*}
    |\taustar| \leq \E\Big[\Big|\frac{\Action_i \Outcome_i}{
  \logistic{\State_i}{\beta}} - \frac{(1 - \Action_i) \Outcome_i}{1 -
  \logistic{\State_i}{\beta}}\Big|\Big]< \widetilde{c}
\end{align*} for  some $(c_2,c_3)$-dependent constant $\widetilde{c}$.  Therefore, we may w.l.o.g. assume $\ParSpace \defn \betaspace\times (-\widetilde{c},\widetilde{c})$, which is a bounded convex open set in $\R^{\ParDim}$.

\paragraph{Verification of Assumption~\ref{ass:tail}.}
Note that the functions \begin{align*}
   \frac{e^x}{1+e^x},~~~ \frac{e^x}{(1+e^x)^2},~~~
   \frac{e^x (1 - e^{x})}{(1 + e^x)^3},~~~ 
    \frac{e^{3x}-4e^{2x}+e^x }{(1 + e^x)^4}
\end{align*} are all bounded by some universal constant. Therefore, it follows 
 the same arguments as for the linear regression problem in Section~\ref{sec:linear_verify} that Assumption~\ref{ass:tail} is satisfied. \\

 
 Although the IPW estimator does not satisfy~\ref{ass:concavity} since $\nabla \PopZfun(\TruePar)$ is non-concave, we will prove that it satisfies the convergence assumption~\ref{ass:convergence}.
\paragraph{Verification of Assumption~\ref{ass:convergence}.}
First, by definition
\begin{align*}
\PopZfun(\TruePar) =\E[ \nabla_{\theta} \ZFun(\Sam_i, \TruePar)] =
\begin{bmatrix}
- \E[ \covar_i \frac{ e^{\smoothtrun}}{(1 + e^{\smoothtrun})^2}\covar_i^{\mytrans} ] & 0_{\ParDim-1}\\
 - \E[\Action_i \Outcome_i \covar_i^{\mytrans}  e^{-\smoothtrun}
+(1 - \Action_i) \Outcome_i \covar_i^{\mytrans} e^{\smoothtrun}]  & -1
\end{bmatrix}=:\begin{bmatrix}
    A & 0_{\ParDim-1}\\
    b^{\mytrans} & -1
\end{bmatrix}.
\end{align*}
Note that under the assumptions, we have
$-A = \E[ \covar_i \frac{ e^{\smoothtrun}}{(1 + e^{\smoothtrun})^2}\covar_i^{\mytrans} ] \succeq \polyshort\IdMat_{\ParDim-1}$ and $ bb^{\mytrans}\preceq  \polyshortprime\IdMat_{\ParDim-1}$ for some  $(c_1,c_2,c_3,c_4,\ConParNorm)$-dependent constants $\polyshort,\polyshortprime>0$. Therefore, it can be verfied that the first part of Assumption~\ref{ass:convergence} is satisfied for some $(c_1,c_2,c_3,c_4,\ConParNorm)$-dependent parameter $\strongconvex>0$.

To establish the second part of Assumption~\ref{ass:convergence}, we proceed by first bounding the estimation error $\enorm{\betahat-\betastar}$ of logistic regression, and then transfering it into a bound on $|\tauhatipw-\taustar|$. Namely, we show
\begin{subequations}
    \begin{align}
    \enorm{\betahat-\betastar}
    &\leq     
    \polyshort \sqrt{ \log \Numobs} \sqrt{
  \frac{\myusedim}{\Numobs}}, \label{eq:ipw_proof_step1}\\
   |\tauhatipw-\taustar|
    &\leq     
    \polyshort \sqrt{ \log \Numobs} \sqrt{
  \frac{\myusedim}{\Numobs}} \label{eq:ipw_proof_step2}
    \end{align}
    with probability at least $1-\delta$ for some $(c_1,c_2,c_3,c_4,\ConParNorm)$-dependent constant $\polyshort>0$.  Putting two results together yields the second part of Assumption~\ref{ass:convergence}.
\end{subequations}\\

\myunder{Proof of the bound on $ \enorm{\betahat-\betastar}$: }
note that $\betahat$ is the MLE of the logistic regression problem consisting of the first $\ParDim-1$-dimension of $\ZFun$.
Applying Lemma~\ref{LemLinErrorBound_verify} to this  problem  (the assumptions for Lemma~\ref{LemLinErrorBound_verify} can be verified following the argument in~\Cref{sec:glm_verify}), we have
 \begin{align}\vecnorm{\betahat -
  \beta}{2}\leq 
\polyshort \sqrt{ \log \Numobs} \sqrt{
  \frac{\myusedim}{\Numobs}}\label{eq:ipw_verify_pf_2}
  \end{align}
  with probability at least $1-\delta$.\\

\myunder{Proof of the bound on $ |\tauhatipw-\taustar|$: }
by~\Cref{lm:empirical_error}, we have with probability at least $1-\delta$
\begin{align*}
\abss{\Exs[\ZFun_\ParDim(\Sam_i, \betahat ,\tauhatipw)]} \le \polyshort\sqrt{ \log
  \Numobs} \sqrt{ \frac{\myusedim}{\Numobs}}.
\end{align*}
The expectation here is only taken with respect to $\Sam_i$, but not
$\betahat $ and $\tauhatipw$.  By a Taylor series expansion, we have
with probability at least $1-\delta$
\begin{align}
&\Bigg|\Exs\big[\int_{0}^1 \nabla_\beta{\ZFun_\ParDim}(\Sam_i, \betastar+t(\betatil-\betastar), \taustar+t(\tautil-\taustar))d t
    \big] (\betahat - \beta)\notag\\
    &+\qquad
    \Exs\big[\int_{0}^1 \nabla_\tau{\ZFun_\ParDim}(\Sam_i, \betastar+t(\betatil-\betastar), \taustar+t(\tautil-\taustar))d t
    \big] (\tauhat - \taustar)  \Bigg|\leq \polyshort \sqrt{ \log
  \Numobs} \sqrt{ \frac{\myusedim}{\Numobs}}.\label{eq:ipw_verify_pf_1}
\end{align}
Note that $ \nabla_\tau{\ZFun}_{\ParDim }(\Sam_i, \betatil, \tautil ) =(0_{\ParDim-1},1)^{\mytrans}$
 for all $\betatil,\tautil$.  
  Moreover,
  \begin{align}
  \enorm{\Exs\big[\int_{0}^1 \nabla_\beta{\ZFun_\ParDim}(\Sam_i, \betastar+t(\betatil-\betastar), \taustar+t(\tautil-\taustar))d t
    \big]}
    \leq 
    \sup_{\Par\in\ParSpace}\opnorm{\E[\nabla_\Par \ZFun(\Sam_i,\Par)]}\leq c \OrParb
    \label{eq:ipw_verify_pf_3}
  \end{align} by Assumption~\ref{ass:tail}. Combining equation~\eqref{eq:ipw_verify_pf_2},~\eqref{eq:ipw_verify_pf_1}~and~\eqref{eq:ipw_verify_pf_3} gives the desired bound on $ |\tauhatipw-\taustar|$.

\section{Empirical process bounds}
\label{SecEmpProcess} 

In this section, we state and prove some results on the suprema of
empirical processes that are used in the proofs of both our main
theorems.  Both results are stated in terms of a universal constant $c
> 0$; some parameter-dependent constant $\polyshort = \poly > 0$; and
using the shorthand $\dimdelta = \ParDim + \log(1/\delta)$. 

\begin{lemma}[Concentration  of operator norm]
\label{lm:op_norm_zfun}
Under condition~\ref{ass:tail}, when the sample size $\numobs\geq \polyshortprime\cdot \Ccerr$,  for any $\Par \in \ParSpace$ and unit vector $\Direc\in\Sphere{\ParDim-1}$, we have
\begin{subequations}
\begin{align}
\label{eq:op_norm_zfun_claim1}    
\opnorm{ \frac{1}{\Numobs} \sum_{i=1}^\Numobs \nabla_\Par
  \ZFun(\Sam_i, \Par) - \E[ \nabla_\Par \ZFun(\Sam_i, \Par)]} 
  & \leq c\OrParb
 \Big( \sqrt{ \frac{\Ccerr}{\Numobs}} + \frac{\dimdelta \log\numobs}{\Numobs} \Big), \\
\label{eq:op_norm_zfun_claim2}
\opnorm{ \frac{1}{\Numobs} \sum_{i=1}^\Numobs \inprod{\Direc}{\nabla^2_\Par
  \ZFun(\Sam_i, \Par)} - \E[ \inprod{\Direc}{\nabla^2_\Par
  \ZFun(\Sam_i, \Par)}]} 
  & \leq
\polyshort\Big( \sqrt{ \frac{\Ccerr\log(\numobs/\delta)}{\Numobs}} + \frac{\dimdelta \log^{5/2}
  (\Numobs/\delta)}{\Numobs} \Big),
 \end{align}
\end{subequations}
with probability at least $1 - \delta$.
\end{lemma}
\noindent See~\Cref{sec:app-proof_op_norm_zfun} for the proof.

\begin{lemma}[Supremum of empirical operator norm]
\label{lm:sup_op_norm_zfun}
Under condition~\ref{ass:tail}, for any fixed unit vector $\Direc \in
\sphere^{\ParDim-1}$, we have
\begin{subequations}
\begin{align}
\label{eq:bounded_grad_zfun_lm}   
\sup_{\Par \in \ParSpace} \opnorm{ \frac{1}{\Numobs}
  \sum_{i=1}^\Numobs \nabla_\Par \ZFun(\Sam_i, \Par)} & \leq c \OrParb
+ \polyshort \Big \{ \sqrt{ \frac{\dimdelta \log(2
    \numobs/\delta)}{\Numobs}} + \frac{\dimdelta \log^2(2
  \numobs/\delta)}{\numobs} \Big \}, \\
\label{eq:bounded_hessian_zfun_lm}
\sup_{\Par \in \ParSpace} \opnorm{ \frac{1}{\Numobs} \sum_{i=1}^\Numobs
  \nabla^2_\Par \ZFun(\Sam_i, \Par)} & \leq c \OrParc\Big(1 + \sqrt{
  \frac{\dimdelta}{\Numobs}} + \frac{(\dimdelta \, (1+\log
  \Numobs))^{\InvOrConc}}{\Numobs} \Big), \\
\label{eq:bounded_third_zfun_lm}
\sup_{\Par \in \ParSpace} \opnorm{ \frac{1}{\Numobs}
  \sum_{i=1}^\Numobs \nabla^3_\Par \ZFun(\Sam_i, \Par)[\Direc]} & \leq
c \OrPard \sqrt{ \log (\Numobs/\delta)} \Big (1 + \sqrt{
  \frac{\dimdelta}{\Numobs}} + \frac{(\dimdelta \, (1+\log
  \Numobs))^{3/2}}{\Numobs} \Big),
\end{align}
\end{subequations}
with probability at least $1 - \delta$.
\end{lemma}
\noindent  In this work, we apply Lemma~\ref{lm:sup_op_norm_zfun} with $\numobs$ being either sufficiently large or equal to $1.$ See~\Cref{sec:app-proof-sup_op_norm_zfun} for the proof of this Lemma.
\\

Our proof of these lemmas make use of a known result on the suprema
of empirical processes, which we state here.  Let $\Fclass$ be a class
of measurable functions $f: \mathcal{X} \rightarrow \real$, zero-mean
under some distribution $\Prob$.  Define the quantities
\begin{align*}
\gamma & \defn \big \|\sup _{f \in \mathcal{F}} |f(X)| \big
\|_{\psi_\alpha} \quad \mbox{and} \quad v^2 \defn \sup _{f \in \mathcal{F}}
\Exs[f^2(X)],
\end{align*}
and given an i.i.d. sequence $\{X_i\}_{i=1}^\numobs$, define the
random variable $Z_\numobs(\Fclass) \defn \sup_{f \in \mathcal{F}}
\left| \frac{1}{\Numobs} \sum_{i=1}^\Numobs f(X_i) \right|$.  With
this set-up, we have the following simplification of Theorem 4 in the
paper~\cite{adamczak2008tail}:
\begin{lemma}
\label{LemAdam}
For any $\alpha>0$,
there is an $\alpha$-dependent constant $c$ such that
\begin{align}
\label{eq:adamczak-concentration}
\Prob \Big[ Z_\numobs(\Fclass) \geq 2 \Exs[Z_\numobs(\Fclass)] + t
  \Big] & \leq \exp \left(- \frac{\Numobs t^2}{4 v^2} \right) + 3 \exp
\left( - \left\{ \frac{c \Numobs t}{ \gamma ( 1 + \log^{1/\alpha} \Numobs)}
\right\}^{\alpha} \right)
\end{align}
\mbox{for any $t > 0$.}
\end{lemma}

The original statement in Adamczak's paper~\cite{adamczak2008tail} has
the quantity $\vecnorm{\max \limits_{i=1, \ldots, n} \sup_{f \in
    \funcClass} |f (\State_i)|}{\psi_\alpha}$ as opposed to the $\gamma (1
+ \log^{1/\alpha} \numobs)$ term in~\cref{eq:adamczak-concentration}. Our
simplified statement follows by applying Pisier's
inequality~\cite{pisier1983some} to the random variables $Y_i =
\sup_{f \in \Fclass} |f(X_i)|$, which guarantees that
\begin{align*}
\vecnorm{\max_{i = 1, \ldots, \numobs} \Outcome_i}{\psi_\alpha} \leq c \big
\{ 1 + \log^{1/\alpha} \numobs \big \} \cdot \max_{i = 1, \ldots, \numobs}
\vecnorm{\Outcome_i}{\psi_\alpha} \; = \; c(1 + \log^{1/\alpha} \numobs) \gamma
\end{align*}for some $\alpha$-dependent constant $c$.


\subsection{Proof of~\Cref{lm:op_norm_zfun}}
\label{sec:app-proof_op_norm_zfun}

\paragraph{Proof of claim~\eqref{eq:op_norm_zfun_claim1}.}
The proof of claim~\eqref{eq:op_norm_zfun_claim1} makes use of a standard discretization method (e.g., see
Chapter 6 in the book~\cite{wainwright2019high} for arguments of this
type); it also makes use of various properties of Orlicz quasi-norm.
Throughout this argument (as well as the proof
of~\Cref{lm:sup_op_norm_zfun} ), we use the standard fact that there
exists a $1/8$-covering set $\Coverset = \{u^1, \ldots, u^M \}$ of the
$\usedim$-dimensional Euclidean unit sphere with cardinality $N \defn
|\Coverset| \leq 17^\usedim$.

At a high-level, the proof consists of two steps: (a) a discretization
argument to reduce the problem to a finite maximum over a family of
zero-mean random variables $Y(u,v)$ indexed by pairs $(u,v)$ of
elements in the cover $\Coverset$; and (b) proving a tail bound on
each $Y(u, v)$, combined with the union bound over the cover.  More precisely, in the discretization step, we prove that
\begin{subequations}
\begin{align}
  \label{EqnDsimple}
  \opnorm{ \frac{1}{\Numobs} \sum_{i=1}^\Numobs \nabla_\Par
  \ZFun(\Sam_i, \Par) - \E[ \nabla_\Par \ZFun(\Sam_i, \Par)]} & \leq 2
\max_{u, v \in \Coverset} |Y(u,v)|,
\end{align}
where $Y(u,v) \defn \frac{1}{\Numobs} \sum_{i=1}^\Numobs \nabla_\Par
\ZFun(\Sam_i, \Par)[\Direc, \Direcb] - \E[ \nabla_\Par \ZFun(\Sam_i,
  \Par)[\Direc, \Direcb]]$.  In the tail bound step, we argue that
\begin{align}
\label{EqnTsimple}  
\Prob \big[ |Y(u,v)| \geq t+\frac{c\OrParb}{\sqrt{\numobs}} \big] & \leq \exp\big(- \frac{c\Numobs t^2}{
  \OrParb^2} \big) + 3 \exp \Big( -  \frac{c\Numobs t}{
  \sigma_{\OrConb} ( \log \Numobs)}  \Big).
\end{align}
\end{subequations}
Combining these two results with the union bound yields the
claim~\eqref{eq:op_norm_zfun_claim1}. \\

\myunder{Proof of the discretization bound~\eqref{EqnDsimple}:}
By the variational definition of the operator norm, we have
\begin{align}
\label{EqnVarSimple}  
T \defn \opnorm{ \frac{1}{\Numobs} \sum_{i=1}^\Numobs \nabla_\Par
  \ZFun(\Sam_i, \Par) - \E[ \nabla_\Par \ZFun(\Sam_i, \Par)]} & =
\sup_{\|u\|_2 = \|v\|_2 = 1} Y(u,v).
\end{align}
By the definition of the $1/8$-cover, for any pair $(u,v)$ of
unit-norm vectors, there exist $u^a, u^b \in \Coverset$ such that
$\max \{ \|\Direca - u^a\|_2, \; \|\Direcb - u^b\|_2, \} \leq
\frac{1}{8}$.  Introducing the shorthand $\Delta^a = u - u^a$, we then
have
\begin{align*}
|Y(u, v)| & = |Y(u^a +\Delta^a, v)| \leq |Y(u^a, v)| + |Y(\Delta^a,
v)| \; \leq \; |Y(u^a, v)| + \tfrac{1}{8} T,
\end{align*}
where the second inequality follows from the fact that $|Y(\Delta^a,
v)| \leq \|\Delta^a\|_2 T$ from the definition~\eqref{EqnVarSimple},
along with the bound $\|\Delta^a\|_2 \leq 1/8$.  A similar argument
applied to $v$ yields $|Y(u,v)| \leq |Y(u^a, u^b)| + \tfrac{1}{4} T$.
Taking suprema and maxima appropriately yields
\begin{align*}
T \; = \; \sup_{u,v} |Y(u,v)| & \leq \max_{u^a, u^b} |Y(u^a, v^b)| +
\tfrac{1}{4} T,
\end{align*}
and re-arranging implies the claimed bound. \\

\myunder{Proof of the tail bound~\eqref{EqnTsimple}:}
Let $Y_0(u,v)\defn \nabla_\Par \ZFun(\Sam, \Par)
  [\Direc, \Direcb] - \E[ \nabla_\Par \ZFun(\Sam, \Par)[\Direc,
      \Direcb]]$. 
Condition~\ref{ass:tail} guarantees that
\begin{align*}
\rnorm{|Y_0(u,v)|}{\psi_1} = \rnorm{| \nabla_\Par \ZFun(\Sam, \Par)
  [\Direc, \Direcb] - \E[ \nabla_\Par \ZFun(\Sam, \Par)[\Direc,
      \Direcb]]|}{\psi_1} & \leq c \OrParb,
\end{align*}
and similarly, we have $\var(Y_0(u,v)) \leq c \OrParb^2$. Moreover, $\E[|Y(u,v)|]\leq\sqrt{\var(Y_0(u,v))}/\sqrt{\numobs}$.   Consequently,
the tail bound~\eqref{EqnTsimple} follows by
applying~\Cref{LemAdam} to the single element function
class $\Funclass = \{Y_0(u,v) \}$.

\paragraph{Proof of claim~\eqref{eq:op_norm_zfun_claim2}.}
The proof of claim~\eqref{eq:op_norm_zfun_claim2} follows from a truncation argument under the matrix operator norm and an application of the matrix Bernstein inequality (see e.g., exercise 6.10 in Wainwright~\cite{wainwright2019high}) to the truncated matrices. We note that the same technique can be used to prove claim~\eqref{eq:op_norm_zfun_claim1}, albeit with additional logarithmic factors.


  More precisely, in the truncation step, letting  $\newtruncate_A{(\Sam_i, \Par)}\defn \nabla^2_\Par
  \ZFun(\Sam_i, \Par)[\Direc]\mathbf{1}_{\{\opnorm{\nabla^2_\Par
  \ZFun(\Sam_i, \Par)[\Direc]}\leq A\}}$ be the truncated matrix with $A = \polyshort\Ccerr\cdot\log^{\InvOrConc}(2\numobs/\delta) $, we prove that 
\begin{subequations}
\begin{align}
  \label{EqnTruncate1}
 \newtruncate_A{(\Sam_i, \Par)} =  \nabla^2_\Par
  \ZFun(\Sam_i, \Par)[\Direc]
\end{align} with probability at least $1-\delta/2$ for all $i\in[\numobs]$, and 
\begin{align}
  \label{EqnTruncate2}
 |\E[\newtruncate_A{(\Sam_i, \Par)}]-\E[\nabla^2_\Par
  \ZFun(\Sam_i, \Par)[\Direc]]| \leq \frac{\polyshort}{\numobs}.  
\end{align}

    We then apply the matrix Bernstein inequality to the i.i.d. matrices $\newtruncate_A(\Sam_i,\Par),i\in[\numobs]$ to obtain
\begin{align}
\label{EqnTruncate3}  
\Prob(\opnorm{ \frac{1}{\Numobs} \sum_{i=1}^\Numobs \newtruncate_A
  (\Sam_i, \Par) - \E[ \newtruncate_A (\Sam_i, \Par)]}\geq t) 
  & \leq c \ParDim\cdot  \exp(-\frac{c\numobs t^2}{c\OrPar^2\ParDim+At})
\end{align}
\end{subequations}

Combining these three results with the union bound and the triangle inequality yields the
claim~\eqref{eq:op_norm_zfun_claim1}. \\

\myunder{Proof of the truncation bound~\eqref{EqnTruncate1}:}
It suffices to show 
\begin{align}
     \opnorm{\nabla^2_\Par
  \ZFun(\Sam_i, \Par)[\Direc]}
    \leq A\label{eq:EqnTruncate1_equiv}
\end{align} with probability at least $1-\delta/2.$

By some basic algebra, we have
\begin{align*}
    \opnorm{\nabla^2_\Par
  \ZFun(\Sam_i, \Par)[\Direc]}
    \leq
    \fronorm{\nabla^2_\Par
  \ZFun(\Sam_i, \Par)[\Direc]}
    =
\sqrt{\sum_{i,j=1}^\ParDim\nabla^2_\Par \ZFun(\Sam_i,\Par)[\Direc,e_i,e_j]^2}
    \leq \ParDim\cdot\max_{i,j\in[\numobs]}|\nabla^2_\Par\ZFun(\Sam_i,\Par)[\Direc,e_i,e_j]|,
\end{align*} Using Assumption~\ref{ass:tail}, properties of sub-Exponential variables and the union bound, we further obtain
\begin{align*}
\Prob(\max_{i,j\in[\numobs]}|\nabla^2_\Par\ZFun(\Sam_i,\Par)[\Direc,e_i,e_j]|\geq t)\leq c\ParDim\cdot\exp\big(-{\polyshort t^{\OrConc}}).
\end{align*}
Combining the last two displays gives
\begin{align}
    \Prob( \opnorm{\nabla^2_\Par \ZFun(\Sam_i,\Par)[\Direc]}\geq t)\leq c\ParDim\cdot\exp\big(-{\polyshort t^\OrConc})\label{eq:EqnTruncate1_equiv_2}
\end{align}
Substituting $t=A$ into equation~\eqref{eq:EqnTruncate1_equiv_2} yields equation~\eqref{eq:EqnTruncate1_equiv}.\\

\myunder{Proof of the truncation bound~\eqref{EqnTruncate2}:}
by the definition of $\newtruncate_A$ and Jensen's inequality, we have
\begin{align*}
\opnorm{\E[\newtruncate_A{(\Sam_i, \Par)}]-\E[\nabla^2_\Par
  \ZFun(\Sam_i, \Par)[\Direc]]} 
  &= \opnorm{\E [\nabla^2_\Par
  \ZFun(\Sam_i, \Par)[\Direc]]\mathbf{1}_{\{\opnorm{\nabla^2_\Par
  \ZFun(\Sam_i, \Par)[\Direc]}>A\}} }\notag
  \\
  &\leq 
  \E [\opnorm{\nabla^2_\Par
  \ZFun(\Sam_i, \Par)[\Direc]}]\mathbf{1}_{\{\opnorm{\nabla^2_\Par
  \ZFun(\Sam_i, \Par)[\Direc]\}}>A}]\\
  &= 
  \int_{A}^\infty \Prob(\opnorm{\nabla^2_\Par
  \ZFun(\Sam_i, \Par)[\Direc]}\geq t) dt.  
\end{align*}
Combining this bound with  equation~\eqref{eq:EqnTruncate1_equiv_2} yields equation~\eqref{EqnTruncate2}. \\

\myunder{Proof of the matrix Bernstein bound~\eqref{EqnTruncate3}:}
 by the variational definition of matrix  operator norm and Cauchy-Schwartz inequality, we have
\begin{align*}
\opnorm{\E[\newtruncate_A(\Sam_i,\Par)^{\mytrans} \newtruncate_A(\Sam_i,\Par)}
&\leq
\sup_{\Direcb,\Direcc\in\Sphere{\ParDim-1}}
\E[\newtruncate_A(\Sam_i,\Par)^{\mytrans} \newtruncate_A(\Sam_i,\Par)][\Direcb,\Direcc]\\
&=
\sup_{\Direcb,\Direcc\in\Sphere{\ParDim-1}}\sum_{j=1}^\ParDim\E[\newtruncate_A(\Sam_i,\Par)[e_j,\Direcb]\cdot \newtruncate_A(\Sam_i,\Par)][e_j,\Direcc]]\\
&\leq 
2 \sum_{j=1}^d \sup_{\Direcb\in\Sphere{\ParDim-1}}\E[\newtruncate_A(\Sam_i,\Par)[e_j,\Direcb]^2]\leq 2 \sum_{j=1}^d \sup_{\Direcb\in\Sphere{\ParDim-1}}\E[\nabla^2\ZFun(\Sam_i,\Par)[\Direc,e_j,\Direcb]^2] .
\end{align*}
Since Assumption~\ref{ass:tail} implies $\E[\nabla^2\ZFun(\Sam_i,\Par)[\Direc,e_j,\Direcb]^2]\leq c\OrPar^2$ for any ${\Direc,\Direcb\in\Sphere{\ParDim-1}}$, it follows that $\opnorm{\E[\newtruncate_A(\Sam_i,\Par)^{\mytrans} \newtruncate_A(\Sam_i,\Par)}\leq c\OrPar^2\cdot \ParDim$. Similarly, we obtain
\begin{subequations}
\begin{align}
    \max\{\opnorm{\E[\newtruncate_A(\Sam_i,\Par)^{\mytrans}
        \newtruncate_A(\Sam_i,\Par)},\opnorm{\E[\newtruncate_A(\Sam_i,\Par)
          \newtruncate_A(\Sam_i,\Par)^{\mytrans}} \} \leq
        c\OrPar^2\cdot d.\label{eq:matrix_bern_1}
\end{align} Moreover, by the definition of $\newtruncate_A$, we have \begin{align}
\opnorm{\newtruncate_A(\Sam_i,\Par)}\leq A. \label{eq:matrix_bern_12}
\end{align}
\end{subequations}
Therefore, applying the matrix Bernstein inequality (see e.g., exercise 6.10 in Wainwright~\cite{wainwright2019high})  to the matrices $\newtruncate_A(\Sam_i,\Par),i\in[\numobs]$ yields
 equation~\eqref{EqnTruncate3}.




\subsection{Proof of~\Cref{lm:sup_op_norm_zfun}}
\label{sec:app-proof-sup_op_norm_zfun}

 We prove each of the three sub-claims in turn, beginning
with the bound~\eqref{eq:bounded_hessian_zfun_lm}.


\paragraph{Proof of the bound~\eqref{eq:bounded_hessian_zfun_lm}.}

We prove the discretization upper bound
\begin{subequations}
  \begin{align}
\label{EqnDbound}    
\sup_{\Par \in \Parspace} \opnorm{\frac{1}{\Numobs} \sum_{i=1}^\Numobs
  \nabla^2_\Par \ZFun(\Sam_i, \Par)[u,v,w]} & \leq 2 \max_{u,v,w \in
  \Coverset} Y_\star(u, v,w),
  \end{align}
\mbox{where $Y_\star(u, v, w) \defn \sup \limits_{\Par \in \ParSpace}
  \frac{1}{\Numobs} \sum_{i=1}^\Numobs \nabla^2_\Par \ZFun(\Sam_i,
  \Par)[u,v,w]$,} combined with the tail bound
\begin{align}
\label{EqnTbound}  
  \Prob \left [Y_\star(u,v,w) \geq c' \sigma + t \right] & \leq \exp
  \big(- \frac{c\Numobs t^2}{ \OrParc^2} \big) + 3 \exp \Big( -
  \Big\{\frac{c \Numobs t}{\OrParc (1+\log^\InvOrConc \Numobs)} \Big\}^{\OrConc} \Big),
\end{align}
\end{subequations}
where $c, c'$ are constants.  The claimed
bound~\eqref{eq:bounded_hessian_zfun_lm} follows by combining these
results with a union bound over the $N \leq 17^\usedim$ elements of
the cover. \\

\myunder{Proof of discretization bound~\eqref{EqnDbound}:} By
definition of the $1/8$-covering set $\Coverset$, for any triple
$(\Direca, \Direcb, \Direcc)$ of vectors in the unit sphere, we can
find a triple $(u^a, u^b, u^c)$ of elements in the cover such that
$\max \big \{ \|\Direca - u^a\|_2, \; \|\Direcb - u^b\|_2, \;
\|\Direcc - u^c\|_2 \big \} \leq \frac{1}{8}$.  Our goal is to bound
the supremum $\sup_{\Par \in \ParSpace} \Term_\theta$, where
$\Term_\theta \defn \opnorm{ \frac{1}{\Numobs} \sum_{i=1}^\Numobs
  \nabla^2_\Par \ZFun(\Sam_i, \Par)}$.  By the variational definition
of the operator norm on tensors, we have
\begin{align*}
\Term_\theta & = \sup_{(u,v,w)} Y_\theta(u, v, w) \qquad \mbox{where
  $Y_\theta(u,v,w) \defn \frac{1}{\Numobs} \sum_{i=1}^\Numobs
  \nabla^2_\Par \ZFun(\Sam_i, \Par)[u,v,w]$.}
\end{align*}
The same discretization argument used previously can be used to show
that
\begin{align*}
  \Term_\theta = \sup_{(u,v,w)} Y_\theta(u,v,w) & \leq \max_{u^a, u^b,
    u^c \in \Coverset} Y_\theta(u^a, u^b, u^c) + \tfrac{3}{8}
  \Term_\theta.
\end{align*}
Note that $Y_\star(u,v,w) = \sup_{\Par \in \Parspace} Y_\theta(u,
v,w)$.  Thus, re-arranging and then taking suprema over $\theta$ on
both sides, we arrive at the claimed bound~\eqref{EqnDbound}. \\

\myunder{Proof of tail bound~\eqref{EqnTbound}:} For the fixed triple
$(u, v, w)$ and each parameter $\Par$, define the zero-mean function
$f_\theta(Z) = \nabla^2_\Par \ZFun(\Sam, \Par)[u,v,w] -
\Exs\Big[\nabla^2_\Par \ZFun(\Sam, \Par)[u,v,w] \Big]$.  We can then
write
\begin{align}
Y_\star(u, v, w) & \leq \sup_{\Par \in \ParSpace} \frac{1}{\Numobs}
\sum_{i=1}^\numobs f_\Par(\Sam_i) + \sup_{\Par \in \ParSpace} \Exs
\big[ \nabla^2_\Par \ZFun(\Sam, \Par)[u,v,w] \big] \notag \\
\label{EqnSstar}
& \leq S_\star + c \sigma, \qquad \mbox{where $S_\star \defn
  \sup_{\Par \in \ParSpace} \frac{1}{\Numobs} \sum_{i=1}^\numobs
  f_\Par(\Sam_i)$,}
\end{align}
and the second inequality follows from
condition~\ref{ass:tail}. 

Thus, we have reduced our problem to obtaining a tail bound on the
random variable $S_\star$, for which we can make use of existing
results on the suprema of empirical processes.  Introducing the
shorthand $\Funclass = \{f_\Par, \: \Par \in \ParSpace \}$, observe
that we have
\begin{align*}
  \vecnorm{\sup_{\Par \in \ParSpace} f_\Par}{\psi_\OrConc} & =
  \vecnorm{\sup_{\Par \in \ParSpace}\Big|\nabla^2_\Par \ZFun(\Sam, \Par)[u,v,w] -
    \Exs\Big[\nabla^2_\Par \ZFun(\Sam, \Par)[u,v,w] \Big]\Big|}{\psi_\OrConc}
  \leq c \OrParc,
\end{align*}
using condition~\ref{ass:tail}.  Similarly, we have
$\E[\sup_{\Par\in \ParSpace} f^2_\theta(\Sam)] \leq c \OrParc^2$.
Consequently, we can apply~\Cref{LemAdam} to assert that
\begin{align*}
  \Prob \left [S_\star \geq 2 \Exs[S_\star] + t \right] & \leq \exp
  \big(- \frac{c\Numobs t^2}{ \OrParc^2} \big) + 3 \exp \Big( -
  \Big\{\frac{c \Numobs t}{\OrParc (1+\log^{\InvOrConc} \Numobs)} \Big\}^{\OrConc} \Big).
\end{align*}
Combined with our earlier inequality~\eqref{EqnSstar}, the tail
bound~\eqref{EqnTbound} follows.


\paragraph{Proof of equation~\eqref{eq:bounded_third_zfun_lm}.}

We use a similar approach based on discretization and tail bounds.  In
particular, we first show that
\begin{subequations}
\begin{align}
 \label{EqnDthree}
\sup_{\Par \in \ParSpace} \opnorm{ \frac{1}{\Numobs} \sum_{i=1}^\Numobs
  \nabla^3_\Par \ZFun(\Sam_i, \Par)[\Direc]} \leq c \sqrt{ \log
  (\Numobs/\delta)} \max_{\Direcb, \Direcc, \Direcd \in \Coverset}
Y(\Direca, \Direcb, \Direcc, \Direcd)
\end{align}
\mbox{with probability at least $1 - \tfrac{\delta}{2}$,} where
\mbox{$Y(\Direc, \Direcb, \Direcc, \Direcd) \defn \frac{1}{\Numobs}
  \sum_{i=1}^{\Numobs} \Factorb(\Sam_i, \Direc, \Direcb, \Direcc,
  \Direcd)$.}  We then show that
\begin{align}
\label{EqnTthree}
\Prob \Big[ Y(\Direc, \Direcb, \Direcc, \Direcd) \geq c \OrPard + t
  \Big] & \leq \exp \big(- \frac{c \Numobs t^2 }{\OrPard^2}) + 3 \exp
\Big( - \big(\frac{c\Numobs t}{ \OrPard(1+ \log^{\InvOrCond} \numobs)} \big)^{\OrCond} \Big).
\end{align}
\end{subequations}
Combining these two inequalities with a union bound over the cover
$\Coverset$ yields the claim. \\

\myunder{Proof of the bound~\eqref{EqnDthree}:}
We have
\begin{align}
\sup_{\Par \in \ParSpace} \opnorm{ \frac{1}{\Numobs} \sum_{i=1}^\Numobs
  \nabla^3_\Par \ZFun(\Sam_i, \Par)[\Direc]} & \stackrel{(i)}{\leq} 2
\max_{\Direcb, \Direcc, \Direcd \in \Coverset} \frac{1}{\Numobs}
\sum_{i=1}^\Numobs \sup_{\Par \in \ParSpace} \abss{ \nabla^3_\Par
  \ZFun(\Sam_i, \Par)[\Direc, \Direcb, \Direcc, \Direcd]} \notag \\
\label{eq:discre_third_eq_2_lm}
& \stackrel{(ii)}{\leq} 2\max_{\Direcb, \Direcc, \Direcd\in \Coverset}
\frac{1}{\Numobs} \sum_{i=1}^\Numobs\Factora(\Sam_i, \Direc)
\Factorb(\Sam_i, \Direc, \Direcb, \Direcc, \Direcd),
\end{align}
where step (i) follows from the same discretization argument; and
step (ii) follows from condition~\ref{ass:tail}.

Applying a union bound over the $\numobs$ sub-Gaussian variables
$\{\Factora(\Sam_i, \Direc) \}_{i=1}^\numobs$ yields
\begin{align*}
\Prob \Big[\max_{i=1, \ldots, \numobs} \Factora(\Sam_i, \Direc) \geq c
  \sqrt{\log (\Numobs/\delta)} \Big] & \leq \delta/2
\end{align*}
for some constant $c>0$. Combining this bound with
equation~\eqref{eq:discre_third_eq_2_lm} yields
the claim~\eqref{EqnDthree}. \\

\myunder{Proof of the bound~\eqref{EqnTthree}:} We prove the tail
bound by applying~\Cref{LemAdam} to the single-element
function class $\Funclass = \{g \}$ where $g \defn \Factorb(\Sam,
\Direc, \Direcb, \Direcc, \Direcd) - \E [\Factorb(\Sam, \Direc,
  \Direcb, \Direcc, \Direcd)]$.  By \mbox{condition~\ref{ass:tail},}
we have
\begin{align*}
\rnorm{g}{\psi_\OrCond} \leq c \OrPard, \quad \var(g) \leq c \OrPard^2, \quad
\mbox{and} \quad \Exs[\Factorb(\Sam, \Direc, \Direcb, \Direcc,
  \Direcd)] \leq c \OrPard,
\end{align*}
so that an application of~\Cref{LemAdam} yields the claimed
bound~\eqref{EqnTthree}.


\paragraph{Proof of equation~\eqref{eq:bounded_grad_zfun_lm}.}

Let $\NewCover(\epsilon)$ be an $\epsilon$-cover of the parameter
space $\ParSpace$; by the $\rho$-boundedness condition, we can find
such a cover with log cardinality $\log \NewCover(\epsilon) \leq
\usedim \, \log \big(1 + \frac{2 \rho}{\epsilon} \big)$.  Our proof is
based on the following two auxiliary claims: first, with the choice
$\eps = (\dimdelta (1 + \log\numobs))^{-2}$, we have
\begin{subequations}
  \begin{align}
\label{EqnDiscretizeEps}    
\sup_{\Par \in \ParSpace} \opnorm{ \frac{1}{\Numobs}
  \sum_{i=1}^\Numobs \nabla_\Par \ZFun(\Sam_i, \Par)} & \leq 2 \left
\{ \max_{\Direc, \Direcb \in \Coverset} \max_{\Par \in
  \NewCover(\epsilon)} Y(\Direc, \Direcb, \Par) + c \OrParc \right \}
\end{align}
with probability at least $1-\delta/2$, 
where $Y(\Direca, \Direcb, \Par) \defn \frac{1}{\Numobs}
\sum_{i=1}^\Numobs \nabla_\Par \ZFun(\Sam_i, \Par)[\Direca, \Direcb]
$.  Second, for any fixed triple $(\Direca, \Direcb, \Par)$, we have
\begin{align}
\label{EqnTailEps}      
\Prob \Big[ Y(\Direc, \Direcb, \Par) \geq c \OrParc + t \Big] & \leq
\exp \big(- \frac{c\Numobs t^2 }{\OrParb^2} \big) + 3 \exp \Big( -
\frac{c \Numobs t}{ \sigma ( 1 + \log \Numobs)} \Big).
\end{align}
\end{subequations}
The claim~\eqref{eq:bounded_grad_zfun_lm} follows by combining these
results with a union bound over the all elements in the Cartesian
product $\Coverset$ with $\NewCover(\epsilon)$. \\

\myunder{Proof of the bound~\eqref{EqnDiscretizeEps}:} The same
discretization argument, over pairs $\Direca, \Direcb$ in the
$1/8$-cover $\Coverset$ yields
\begin{subequations}
\begin{align}
\label{EqnOlympicsOne}  
  \opnorm{ \frac{1}{\Numobs} \sum_{i=1}^\Numobs \nabla_\Par
    \ZFun(\Sam_i, \Par)} & \leq 2 \max_{\Direca, \Direcb \in
    \Coverset} Y(\Direca, \Direcb, \Par).
\end{align}
Next we write $\Par = \Par' + \Delta$ for some element $\Par' \in
\NewCover(\epsilon)$ with $\|\Delta\|_2 \leq \epsilon$.  We then have
\begin{align}
  Y(\Direca, \Direcb, \Par) & \leq Y(\Direca, \Direcb, \Par') +
  \sup_{\|\theta_1 - \theta_2\|_2 \leq \epsilon} \Big|
  \frac{1}{\numobs} \sum_{i=1}^\numobs \big[\nabla \Zfun(\Sam_i,
    \theta_1) - \nabla\Zfun(\Sam_i, \theta_2) \big][u,v] \Big| \notag \\
  & \stackrel{(i)}{\leq} Y(\Direca, \Direcb, \Par') + \epsilon \;
  \sup_{\Par \in \ParSpace} \opnorm{ \frac{1}{\numobs}
    \sum_{i=1}^\numobs \nabla^2 \Zfun(\Sam_i, \Par)} \notag \\
\label{EqnOlympicsTwo}  
  & \stackrel{(ii)}{\leq} Y(\Direca, \Direcb, \Par') + c \sigma
\end{align}
\end{subequations}
with probability at least $1-\delta/2$, 
where step (i) follows by a Taylor series expansion; and step (ii)
follows from our previous bound~\eqref{eq:bounded_hessian_zfun_lm}
combined with our choice of $\epsilon$.  Taking suprema over $\theta$
in the bound~\eqref{EqnOlympicsOne} and combining with
inequality~\eqref{EqnOlympicsTwo}, we find that
\begin{align*}
\sup_{\Par \in \ParSpace} \opnorm{ \frac{1}{\Numobs}
  \sum_{i=1}^\Numobs \nabla_\Par \ZFun(\Sam_i, \Par)} & \leq 2 \left
\{ \max_{\Direca, \Direcb \in \Coverset} \max_{\Par' \in
  \NewCover(\epsilon)} Y(\Direca, \Direcb, \Par) + c \sigma \right \},
\end{align*}
which establishes the 
claim~\eqref{EqnDiscretizeEps}. \\

\myunder{Proof of the tail bound~\eqref{EqnTailEps}:} Let
$\Ytil(\Direca, \Direc, \Par) = Y(\Direca, \Direcb, \Par) -
\Exs[Y(\Direca, \Direcb, \Par)]$ be the recentered random variable.
Since $\Exs[Y(\Direca, \Direcb, \Par)] \leq c \sigma$ by
condition~\ref{ass:tail}, it suffices to prove a tail bound of the
form~\eqref{EqnTailEps} on the recentered version.

Using $\Ytil$ as a shorthand, we do so by
applying~\Cref{LemAdam} to the single-element function
class $\Funclass = \{\Ytil\}$.  By \mbox{condition~\ref{ass:tail},} we
have
\begin{align*}
\rnorm{\Ytil}{\psi_1} \leq c \OrPard, \quad \mbox{and} \quad \var(\Ytil)
\leq c \sigma^2,
\end{align*}
so that the claimed tail bound~\eqref{EqnTailEps} follows
from~\Cref{LemAdam}.



\end{document}